\documentclass[12pt]{amsart}
\usepackage{graphicx} 
\usepackage{tikz}
\usepackage[english]{babel}
\usepackage[letterpaper,top=2cm,bottom=2cm,left=3cm,right=3cm,marginparwidth=1.75cm]{geometry}
\usepackage{amsmath, amssymb, latexsym}
\usepackage{mathtools}
\usepackage{graphicx}
\usepackage[colorlinks=true, allcolors=blue]{hyperref}
\usepackage{amsfonts}
\usepackage{amsthm}
\usepackage{color}
\usepackage{indentfirst}
\usepackage{comment}
\usepackage{enumerate}
\usepackage{enumerate,xcolor}
\usepackage{amsrefs}

\numberwithin{equation}{section}

\theoremstyle{plain}
\newtheorem{theorem}{Theorem}[section]
\newtheorem{proposition}[theorem]{Proposition}
\newtheorem{lemma}[theorem]{Lemma}
\newtheorem{corollary}[theorem]{Corollary}

\theoremstyle{definition}
\newtheorem{definition}[theorem]{Definition}
\newtheorem{example}[theorem]{Example}
\newtheorem{Conjecture}[theorem]{Conjecture}
\newtheorem{remark}[theorem]{Remark}

\DeclareMathOperator{\Aut}{Aut}

\DeclareMathOperator{\ad}{ad}
\DeclareMathOperator{\Ad}{Ad}
\DeclareMathOperator{\Inn}{Inn}
\DeclareMathOperator{\id}{id}

\DeclareMathOperator{\Span}{span}
\newcommand{\B}{\mathcal{B}}

\newcommand{\N}{{\mathbb N}}
\newcommand\GK{\operatorname{GKdim}}
\newcommand\Mr{\operatorname{M}}
\newcommand\cB{\mathcal{B}}
\newcommand{\J}{\mathbf J}
\newcommand{\Id}{\operatorname{Id}}
\newcommand{\fe}{\mathfrak{e}}

\newcommand{\kk}{\Bbbk}
\DeclareMathOperator{\Hom}{Hom}

\title[The automorphism group of the super Jordan plane]{The automorphism group of the super Jordan plane}

\author[W. de Araujo]{Willian de Araujo}
\address{\textnormal{Willian de Araujo}\newline
	Universidade Tecnol\'ogica Federal do Paran\'a,
	Toledo, PR, Brazil}
\email{wilianmat@yahoo.com.br}

\author[D. Bagio]{Dirceu Bagio}
\address{\textnormal{Dirceu Bagio}\newline
	Universidade Federal de Santa Catarina,
	Florian\'opolis, SC, Brazil}
\email{d.bagio@ufsc.br}

\author[R. Baltazar]{Rene Baltazar}
\address{\textnormal{Rene Baltazar}\newline
	Universidade Federal do Rio Grande,
	Santo Ant\^onio da Patrulha, RS, Brazil}
\email{renebaltazar.furg@gmail.com}

\author[A. Santana]{Adriano Santana}
\address{\textnormal{Adriano Santana}\newline
	Universidade Tecnol\'ogica Federal do Paran\'a,
	Toledo, PR, Brazil}
\email{adrianosantana@utfpr.edu.br}

\author[R. Vinciguerra]{Robson Vinciguerra}
\address{\textnormal{Robson Vinciguerra}\newline
	Universidade Tecnol\'ogica Federal do Paran\'a,
	Toledo, PR, Brazil}
\email{robsonwv@gmail.com}


\thanks{D. Bagio was partially supported by Funda\c{c}\~ao
	de Amparo \`a Pesquisa e Inova\c{c}\~ao do Estado de
	Santa Catarina (FAPESC), Call No.~21/2024.}

\thanks{R. Baltazar was partially supported by
Funda\c{c}\~ao de Amparo \`a Pesquisa do Estado
do Rio Grande do Sul (FAPERGS), Project
No~24/2551-0001560-4.}

\begin{document}
	
	\begin{abstract}
		We complete the description of the automorphism group of the super Jordan plane $\mathcal B$ over a field of characteristic zero, extending the classification of
		super-automorphisms obtained by N. Andruskiewitsch and F. Dumas in \cite{AD}. We characterize a subgroup $\mathbb L$ stabilizing two one-dimensional vector spaces and obtain two semidirect product decompositions of it. This subgroup contains all super-automorphisms and its structure motivates our main result. Precisely, we prove that $\operatorname{Aut}\mathcal B
		=\operatorname{Inn}\mathcal B\rtimes
		\operatorname{Aut}_{\mathrm{super}}\mathcal B$,
		where $\operatorname{Inn}\mathcal B$ and
		$\operatorname{Aut}_{\mathrm{super}}\mathcal B$ denote,
		respectively, the groups of inner automorphisms and
		super-automorphisms of $\mathcal B$. We also characterize the nilpotent elements of $\mathcal B$ and prove that every nonzero nilpotent has nilpotency index two. Furthermore, we give a geometric interpretation of three distinguished families of nilpotent elements. After degree truncation and projectivization, two of these families yield projective linear subspaces, while the third yields quadratic
		Veronese varieties of arbitrary dimension.
	\end{abstract}
	
	\maketitle


	\section{Introduction}
	
	The study of automorphism groups of algebras presents significant challenges, even for commutative polynomial rings. A notable example is Nagata's conjecture,	formulated in 1972, concerning the wildness of an automorphism of $\mathbb C[x_1,x_2,x_3]$. Shestakov and Umirbaev  proved the conjecture in 2004 \cite{ShestakovUmirbaev}. Automorphisms and derivations of
	noncommutative algebras, particularly quantum
	algebras and related structures, have also attracted
	considerable attention; see, for example,
	\cites{Alev,Shirikov,Samuel,vivas}. In characteristic zero,
	exponentials of locally nilpotent derivations provide a further source of automorphisms. A recent approach, developed by Ceken, Palmieri, Wang, and Zhang~\cites{CPWZ2015,CPWZ2016}, uses discriminants to determine automorphism groups of several classes of quantized algebras, including tensor products of quantum Weyl algebras and some skew polynomial rings.
	
	Throughout, $\Bbbk$ denotes a field of characteristic zero. The Jordan plane is a basic example for which the automorphism group admits an explicit description. In the notation used below, it is the algebra
	\[
	\mathbf J
	=\Bbbk\langle y_1,y_2\rangle/(y_1y_2-y_2y_1-y_2^2).
	\]
	It is a quadratic, graded domain of Gelfand--Kirillov dimension two, and its representation theory was studied in \cite{Na}. Shirikov~\cite{Shirikov} showed that every automorphism of $\mathbf J$ is determined by
	\[
	y_1\longmapsto\alpha y_1+p(y_2),\quad
	y_2\longmapsto\alpha y_2,
	\quad \alpha\in\Bbbk^\times,\quad p\in\Bbbk[t].
	\]
	Thus every automorphism preserves the one-dimensional vector space $\Span\{y_2\}$, and the image of $y_1$ is determined by the same scalar $\alpha$ and a polynomial in $y_2$. Since $\mathbf J$ is a connected graded domain, its only invertible elements are the nonzero scalars. In particular, its inner automorphism group is trivial.
	
	The Jordan plane also occurs in the study of Nichols algebras and of pointed Hopf algebras of finite Gelfand--Kirillov dimension \cites{aah-triang,aah-jordan}. In this setting, Andruskiewitsch, Angiono and Heckenberger introduced the \emph{super Jordan plane} \cite{aah-triang}. It is the associative algebra
	\[
	\B=\Bbbk\langle x,y\rangle/(x^2,\;yw-wy-xw),
	\quad w:=xy+yx.
	\]
	Equivalently, the second defining relation is $y^2x-xy^2-xyx=0$. As is the case for the Jordan plane, $\B$ is a Nichols algebra associated with a two-dimensional braided vector space of block type and has Gelfand--Kirillov dimension two. Writing $z:=yx$,
	a PBW basis of $\B$ is given by
	\[
	\bigl\{x^\epsilon z^a y^b:
	\epsilon\in\{0,1\},\ a,b\in\mathbb N_0\bigr\}.
	\]
	
	The Jordan plane and the super Jordan plane are also
	closely related as associative algebras: the map $y_1\mapsto y^2$, $y_2\mapsto w$ identifies the Jordan plane with the subalgebra of $\B$ generated by $y^2$ and $w$~\cite{AD}.
	There is, however, a difference that matters for automorphisms:
	the relation $x^2=0$ gives $\B$ both nonzero nilpotent elements
	and nonscalar invertible elements.
	
	In this paper, we distinguish two subalgebras of $\B$
	isomorphic to the Jordan plane. Under the preceding embedding,
	we identify $\mathbf J$ with $\Bbbk\langle y^2,w\rangle$.
	We denote the other subalgebra by
	\[
	\mathcal J=\Bbbk\langle u,w\rangle,
	\qquad u=y^2-xy.
	\]
	By Lemma~\ref{lemma-centralizador-x1},
	$\mathcal J$ is contained in the centralizer $C_{\B}(x)$ of $x$ in $\B$.
	The subalgebra $\mathbf J$ is not contained in $C_{\B}(x)$,
	since $[y^2,x]=xz\neq0$. In Section~4 we construct an embedding of $\mathcal B$ into the matrix algebra $\Mr_2(\mathbf J)$; we use the generators
	$s$ and $w$ of $\mathbf J$, with $s$ corresponding
	to $y^2$.
	
	Several aspects of the structure of $\B$ have already been
	investigated. Over an algebraically closed field of characteristic zero, Andruskiewitsch, Bagio, Della Flora and Fl\^ores proved that every finite-dimensional simple $\B$-module is one-dimensional and classified the indecomposable modules of dimensions two and three~\cite{ABDF}. Under the same hypotheses, Reca and Solotar computed the Hochschild homology and cohomology of $\B$ and related the Lie algebra structure on its first Hochschild cohomology to the Virasoro algebra~\cite{RS}. More recently, Lu studied the prime spectrum and the representation theory of $\B$. In particular, he realized the localization of $\B$ at the powers of $w$ as a $2\times2$ matrix algebra over a localization of the first Weyl algebra~\cite{Lu2026SJP}.
	
	For the automorphism problem, our starting point is the work of Andruskiewitsch and Dumas~\cite{AD}. 
	Since the defining relations are homogeneous with
	$\deg x=\deg y=1$, the algebra $\B$ is naturally
	$\mathbb N_0$-graded. Consequently, $\B$ is
	$\mathbb Z_2$-graded with
	
	\[
	\B=\B_{\bar0}\oplus\B_{\bar1},\quad
	\B_{\bar0}=\bigoplus_{j\geq0}\B_{2j},\quad
	\B_{\bar1}=\bigoplus_{j\geq0}\B_{2j+1}.
	\]
	A \emph{super-automorphism} is an algebra automorphism that
	preserves $\B_{\bar0}$ and $\B_{\bar1}$. Thus it preserves
	parity, but may send a homogeneous element to a sum of terms
	of different degrees having the same parity. By \cite[Lemma~5.2 and Theorem~5.3]{AD}, the super-automorphisms of $\B$ are precisely the maps \begin{equation}\label{eq:intro-super-aut}
		\sigma_{p,\mu}(x)=\mu x,\quad
		\sigma_{p,\mu}(y)=\mu y+p(w)x,
		\quad p\in\Bbbk[t],\quad \mu\in\Bbbk^\times.
	\end{equation}
	This determines the automorphisms of $\B$ as a superalgebra.
	Our main result describes the full group of algebra automorphisms
	over an arbitrary field of characteristic zero in terms of inner
	automorphisms and the subgroup in \eqref{eq:intro-super-aut}. Precisely, we prove the following.
	
	\medskip
	\noindent\textbf{Theorem A.}
	\emph{Let $\B$ be the super Jordan plane over a field $\Bbbk$
		of characteristic zero. Then, $\Aut\B=\Inn\B\rtimes\Aut_{\mathrm{super}}\B$}.

	In particular, the restriction of the quotient map
	$\Aut\B\to\operatorname{Out}\B$ gives an isomorphism
	$\Aut_{\mathrm{super}}\B\simeq\operatorname{Out}\B$. Thus, every automorphism of $\B$ can be written as $\Ad_U\circ \,\sigma_{p,\mu}$, where $U$ is an invertible element of $\B$ and $\Ad_U(a)=UaU^{-1}$; the inner automorphism $\Ad_U$ and the pair $(p,\mu)$ are uniquely determined. The full group is strictly larger than the group of super-automorphisms. For instance, conjugation by $1+x$ sends $y$ to $y+xy-yx-xyx$, whose even component $xy-yx$ is nonzero. The proof uses the two-sided ideal $I=(x)$ of $\B$. We first determine the induced actions on $\B/I\simeq\Bbbk[y]$ and $I/I^2$. We then embed $\B$ in the $2\times 2$ matrix ring $\Mr_2(\mathbf J)$ over the Jordan plane and show that every automorphism of $\B$ extends uniquely to this ring. We begin by establishing the result over an
	algebraically closed field of characteristic zero. In Lemma~\ref{lem:picard-localization} of
	Section~4.4, we apply a result of Cannings and
	Holland~\cite[Corollary~1.13]{CH} to the ring of differential operators on $\Bbbk[w,w^{-1}]$. Since $\Bbbk[w,w^{-1}]$ is a principal ideal domain, we deduce that every invertible $\Bbbk$-central bimodule over this ring of operators is free of rank one as a right module. Proposition~\ref{prop:standard} uses this freeness to factor the extensions of automorphisms of $\B$,
	leading to the decomposition in Theorem~\ref{thm:main} over an algebraically closed field. The application of Cannings and Holland is the only step in the proof of Theorem~A that requires algebraic closedness. To obtain Theorem~A over an arbitrary field of characteristic zero, we extend scalars to the algebraic closure and descend both factors of the resulting decomposition, as shown in Corollary~\ref{cor:descent}. The results on the subgroup $\mathbb L$ and the algebraic descriptions of nilpotent elements also hold over any field of characteristic zero. For the geometric discussion in Section~5, we assume that the ground field is algebraically closed.

 	Before proving Theorem~A, we determine a subgroup of $\Aut\B$ for which both factors in the decomposition can be described explicitly. Let $\mathbb L$ consist of the automorphisms that stabilize both $\Span\{x\}$ and $\Span\{w\}$. Its elements are precisely the maps
	\begin{equation}\label{eq:intro-L}
		\begin{aligned}
			\sigma_{p,q,\mu}(x)&=\mu x,\\
			\sigma_{p,q,\mu}(y)&=\mu y+xp(z)+xq(z)y-q(z)z,
		\end{aligned}
		\quad\quad p,q\in\Bbbk[t],\quad \mu\in\Bbbk^\times.
	\end{equation}
	The parameters are unique and $q=0$ gives precisely the
	super-automorphisms in \eqref{eq:intro-super-aut}, since
	$xp(z)=p(w)x$. We determine the group law for $\mathbb L$ and show that
	$\mathbb L\simeq\mathbb T\rtimes\Bbbk^\times$. Here, $\mathbb T$ is the abelian group $\Bbbk[t]\times\Bbbk[t]$ with operation $(p,q)(p',q')=(p+p'-2tqq',\,q+q')$, $p,q,p',q'\in \Bbbk[t]$. We also identify the subgroup of inner
	automorphisms in $\mathbb L$:
	\[
	N:=\Inn\B\cap\mathbb L
	=\{\sigma_{-tq^2,q,1}:q\in\Bbbk[t]\}.
	\]
	It follows from Corollary \ref{cor-quotient-L-inner} that $\mathbb L=N\rtimes\Aut_{\mathrm{super}}\B$. Moreover, each $\sigma_{p,q,1}$ is the exponential of an explicit locally nilpotent derivation, extending the construction for super-automorphisms in \cite{AD}.
	
	We also investigate the nilpotent elements of $\B$
	and the geometry of their families. Using Lu's matrix description, we prove that every nilpotent element of $\B$ has square zero and give a parametrization of all such elements without any assumption on their parity. For $m\geq0$, we consider the projective square-zero locus
	\[
	\overline{\mathcal Z}_m
	=\{[a]\in\mathbb P(\B_{\leq 2m+1}):a^2=0\},
	\quad
	\B_{\leq d}=\bigoplus_{j=0}^{d}\B_j.
	\]
	This locus is defined by homogeneous quadratic equations. For $m\geq1$, the projective classes of odd
	square-zero elements form the disjoint union
	of two projective linear subspaces,
	$\mathbb P(S_m)$ and $\mathbb P(H_m)$.
	We also construct a third family of square-zero
	elements using polynomials in
	$V_m=\Bbbk[y]_{\leq m}$. We prove that
	\[
	\mathbb P(V_m)\longrightarrow\overline{\mathcal Z}_m,\quad
	[p]\longmapsto[Q_m(p)],\quad Q_m(p)=p(-y)xp(y),
	\]
	is a closed embedding whose image is projectively equivalent to the quadratic Veronese variety of $\mathbb P^m$. For $m\geq1$, its points in $\mathbb P(S_m)$ correspond to even polynomials $p$ and those in $\mathbb P(H_m)$ to odd polynomials. For $m=1$, the entire locus is the union of the plane $\mathbb P(S_1)$ and the Veronese conic, which meet only at $[x]$. The degree-five element $(1+y)xz(1-y)$ has square zero and belongs to neither linear family; its projective class also lies outside the Veronese family.
	
	The paper is organized as follows. In Section~2, we
	recall the basic properties of the Jordan and super
	Jordan planes and establish some commutation formulas
	used throughout the paper. In Section~3, we determine
	the structure of $\mathbb L$ and describe its relation
	to inner automorphisms, super-automorphisms and locally
	nilpotent derivations. Section~4 is devoted to the proof
	of Theorem~A, first over algebraically closed fields
	and then over arbitrary fields of characteristic zero.
	Finally, in Section~5, we obtain an explicit parametrization of all nilpotent
	elements of $\B$ and study the
	geometry of the square-zero loci, with particular
	emphasis on their maximal linear subspaces and
	Veronese varieties.

	\subsection*{Notations and Conventions}
	All algebras are associative and unital over a field $\Bbbk$ of characteristic zero and all
	algebra automorphisms are $\Bbbk$-linear. We write
	$\N=\{1,2,\ldots\}$ and $\N_0=\N\cup\{0\}$. For an $\mathbb N_0$-graded algebra
	$R=\bigoplus_{n\geq0}R_n$, we write
	$R_{\leq n}=\bigoplus_{j=0}^n R_j$.
	For a nonzero element $a=\sum_{j\geq0}a_j$,
	with $a_j\in R_j$, we set
	$\deg a=\max\{j:a_j\neq0\}$. We set $\deg 0=-\infty$.
		We write $\Span(X)$ for the
	$\Bbbk$-linear span of a set $X$. For a list of vectors $v_1,\ldots,v_r$, we also write
	$\Span\{v_1,\ldots,v_r\}$. The group of invertible elements (or units) of an algebra
	$R$ is denoted by $R^\times$. For $n\in\N_0$,
	we use the notation $(n)_0=1$ and
	$(n)_i=n(n+1)\cdots(n+i-1)$ for $i\geq1$. We denote
	by $\Ad_u$, $u\in R^\times$, the inner automorphism of $R$ given by $\Ad_u(a)=uau^{-1}$, $a\in R$. We write $\Inn R$ for the group of inner automorphisms of $R$ and $\Aut R$ for the group of automorphisms of $R$. We denote the
	algebra of $n\times n$ matrices over $R$ by $\Mr_n(R)$ and its identity matrix by $\Id_n$. We denote by $e_{ij}$ the elementary matrix with $1$ in
	the $(i,j)$-entry and $0$ elsewhere. We write $\GK R$ for the Gelfand--Kirillov dimension of  $R$. We write $[a,b]=ab-ba$ for the commutator and
	$\ad(a)$ for the inner derivation given by
	$\ad(a)(b)=[a,b]$, for $a,b\in R$. For a subset $X$ of an algebra $R$, we write
	$C_R(X)=\{a\in R:[a,x]=0\text{ for all }x\in X\}$
	for its centralizer and $Z(R)$ for the center of $R$.
	If $X=\{x\}$, we write $C_R(x)$ instead of $C_R(X)$. For a subset $X$ of an algebra $R$, we denote by
	$(X)$ the two-sided ideal of $R$ generated by $X$. For a regular normal element $s$ of an algebra $R$,
	we denote the localization of $R$ at the powers of
	$s$ by $R[s^{-1}]$. When it exists, the classical
	ring of quotients of $R$ is denoted by $Q(R)$.
	If $R$ is an Ore domain, we also write
	$\operatorname{Frac}(R)$ for its quotient division ring.

	\section{Jordan and super Jordan planes}
	
	\subsection{Jordan plane}
	We begin by recalling the definition of the Jordan plane
	and some relevant results. For further details, we refer
	the reader to \cite{Na}.

	\begin{definition}\label{def:Jordan-plane}
		The Jordan plane is the associative $\Bbbk$-algebra $\mathbf J$
		generated by $y_1$ and $y_2$ subject to the single quadratic relation
		\[
		y_1y_2-y_2y_1-y_2^2=0.
		\]
	\end{definition}
	
It is well known that $\mathbf J$ is an
$\N_0$-graded domain with $\deg y_1=\deg y_2=1$.
The set $\{y_1^ay_2^b:a,b\in\N_0\}$ is a basis
of $\mathbf J$. Moreover, $\mathbf J$ is a Koszul
algebra. Its prime and primitive ideals and its
representation theory were studied in \cite{Na}.
We recall below the description of its automorphism
group (see \cite[Theorem~3.1]{Shirikov}).

	\begin{theorem}\label{auto-Jordan-plane}
		The automorphisms of $\mathbf J$
		are of the form
		\[y_1 \mapsto \alpha y_1 + p(y_2), \quad y_2 \mapsto \alpha y_2, \]
		where $\alpha \in \Bbbk^\times$ and $p \in \Bbbk[t]$. Then, $\operatorname{Aut}\mathbf J\simeq \Bbbk[y_2]\rtimes \Bbbk^{\times}$. 
	\end{theorem}
	
	\subsection{Super Jordan plane}
	We now recall the notion of the super Jordan plane, which will be the object of study in this work.
	
	\begin{definition}
		The super Jordan plane is the associative $\Bbbk$-algebra $\B$
		generated by $x$ and $y$ subject to the relations
		\begin{align}
			x^2 &= 0, \label{eq:rels-B(V(-1,2))-1}\\
			yw-wy-xw &= 0, \label{eq:rels-B(V(-1,2))-2}
		\end{align}
		where $w=xy+yx$.
	\end{definition}
	The algebra $\mathcal{B}$ was introduced in \cite{aah-triang}. This algebra is the Nichols algebra of a $2$-dimensional braided vector space. 
	In what follows, we denote
	\[z:=yx.\]
	The following result was proved in  \cite{aah-triang}.
	
	\begin{proposition} \label{pr:sjp}The algebra $\cB$ is a Nichols algebra with $\GK \cB = 2$. The set of monomials
		\begin{equation}\label{basisPBW2}
			\{x^a z^b y^c , a\in\{0,1\},\; b,c\in\mathbb N_0\}
		\end{equation}
		is a Poincaré–Birkhoff–Witt basis for $\B$.
		\qed
	\end{proposition}
	Since $x$ is a nonzero PBW monomial and $x^2=0$, the algebra $\B$
	is not a domain. The defining relations give the following identities:
	\begin{align}
		wx &= xw,             \label{eq:relations1}\\
		y^2x &= xy^2+xz,      \label{eq:relations2}\\
		wy^2 &= (y^2-w)w,     \label{eq:relations3}\\
		zx &= 0.              \label{eq:zx-zero}
	\end{align}
	Indeed, since $w=xy+yx$, $z=yx$, and $x^2=0$, we have
	$wx=xw=xz$ and $zx=yx^2=0$. Expanding $yw-wy$ gives
	$y^2x-xy^2$, so the defining relation
	\eqref{eq:rels-B(V(-1,2))-2} is equivalent to
	\eqref{eq:relations2}. Writing the same defining relation as
	$wy=(y-x)w$, we obtain
	\[
	wy^2=(y-x)wy=(y-x)^2w=(y^2-w)w,
	\]
	which proves \eqref{eq:relations3}. Equation~\eqref{eq:relations3} can be rewritten as $y^2w-wy^2=w^2$. Thus, $y_1\mapsto y^2$ and $y_2\mapsto w$ define a surjective algebra homomorphism from the Jordan plane $\mathbf J$ onto the subalgebra of $\B$ generated by $y^2$ and $w$. To see that this map is injective, note that the monomials $y_2^a y_1^b$, with $a,b\geq0$, span $\mathbf J$ by its defining relation. Their images are $w^a y^{2b}$. Since $wx=xw=xz$, we have $xw^a=w^ax=xz^a$ for every $a\geq0$. Hence, left multiplication by $x$ sends these images to the distinct PBW monomials $xz^a y^{2b}$ of Proposition~\ref{pr:sjp}. It follows that the images are linearly independent and the homomorphism is an isomorphism.
	We henceforth identify $\mathbf J$ with its image
	$\Bbbk\langle y^2,w\rangle$ under this embedding.
	
	\begin{remark}
		The algebra $\B$ is $\N_0$-graded. Indeed, we put $\deg x=\deg y=1$. Then, the monomial $x^a z^b y^c$ has degree $a+2b+c$ and $\B=\bigoplus_{n\in \N_0}\B_n$ with $\B_n$ the homogeneous component of degree $n$.
	\end{remark}

	\subsection{Commutation formulas and PBW products}
	We record the following commutation formulas for later use.
	They can be proved by induction from the defining relations of $\B$. After reindexing the sums, identities \eqref{eqpropy2yx} and
	\eqref{eqpropy1yx} agree with formulas~(2.3) and~(2.4) of
	\cite{RS}, respectively.

	\begin{lemma}\label{propy2yx}
		Let $n\in\N$ and $k\in\N_0$. Then
		\begin{align}
			xy^{2k}x
			&=0, \label{x-square}\\[.4em]
			yz^n
			&=xz^{n-1}y^2+nxz^n, \label{propyyx}\\
			y^{2k}z^n
			&=\sum_{i=0}^k \binom{k}{i}
			(n)_i\,
			z^{n+i}y^{2(k-i)}. \label{eqpropy2yx}
		\end{align}
		If $k\geq1$, we also have
		\begin{equation}\label{eqpropy1yx}
			y^{2k-1}z^n
			=\sum_{i=0}^k \binom{k}{i}
			(n)_i\,
			xz^{n+i-1}y^{2(k-i)}.
		\end{equation}
		\qed
	\end{lemma}
	
	The following lemma gives the expansion of the product of two monomials with respect to the PBW basis of $\B$.
	
	\begin{lemma}\label{prop:product}
		Let $a,d\in\{0,1\}$ and $b,c,e,f\in\mathbb N_0$, and consider 
		\[
		Y:=(x^az^by^c)(x^dz^ey^f).
		\]
		Then, $Y$ is given by the following cases:
		
		\[Y=
		\begin{cases}
			0,
			& d=1,\ c=2k,\ a+b>0, \\[.4em]
			
			0,
			& d=0,\ e>0,\ c=2k+1,\ a+b>0, \\[.4em]
			
			x^az^b y^{c+f},
			& d=0,\ e=0, \\[.4em]
			
			\displaystyle
			\sum_{i=0}^k
			\binom{k}{i}(e+1)_i\,
			x^az^{b+e+1+i}y^{2(k-i)+f},
			& d=1,\ c=2k+1, \\[1em]
			
			\displaystyle
			\sum_{i=0}^k
			\binom{k}{i}(e+1)_i\,
			xz^{e+i}y^{2(k-i)+f},
			& d=1,\ c=2k,\ a+b=0, \\[1em]
			
			\displaystyle
			\sum_{i=0}^k
			\binom{k}{i}(e)_i\,
			x^az^{b+e+i}y^{2(k-i)+f},
			& d=0,\ e>0,\ c=2k, \\[1em]
			
			\displaystyle
			\sum_{i=0}^{k+1}
			\binom{k+1}{i}(e)_i\,
			xz^{e+i-1}y^{2(k+1-i)+f},
			& d=0,\ e>0,\ c=2k+1,\ a+b=0.
		\end{cases}\]
	\end{lemma}
	
	\begin{proof}
		For \(d=1,\ c=2k\) and \(a+b>0\), notice that
		\[Y=\big(x^a z^b y^{2k}\big)\big(xz^e y^{f}\big)=\big(x^a z^b y^{2k}x\big)\big(z^e y^{f}\big).\]
		Since $a+b>0$, the monomial $xy^{2k}x$ appears in $Y$. Thus, \eqref{x-square} implies that $Y=0$. Now we consider the case $d=0,\ e>0$ and $c=2k$. We have that
		\begin{align*}
			Y&=x^a z^b y^{2k}z^e y^{f}\\
			&\overset{\mathclap{\eqref{eqpropy2yx}}}{=}x^a z^b\Big[\sum_{i=0}^k  \binom{k}{i}\frac{(e-1+i)!}{(e-1)!} z^{e+i}  y^{2(k-i)}\Big] y^{f}\\
			&=\sum_{i=0}^k  \binom{k}{i}\frac{(e-1+i)!}{(e-1)!}x^a  z^{b+e+i}y^{2(k-i)+f}.
		\end{align*}
		The remaining cases follow similarly.
	\end{proof}
	
	\begin{proposition}\label{sonhotodoaluno}
		For every $n>0$, we have that
		\begin{equation}\label{eq:x21_power}
			w^n = (xy)^n +(yx)^n = (xy)^n +z^n.
		\end{equation}
		In particular, if  $a\in\{0,1\}$,  $b\in \N$ and $c\in \N_0$, then
		\begin{equation}\label{eq:rewrite_PBW}
			x^a w^b y^c
			=
			x^{a+1} z^{b-1} y^{c+1}+x^{a} z^{b} y^{c}.
		\end{equation}
	\end{proposition}
	
	\begin{proof}
		We proceed by induction on $n$. By definition of $w$, the identity \eqref{eq:x21_power} is true for $n=1$. Suppose that $n\geq 2$ and that \eqref{eq:x21_power} holds for $n-1$. Then
		\begin{align*}
			w^n&=(xy+yx)\big((xy)^{n-1}+z^{n-1}\big)=(xy)^{n}+z^n+xy^2xz^{n-2}+yx^2y(xy)^{n-2}.
		\end{align*}
		Since $x^2=0$ and $xy^2x=0$ by \eqref{x-square}, the identity \eqref{eq:x21_power} is proved. Now, suppose that $a\in\{0,1\}$,  $b\in \N$ and $c\in \N_0$. Then, by \eqref{eq:x21_power}, we have that
		\begin{align*}
			x^a w^b y^c&=x^a\big((xy)^b+z^b\big)y^c=x^a(xy)^by^c+x^az^by^c=x^{a+1}z^{b-1}y^{c+1}+x^az^by^c.
		\end{align*}
	\end{proof}

\section{The subgroup \texorpdfstring{$\mathbb L$}{L}}
	In this section, we characterize the subgroup $\mathbb L$ of $\Aut\B$ consisting of automorphisms that stabilize both one-dimensional vector spaces $\Span\{x\}$ and $\Span\{w\}$. We also describe its structure in terms of semidirect products and locally nilpotent derivations.

	\subsection{Automorphisms stabilizing \texorpdfstring{$\Span\{x\}$}{span(x)} and \texorpdfstring{$\Span\{w\}$}{span(w)}}
	In this subsection, we calculate the automorphisms of the super Jordan plane that stabilize $\Span \{x\}$ and $\Span \{ w\}$. Precisely, we describe the class of automorphisms $\sigma\in \Aut\B$ such that $\sigma(x)=\lambda_1x$ and $\sigma(w)=\lambda_2w$, for some $\lambda_1,\lambda_2\in \Bbbk^{\times}$. We obtain some auxiliary results and then we describe in Lemmas \ref{oddcomponent} and \ref{evencomponent} the homogeneous solutions of the following equations:
	\begin{align}\label{solution-u}
		xu+ux=0, & \quad [u,w]=0.
	\end{align}
	
	We first compute the commutators of certain PBW monomials
	with $w$.
	
	\begin{lemma}\label{commidentities}
		Let $b\in \N_0$ and $k\in \N$. Then the following identities hold:
		\begin{align}
			[xz^b y^{2k},w]&= \sum_{i=1}^{k}\binom{k}{i}i!xz^{b+i+1}y^{2(k-i)}, \label{commevenodd}  \\
			[xz^b y^{2k+1},w]&=\sum_{i=1}^{k}\binom{k}{i}i!xz^{b+i+1}y^{2(k-i)+1}.\label{commodd}
		\end{align}
		Moreover, if $b\in \N$ then we have that:
		\begin{align}
			[z^b y^{2k},w]&=\sum_{i=1}^{k}\binom{k}{i}i!z^{b+i+1}y^{2(k-i)},   \label{commeven}\\
			[z^b y^{2k+1},w]&=\sum_{i=1}^{k}\binom{k}{i}i!z^{b+i+1}y^{2(k-i)+1}.   \label{comm-new}
		\end{align}
	\end{lemma}
	
	\begin{proof}
		Firstly, we notice that
		\begin{align*}
			[y^{2k},xy]
			&=y^{2k}xy-xy^{2k+1}
			=y^{2k-1}(yx)y-xy^{2k+1}\\
			&\overset{\mathclap{\eqref{eqpropy1yx}}}{=}
			\, \Bigl(\sum_{i=0}^k \binom{k}{i}i!xz^{i}y^{2(k-i)}\Bigr)y
			-xy^{2k+1}\\
			&=\sum_{i=1}^k \binom{k}{i}i!xz^{i}y^{2(k-i)+1}.
		\end{align*}
		Consequently, we have
		\begin{align}\label{eq-colchete-1}
			[y^{2k},w] \, \,
			&  \overset{\mathclap{\eqref{eqpropy2yx}}}{=}
			\, \,  \sum_{i=1}^{k}\binom{k}{i}i!z^{i+1}y^{2(k-i)}
			+\sum_{i=1}^k \binom{k}{i}i!xz^{i}y^{2(k-i)+1}.
		\end{align}
		Since $[y,w]=y^2x-xy^2=xyx=xw$, the commutator identity gives
		\begin{align}
			[y^{2k+1},w]
			&=[y^{2k},w]y+y^{2k}xw \notag\\
			&=\sum_{i=1}^{k}\binom{k}{i}i!z^{i+1}y^{2(k-i)+1}
			+\sum_{i=1}^{k}\binom{k}{i}i!xz^iy^{2(k-i)+2}
			+y^{2k}xw. \label{eq-colchete-2}
		\end{align}
		
		By \eqref{eqpropy1yx}, the PBW expansion of
		$y^{2k}x=y^{2k-1}(yx)$ has a left factor $x$ in every term.
		Thus, $y^{2k}xw\in x\mathcal B$. Multiplying \eqref{eq-colchete-1} and \eqref{eq-colchete-2}
		on the left by $xz^b$ annihilates the second sum in both
		identities, since $xz^bx=0$. For the same reason, it also
		annihilates the additional term $y^{2k}xw$ in \eqref{eq-colchete-2}. Since $[x,w]=[z,w]=0$, we have
		$[xz^by^c,w]=xz^b[y^c,w]$ for all $c\in\N_0$.
		This yields \eqref{commevenodd} and \eqref{commodd}. Now suppose that $b\geq1$. Since $[z^by^c,w]=z^b[y^c,w]$,
		we multiply \eqref{eq-colchete-1} and \eqref{eq-colchete-2}
		on the left by $z^b$. As $z^bx=0$, the second sum in each
		identity vanishes, as does the term $y^{2k}xw$
		in \eqref{eq-colchete-2}. This proves \eqref{commeven}
		and \eqref{comm-new}, respectively.
	\end{proof}
	
	\begin{lemma} \label{wilian}
		Let $v\in \mathcal{B}$ and $\gamma\in\Bbbk$.
		The following conditions are equivalent:
		\begin{enumerate}[\rm (i)]
			\item
			$xv + vx = \gamma w$; \vspace{.1cm}
			
			\item
			there exists $u\in \mathcal{B}$ such that $v = \gamma y + u$, with $xu + ux = 0$.
		\end{enumerate}
	\end{lemma}
	
	\begin{proof}
		Notice that if $\gamma=0$ then the equivalence is trivial. Suppose that $\gamma\neq 0$.
		The implication (ii) $\Rightarrow$ (i) is immediate. Conversely, write $v$ as a sum of homogeneous components, that is, $v=P_0+P_1+\cdots+P_m$
		where each $P_i$ is homogeneous of degree $i$. Now, $xv+vx=\sum_{i=0}^m (xP_i+P_ix)$. Each summand $xP_i+P_ix$ is either zero or homogeneous of degree $i+1$. Since $\gamma w$ has degree two, if $\gamma \neq 0$, it follows that
		\[
		xP_i+P_ix=0,
		\quad\text{for every } i\neq 1.
		\]
		Thus, $xv+vx=xP_1+P_1x$. Since $P_1$ has degree one, $P_1=\alpha x+\beta y$ for some $\alpha,\beta\in\Bbbk$. Therefore, $xP_1+P_1x = \beta(xy+yx) = \beta w$. Hence, $\beta = \gamma$. Moreover, $u := (P_0 + P_2 + \cdots + P_m) + \alpha x$ satisfies the conditions required in (ii).
	\end{proof}
	
	\begin{lemma}\label{zero}
		Let $P_m=\sum_{b,c\geq0}\lambda_{b,c}xz^by^c\in\mathcal B$
		be a nonzero homogeneous element of degree $m$ such that
		$c$ is odd whenever $\lambda_{b,c}\neq0$.
		Then, $P_mx\neq0$.
	\end{lemma}
	
	\begin{proof}
		Write $
		P_m=\sum_{j=1}^l \lambda_{b_j,c_j}xz^{b_j}y^{c_j},
		$
		where $l\geq1$, all coefficients are nonzero and the pairs
		$(b_j,c_j)$ are distinct. Since $P_m$ is homogeneous of
		degree $m$, we have $1+2b_j+c_j=m$ for every $j$.
		Consequently, the exponents $c_j$ are distinct, and we may
		order them so that $c_1<\cdots<c_l$.
		By hypothesis, $c_j=2k_j+1$, with
		$0\leq k_1<\cdots<k_l$. Since $xz^{b_j}y^{c_j}x=xz^{b_j}y^{2k_j+1}x=xz^{b_j}y^{2k_j}z$, it follows from \eqref{eqpropy2yx}, with $n=1$, that
		\[
		xz^{b_j}y^{c_j}x
		=
		\sum_{i=0}^{k_j}
		\binom{k_j}{i}i!\,
		xz^{b_j+i+1}y^{2(k_j-i)}.
		\]
		Thus,
		$
		P_mx
		=
		\sum_{j=1}^l \lambda_{b_j,c_j}
		\sum_{i=0}^{k_j}
		\binom{k_j}{i}i!\,
		xz^{b_j+i+1}y^{2(k_j-i)}.
		$
		Since $k_1<\cdots<k_l$, the PBW monomial
		$xz^{b_l+1}y^{2k_l}$ occurs only for $j=l$ and $i=0$,
		with coefficient $\lambda_{b_l,c_l}\neq0$.
		Thus, $P_mx=\lambda_{b_l,c_l}xz^{b_l+1}y^{2k_l}+Q$, where $Q=0$ or $Q$ is homogeneous of degree $m+1$,
		and the powers of $y$ in the monomials appearing in $Q$
		are strictly smaller than $2k_l$.
		By the linear independence of the PBW monomials,
		$P_mx\neq0$, as desired.
	\end{proof}
	
	\begin{lemma}\label{oddcomponent}
		Let $P_{2n+1} \in \mathcal{B}$ be a homogeneous element of degree $2n+1$ and assume that $xP_{2n+1}+P_{2n+1}x=0$ and $[P_{2n+1},w]=0$. Then, there exists $a_n\in\Bbbk$ such that $P_{2n+1}=a_nxz^n$.
	\end{lemma}
	
	\begin{proof}
		By the PBW basis, every homogeneous element of degree $2n+1$ can be written as
		\[
		P_{2n+1}
		=
		\sum_{k=0}^{n}\lambda_kxz^{n-k}y^{2k}
		+
		\sum_{k=0}^{n}\mu_kz^{n-k}y^{2k+1}.
		\]
		We claim that all coefficients $\mu_k$ are zero and also $\lambda_k=0$ for all $k \geq 1$. Indeed, assume that some $\mu_k\neq 0$, and choose $k_0$ maximal with this property. Observe that
		\begin{align*}
			z^{n-k_0}y^{2k_0+1}x&=\,\,z^{n-k_0}y^{2k_0}z\,\,\overset{\mathclap{\eqref{eqpropy2yx}}}{=} \,\,
			z^{n-k_0} \sum_{i=0}^{k_0}  \binom{k_0}{i}i!z^{i+1}y^{2(k_0-i)}\\
			&= \,\,
			z^{n-k_0+1}y^{2k_0}+\sum_{i=1}^{k_0}  \binom{k_0}{i}i!z^{n-k_0+i+1}y^{2(k_0-i)}.
		\end{align*}
		
		On the other hand, $xz^{n-k}y^{2k}x=0$ for every $0\leq k\leq n$. Then, the terms involving the coefficients $\lambda_k$ do not contribute to $P_{2n+1}x$. Moreover, the PBW monomial $z^{n-k_0+1}y^{2k_0}$ arising from $z^{n-k_0}y^{2k_0+1}x$ cannot arise from a term corresponding to $\mu_k$ with $k<k_0$. Indeed, a monomial arising from $z^{n-k}y^{2k+1}x$ has the form
		\[
		z^{n-k+i+1}y^{2(k-i)},
		\quad 0\leq i\leq k,
		\]
		and equality with $z^{n-k_0+1}y^{2k_0}$ would imply
		$k-i=k_0$, which is impossible if $k<k_0$. By the maximality of $k_0$, we have $\mu_k=0$, for $k>k_0$. Thus, this monomial does not appear in $xP_{2n+1}$, since all monomials occurring there start with $x$. It follows from $xP_{2n+1}+P_{2n+1}x=0$, that $\mu_{k_0}=0$, which is a contradiction. Therefore, $\mu_k=0$, for every $0\leq k\leq n$.\vspace{.1cm}
		
		Now assume that some $\lambda_k\neq 0$ with $k\ge 1$, and choose $k_0$ maximal with this property. By \eqref{commevenodd} we have that
		$
			[xz^{n-k_0}y^{2k_0},w]= \sum_{i=1}^{k_0}\binom{k_0}{i}i!xz^{n-k_0+i+1}y^{2(k_0-i)},
	$
		which contains the element $k_0xz^{n-k_0+2}y^{2k_0-2}$. From $[P_{2n+1},w]=0$ and using the maximality of $k_0$, we conclude that $\lambda_{k_0}=0$ and we again have a contradiction. Therefore, $P_{2n+1}=\lambda_0xz^n$ as claimed.
	\end{proof}
	
	\begin{lemma}\label{evencomponent}
		Let $P_{2n}\in \mathcal{B}$ be a homogeneous element of degree $2n$, with $n\ge 1$. Assume that $xP_{2n}+P_{2n}x=0$
		and $[P_{2n},w]=0$. Then, there exists $b_n\in\Bbbk$ such that $P_{2n}=b_n(xz^{n-1}y-z^n)$.
	\end{lemma}
	
	\begin{proof}
		By the PBW basis, every homogeneous element of degree $2n$ can be written as
		$
		P_{2n}
		=
		\sum_{k=0}^{n}\lambda_kz^{n-k}y^{2k}
		+
		\sum_{k=0}^{n-1}\mu_kxz^{n-k-1}y^{2k+1}.
		$
		We first note that $\lambda_n=0$. Indeed, by \eqref{eqpropy2yx} and \eqref{eqpropy1yx}, we have
		$
		[y^{2n},w]=\sum_{i=1}^{n}\binom{n}{i}i!z^{i+1}y^{2(n-i)}+
		\sum_{i=1}^{n}\binom{n}{i}i!xz^i y^{2(n-i)+1}.
		$
		Observe that the monomial $xzy^{2n-1}$ appears in $[y^{2n},w]$. We claim that this monomial does not appear again in the bracket $[P_{2n},w]$. Indeed, the commutators $[z^{n-k}y^{2k},w]$, with $0\le k\le n-1$, produce PBW monomials with even powers of $y$. On the other hand, the PBW monomials occurring in $[xz^{n-k-1}y^{2k+1},w]$, with $0\le k\le n-1$, have highest power of $y$ equal to $2k-1\le 2n-3$. Since $[P_{2n},w]=0$, it follows that $\lambda_n=0$. \vspace{.1cm}
		
		Now, assume that some $\lambda_k\neq 0$ with $1\le k\le n-1$, and let $k_0$ be maximal with this property. Since $n-k_0\ge 1$, applying \eqref{commeven} we obtain
		\[
		[z^{n-k_0}y^{2k_0},w]
		=
		\sum_{i=1}^{k_0}\binom{k_0}{i}i!z^{n-k_0+i+1}y^{2(k_0-i)}.
		\]
		Hence, $k_0z^{n-k_0+2}y^{2k_0-2}$ appears in $[z^{n-k_0}y^{2k_0},w]$ and, by maximality of $k_0$,
		the element $k_0z^{n-k_0+2}y^{2k_0-2}$ does not appear in any other bracket of the form $[z^{n-k}y^{2k},w]$. Also, all the commutators
		$[xz^{n-k-1}y^{2k+1},w]$ begins with $x$. Thus, using that $[P_{2n},w]=0$, we conclude that $\lambda_{k_0}=0$ which is a contradiction. Therefore, $\lambda_k=0$, for all $1\le k\le n$.\vspace{.1cm}
		
		Suppose now that some $\mu_k\neq 0$ with $1\le k\le n-1$, and let $k_0$ be maximal with this property. By identity \eqref{commodd},
		$
		[xz^{n-k_0-1}y^{2k_0+1},w]
		=
		\sum_{i=1}^{k_0}\binom{k_0}{i}i!
		xz^{n-k_0+i}y^{2(k_0-i)+1}.
		$
		Then $k_0xz^{n-k_0+1}y^{2k_0-1}$ appears in $[xz^{n-k_0-1}y^{2k_0+1},w]$. By maximality of $k_0$,
		the element $k_0xz^{n-k_0+1}y^{2k_0-1}$ does not appear in any other bracket $[xz^{n-k-1}y^{2k+1},w]$.
		On the other hand,  the PBW monomials occurring in
		$[z^{n-k}y^{2k},w]$ do not begin with $x$. Hence, $[P_{2n},w]=0$ implies that $\mu_{k_0}=0$, a contradiction.\vspace{.1cm}
		
		By the preceding arguments,
		 $P_{2n}=\lambda_0 z^n+\mu_0xz^{n-1}y$. Thus, $xP_{2n}+P_{2n}x=0$ implies that
		\[
		0=xP_{2n}+P_{2n}x
		=
		\lambda_0xz^n+\mu_0xz^{n-1}yx
		=
		(\lambda_0+\mu_0)xz^n.
		\]
		Then $\lambda_0=-\mu_0$ and the result is proved.
	\end{proof}
	
	The next result characterizes the elements of $\B$ that satisfy a stronger version of \eqref{solution-u}.
	
	\begin{proposition}\label{supercase}
		Let $u\in \mathcal{B}$. The following statements are equivalent:
		\begin{enumerate}[\rm (i)]
			\item $xu=0=ux$ and $[u,w]=0$, \vspace{.1cm}
			\item  there exists a unique polynomial $p(t)\in\Bbbk[t]$ such that $u=xp(z)$.
		\end{enumerate}
	\end{proposition}
	
	\begin{proof}
		To prove (i) $\Rightarrow$ (ii), we first observe that, since
		$\mathcal B$ is $\mathbb N_0$-graded, we can write
		$u=\sum_{m\geq0}P_m$, where $P_m\in\mathcal B_m$ and only
		finitely many components are nonzero. Since $w$ has degree $2$,
		we have $[P_m,w]\in\mathcal B_{m+2}$ for every $m\geq0$.
		Thus, by the directness of the grading, $[u,w]=0$ implies
		that $[P_m,w]=0$ for all $m\geq0$.
		Similarly, since $x$ has degree $1$, the equalities
		$xu=0=ux$ imply that $xP_m=0=P_mx$ for all $m\geq0$.
		Using the PBW basis and $xP_m=0$, we conclude that $P_m$
		is a linear combination of monomials of the form $xz^by^c$
		with $2b+c+1=m$. Indeed, left multiplication by $x$
		annihilates the PBW monomials beginning with $x$ and sends
		the remaining monomials $z^by^c$ to the distinct PBW
		monomials $xz^by^c$.
		Suppose that $m$ is even. Every monomial $xz^by^c$
		occurring in the PBW expansion of $P_m$ satisfies
		$2b+c+1=m$, so $c$ is odd.
		If $P_m\neq0$, Lemma~\ref{zero} would give $P_mx\neq0$,
		contrary to $P_mx=0$. Therefore, $P_m=0$ for every even $m$,
		and hence
		$u=\sum_{n\geq0}P_{2n+1}$.
		Since $xP_{2n+1}+P_{2n+1}x=0$ and $[P_{2n+1},w]=0$,
		Lemma~\ref{oddcomponent} gives
		$P_{2n+1}=a_nxz^n$ for some $a_n\in\Bbbk$.
		Therefore, $u=xp(z)$, where
		$p(t)=\sum_{n\geq0}a_nt^n\in\Bbbk[t]$.
		The uniqueness of $p(t)$ follows from the linear independence
		of the PBW monomials $xz^n$.
		
		Conversely, the identities $x^2=zx=0$ and
		$[x,w]=[z,w]=0$ show that every element $u=xp(z)$,
		with $p(t)\in\Bbbk[t]$, satisfies
		$xu=0=ux$ and $[u,w]=0$.
		Thus, (ii) implies (i).
	\end{proof}
	
	Given $\mu \in \Bbbk^{\times}$ and $p(t),q(t)\in \Bbbk[t]$ we consider the map
	$\sigma_{p,q,\mu}$  defined on $\B$ by

	\begin{equation}\label{automorfismos-Adriano}
		\left\{
		\begin{aligned}
			\sigma_{p,q,\mu}(x)&=\mu x, \\[.2em]
			\sigma_{p,q,\mu}(y)&=\mu y+xp(z)+ xq(z)y-q(z)z.
		\end{aligned}
		\right.
	\end{equation}
	A direct computation shows that $\sigma_{p,q,\mu}$ is an automorphism of $\B$, for all $p,q\in \Bbbk[t]$ and $\mu\in \Bbbk^\times$. We will denote 
	\begin{align}\label{subgroup-automorfismos-B}
		\mathbb{L}:=\{\sigma_{p,q,\mu}:p(t), q(t)\in \Bbbk[t],\mu\in \Bbbk^\times\}.   
	\end{align}
	\begin{remark}\label{uniqueness}
		The parameters of an element in $\mathbb{L}$ are uniquely determined. Indeed, suppose that $\sigma_{p,q,\mu}=\sigma_{p',q',\mu'}$, with $p,p',q,q'\in \Bbbk[t]$ and $\mu,\mu'\in \Bbbk^\times$. Then, by evaluating at $x$, we obtain $\mu=\mu'$. Comparing the images of $y$ and using the PBW basis of $\mathcal B$,
		it follows that $p=p'$ and $q=q'$.
	\end{remark}
	
	The next result gives an intrinsic description of the subgroup $\mathbb L$, showing that it is determined by the preservation of two distinguished
	one-dimensional subspaces of $\mathcal B$. 
	
	\begin{theorem}\label{generalaut}
		The subgroup $\mathbb L$ consists precisely of those automorphisms
		of $\B$ that stabilize the one-dimensional subspaces
		$\Span\{x\}$ and $\Span\{w\}$.
	\end{theorem}
	
	\begin{proof}
		Suppose that $\sigma\in\Aut\B$ stabilizes $\Span\{x\}$ and
		$\Span\{w\}$. Then there exist $\mu,\lambda\in\Bbbk^\times$
		such that $\sigma(x)=\mu x$ and $\sigma(w)=\lambda w$.
		Let $v=\sigma(y)$. Since
		$\sigma(w)=\mu(xv+vx)=\lambda w$, we obtain
		$xv+vx=\lambda\mu^{-1}w$.
		By Lemma~\ref{wilian}, there exists $u\in\mathcal B$ such that
		\[
		v=\lambda\mu^{-1}y+u,
		\qquad xu+ux=0.
		\]
		Applying $\sigma$ to the relation $[y,w]=xw$ gives $\lambda[v,w]=\sigma([y,w])=\sigma(xw)=\mu\lambda xw$. Since $\lambda\neq0$, we have $[v,w]=\mu xw$. Therefore,
		$
		[u,w]
		=[v,w]-\lambda\mu^{-1}[y,w]
		=(\mu^2-\lambda)\mu^{-1}xw.
		$
		Write $u=\sum_{m\geq0}P_m$, where $P_m\in\mathcal B_m$
		and only finitely many components are nonzero.
		Since $xu+ux=0$ and
		$xP_m+P_mx\in\mathcal B_{m+1}$, the directness of the grading
		gives $xP_m+P_mx=0$ for every $m\geq0$.
		Moreover, $[P_m,w]\in\mathcal B_{m+2}$, while
		$(\mu^2-\lambda)\mu^{-1}xw\in\mathcal B_3$.
		Comparing homogeneous components, we obtain $[P_m,w]=0$, with $m\neq1$, and $[P_1,w]=(\mu^2-\lambda)\mu^{-1}xw$. Since $P_1\in\mathcal B_1$, we may write
		$P_1=ax+by$ for some $a,b\in\Bbbk$.
		The equality $xP_1+P_1x=bw=0$ implies $b=0$, since $w\neq0$.
		Thus $P_1=ax$, and $[x,w]=0$ gives $[P_1,w]=0$.
		Consequently, $(\mu^2-\lambda)\mu^{-1}xw=0$. Since $xw=xz$ is a nonzero PBW monomial, we conclude that
		$\lambda=\mu^2$.
			It follows that $u=\sigma(y)-\mu y$ and that
		$[P_m,w]=0$ and $xP_m+P_mx=0$ for every $m\geq0$. Since $P_0\in\mathcal B_0=\Bbbk$, write $P_0=c$. Then, $xP_0+P_0x=2cx=0$, and characteristic zero gives $c=0$. Hence, $P_0=0$. Applying Lemma~\ref{oddcomponent} to the odd components and
		Lemma~\ref{evencomponent} to the positive even components,
		we obtain $P_{2n+1}=a_nxz^n$, for $n\geq0$ and $P_{2n}=b_n(xz^{n-1}y-z^n)$ for $n\geq1$, for suitable scalars $a_n,b_n\in\Bbbk$. Thus, $u=xp(z)+xq(z)y-q(z)z$, where
		\[
		p(t)=\sum_{n\geq0}a_nt^n,
		\quad
		q(t)=\sum_{n\geq1}b_nt^{n-1}.
		\]
		Both sums are finite, so $p,q\in\Bbbk[t]$.
		Therefore, $\sigma=\sigma_{p,q,\mu}\in\mathbb L$. Conversely, direct substitution gives
		$\sigma_{p,q,\mu}(x)=\mu x$ and
		$\sigma_{p,q,\mu}(w)=\mu^2w$.
		Hence, every element of $\mathbb L$ stabilizes both
		$\Span\{x\}$ and $\Span\{w\}$.
	\end{proof}
	
	\begin{remark} 
		The automorphisms in $\mathbb L$ do not, in general, preserve the Jordan
		plane subalgebra of $\mathcal B$. Indeed, take $p=0$, $q=1$ and $\mu=1$. Then
		$
		\sigma_{0,1,1}(x)=x$ and $\sigma_{0,1,1}(y)=y+xy-yx$.	Moreover, $\sigma_{0,1,1}(w)=w$ and
		$\sigma_{0,1,1}(y^2)=y^2+w^2-xw$. The term $xw$ has odd degree and is nonzero. Since the subalgebra generated
		by $y^2$ and $w$ is concentrated in even degrees, this term does not belong
		to that subalgebra. Therefore $\sigma_{0,1,1}$ does not restrict to an automorphism
		of the Jordan plane $\mathbf J$. On the other hand, the subfamily $q=0$ does preserve the Jordan plane. Indeed, for $\sigma_{p,0,\mu}$ one has $\sigma_{p,0,\mu}(w)=\mu^2w$ and $\sigma_{p,0,\mu}(y^2)=\mu^2y^2+\mu wp(w)$. Thus, the restriction of $\sigma_{p,0,\mu}$ to $\mathbf J$ has the form described in Theorem \ref{auto-Jordan-plane}.
	\end{remark}

	\subsection{The group structure of \texorpdfstring{$\mathbb L$}{L}} In this subsection, we realize $\mathbb L$ as a semidirect product built from
	a polynomial group $\mathbb T$ and the natural rescaling action of $\Bbbk^\times$. We recall from \cite{AD} that $\B$ has the $\mathbb{Z}_2$-grading
	\[
	\B_{\bar 0}=\bigoplus_{n\ \mathrm{even}} \B_n,
	\quad\qquad\quad
	\B_{\bar 1}=\bigoplus_{n\ \mathrm{odd}} \B_n.
	\]
	Following \cite{AD}, an automorphism $\sigma\in \Aut \B$ is called a {\it super-automorphism} of $\B$ if $\sigma(\B_{\bar 0})\subset \B_{\bar 0}$ and $\sigma(\B_{\bar 1})\subset \B_{\bar 1}$. Denote the group of super-automorphisms of $\B$ by $\Aut_{\mathrm{super}}\B$. Given $\mu \in \Bbbk^\times$ and $p(t)\in \Bbbk[t]$, consider $\sigma_{p,\mu}$ defined by
	\begin{equation}\label{super-AD}
		\sigma_{p,\mu}(x)=\mu x
		\quad\text{and}\quad
		\sigma_{p,\mu}(y)=\mu y+p(w)x.
	\end{equation}
	It follows from Lemma 5.2 of \cite{AD} that $\sigma_{p,\mu}$ is a super-automorphism of $\B$. Moreover, the group $\Aut_{\mathrm{super}}\B$ was computed in \cite{AD}. For the reader's convenience, we recall it below.
	
	\begin{theorem}[\cite{AD}, Theorem~5.3] The group $\Aut_{\mathrm{super}}\B$ is given by the elements $\sigma_{p,\mu}$ defined in \eqref{super-AD}, with $\mu\in \Bbbk^\times$ and $p(t)\in \Bbbk[t]$. \qed
	\end{theorem}

	\begin{corollary} 
		The group $\Aut_{\mathrm{super}}\B$ is given by the elements $\sigma_{u,\mu}$ defined by
		\begin{equation*}
			\sigma_{u,\mu}(x)=\mu x
			\quad\text{and}\quad
			\sigma_{u,\mu}(y)=\mu y+u,
		\end{equation*}    
		where $\mu\in\Bbbk^\times$ and $u\in\B$ satisfies the
		equivalent conditions of Proposition~\ref{supercase}.
	\end{corollary}
	\begin{proof}
		Given $p(w)\in \Bbbk[w]$, notice that $p(w)x\,\,\overset{\mathclap{\eqref{eq:relations1}}}{=}\,\,xp(w)\,\,\overset{\mathclap{\eqref{eq:x21_power}}}{=}\,\,xp(z)$.
		Then, the result is an immediate consequence of the previous theorem.
	\end{proof}

	\begin{remark}\label{obs-dumas-nicolas}
		Notice that $\Aut_{\mathrm{super}}\B\subset \mathbb{L}$. Indeed, using \eqref{eq:x21_power}, we see that $\sigma_{p,\mu}=\sigma_{p,0,\mu}$, for all $p(t)\in \Bbbk[t]$ and $\mu\in \Bbbk^\times$. 
	\end{remark}
	
	Let $\sigma,\tau\in \mathbb{L}$. By Theorem \ref{generalaut}, $\sigma\tau^{-1}\in \mathbb{L}$ and consequently $\mathbb{L}$ is a subgroup of $\Aut \B$. In order to give the composition law of $\mathbb{L}$, we prove the following auxiliary result.
	
	\begin{lemma} \label{eqs12}We have the following identities:
		\begin{align}
			\label{eqs1}
			\sigma_{p,q,\mu}(z^n)&=\mu^{2n}z^n+\mu^{2n-1}xq(z)z^n, && n\ge 1, \\[.3em]
			\label{eqs2}
			x\sigma_{p,q,\mu}(z^n)&=\mu^{2n}xz^n, &&n\ge 0.
		\end{align}
	\end{lemma}
	
	\begin{proof}
		We will prove \eqref{eqs1} and \eqref{eqs2} by induction on $n$. For $n=1$, notice that
		\begin{align*}
			\sigma_{p,q,\mu}(z)
			&=\big(\mu y+xp(z)+xq(z)y-q(z)z\big)\mu x=
			\mu^2z+\mu xq(z)z,
		\end{align*}
		since $xp(z)x=0$ and $x^2=0$. This proves
		\eqref{eqs1} for $n=1$. Multiplying \eqref{eqs1} on the left by $x$, we obtain that \eqref{eqs2} also holds for $n=1$. The case $n=0$ in \eqref{eqs2} is trivial. Assume that \eqref{eqs1} and \eqref{eqs2} hold for $n\ge 1$. Then,
		\begin{align*}
			\sigma_{p,q,\mu}(z^{n+1})
			&=
			\sigma_{p,q,\mu}(z^n)\sigma_{p,q,\mu}(z)=
			\big(\mu^{2n}z^n+\mu^{2n-1}xq(z)z^n\big)
			\big(\mu^2z+\mu xq(z)z\big).
		\end{align*}
		Since $z^n x=0$, we obtain
		$\sigma_{p,q,\mu}\big(z^{n+1}\big)
		=
		\mu^{2n+2}z^{n+1}
		+\mu^{2n+1}xq(z)z^{n+1},
		$
		which proves \eqref{eqs1} for $n+1$. Multiplying on the left by $x$, we also obtain
		$
		x\sigma_{p,q,\mu}(z^{n+1})
		=
		\mu^{2n+2}xz^{n+1},
		$
		which proves \eqref{eqs2} for $n+1$, as desired. 
	\end{proof}
	
	\begin{proposition}\label{composition}
		Let $\mathbb{L}$ be as in \eqref{subgroup-automorfismos-B}. The following statements hold: \vspace{.1cm}
		\begin{enumerate} [\rm (i)]
			\item if $\sigma_{p,q,\mu},\sigma_{r,s,\nu} \in \mathbb{L}$, then $\sigma_{p,q,\mu}\sigma_{r,s,\nu}=\sigma_{\bar p,\bar q,\mu\nu}$,
			where \vspace{.2cm}
			\begin{equation*}
				\left\{
				\begin{aligned}
					\bar p(t)&=\nu p(t)+\mu r(\mu^2 t)-2\mu tq(t)s(\mu^2 t),\\
					\bar q(t)&=\nu q(t)+\mu^2 s(\mu^2 t),
				\end{aligned}
				\right.
			\end{equation*}
			\item every $\sigma_{p,q,\mu}\in \mathbb{L}$ is invertible with inverse given by $\sigma^{-1}_{p,q,\mu}=\sigma_{r,s,\mu^{-1}}$ where
			\begin{equation*}
				\left\{
				\begin{aligned}
					r(t)&=-\mu^{-2}p(\mu^{-2}t)-2\mu^{-5}tq^2(\mu^{-2}t),\\
					s(t)&=-\mu^{-3}q(\mu^{-2}t).
				\end{aligned}
				\right.
			\end{equation*}
		\end{enumerate}
	\end{proposition}
	
	\begin{proof}
		Let $f(t)=\sum_{i=0}^m \lambda_i t^i\in \Bbbk[t]$. By \eqref{eqs2}, $\sigma_{p,q,\mu}(xf(z))=\mu xf(\mu^2z)$.
		On the other hand, we have by \eqref{eqs1} that
		\begin{align*}
			\sigma_{p,q,\mu}(f(z)z)
			&=
			\sum_{i=0}^m \lambda_i\,\sigma_{p,q,\mu}(z^{i+1})=
			\sum_{i=0}^m \lambda_i
			(
			\mu^{2i+2}z^{i+1}
			+\mu^{2i+1}xq(z)z^{i+1}
			)\\
			&=
			\mu^2f(\mu^2z)z
			+\mu xq(z)f(\mu^2z)z.
		\end{align*}
		Therefore,
		\begin{align*}
			\sigma_{p,q,\mu}\big([x, f(z)y]\big)&=
			\sigma_{p,q,\mu}(xf(z))\sigma_{p,q,\mu}(y)
			-\sigma_{p,q,\mu}(f(z)z)\\[.4em]
			&=
			\mu xf(\mu^2z)
			(\mu y+xp(z)+xq(z)y-q(z)z)\\[.4em]
			&
			-\mu^2f(\mu^2z)z
			-\mu xq(z)f(\mu^2z)z\\[.4em]
			&\overset{(\ast)}{=}\mu^2xf(\mu^2z)y
			-\mu xf(\mu^2z)q(z)z\\
			&-\mu^2f(\mu^2z)z
			-\mu xq(z)f(\mu^2z)z \\[.4em]
			&=\mu^2\big([x,f(\mu^2z)y]\big)-2\mu xzf(\mu^2z)q(z),
		\end{align*}
		where $(\ast)$ holds because $xf(\mu^2z)x=0$. Hence, we have
		\begin{align*}
			\sigma_{p,q,\mu} \sigma_{r,s,\nu}(y)&= \nu \sigma_{p,q,\mu}(y)
			+\sigma_{p,q,\mu}\big(xr(z)\big)
			+\sigma_{p,q,\mu}\big([x,s(z)y]\big) \\
			&=\nu\big(\mu y+xp(z)+[x,q(z)y]\big)+\mu xr(\mu^2z)\\
			&+\mu^2\big([x,s(\mu^2z)y]\big)-2\mu xz\,s(\mu^2z)q(z)\\
			&=\mu\nu y+x\big(\nu p(z)+\mu r(\mu^2z)-2\mu zq(z)s(\mu^2z)\big)\\
			&+x\big(\nu q(z)+\mu^2 s(\mu^2z)\big)y-\big(\nu q(z)+\mu^2 s(\mu^2z)\big)z,
		\end{align*}
		and (i) is proved. Notice that (ii) is an immediate consequence of (i).
	\end{proof}
	
	Let $A_1,A_2$ be abelian groups and $c:A_2\times A_2\to A_1$ an additive $2$-cocycle, that is,
	\[
	c(u,v)+c(u+v,w)=c(v,w)+c(u,v+w),
	\quad u,v,w\in A_2.
	\]
	We recall that $A_1\times_c A_2$ is a group with product given by
	\[
	(a,u)\cdot(b,v)=\big(a+b+c(u,v),\, u+v\big), \quad a,b\in A_1,\,u,v\in A_2.
	\]
	In particular, take $A_1=A_2=\Bbbk[t]$ and
	$c(q,s)=-2tqs$ for all $q,s\in\Bbbk[t]$. Since this cocycle is symmetric, the group
	$\mathbb T:=\Bbbk[t]\times_c\Bbbk[t]$, with product
	\begin{equation}\label{group-T}
		(p,q)(r,s)=(p+r-2tqs,q+s),
	\end{equation}
	is abelian. Moreover, $\Bbbk^\times$ acts on $\mathbb{T}$ by
	\begin{align}\label{action-field}
		\mu\cdot (p(t),q(t))=(p(\mu^2t),\mu q(\mu^2t)),\quad p,q\in \Bbbk[t],\,\,\mu\in \Bbbk^\times.
	\end{align}
	
	We are now ready to prove the following result.
	
	\begin{proposition}\label{subgr-autom}
		Let $\mathbb L$ and $\mathbb T$ be as above. Then $
		\mathbb L\simeq\mathbb T\rtimes\Bbbk^\times$
		via the group isomorphism
		$
		\psi:\mathbb T\rtimes\Bbbk^\times
		\to\mathbb L$ given by 
		\[\psi((p,q),\mu)=\sigma_{\mu p,\mu q,\mu},\qquad p,q\in\Bbbk[t],\,\,\, \mu\in\Bbbk^\times.\]
	\end{proposition}
	\begin{proof}
		Since $\Bbbk^\times$ acts on the abelian group $\mathbb{T}$ (see \eqref{action-field}), we can consider the semidirect product $\mathbb{T}\rtimes \Bbbk^\times$ and a direct computation shows that $\psi$ defined in the statement of the proposition is a group morphism. Finally, we notice that $\psi$ is invertible. Precisely, the inverse map is given by
		\[
		\psi^{-1}(\sigma_{a,b,\mu})
		=
		\big(
		(\mu^{-1}a,\mu^{-1}b),\mu\big), \quad a,b\in \Bbbk[t],\,\,\mu \in \Bbbk^\times,
		\]
		and the result follows.
	\end{proof}
	
	\begin{remark}\label{obs-T}
		Let $\mathbb{T}$ be the group defined in \eqref{group-T}. Note that
		\[
		\psi_0:\mathbb{T}\to \big\{\sigma_{p,q,1}\,:\, p,q\in \Bbbk[t]\big\},\quad (p,q)\mapsto \sigma_{p,q,1},
		\]
		is a group isomorphism. We shall see later that every automorphism $\sigma_{p,q,1}$ can be written as a composition of exponentials of locally nilpotent derivations.
	\end{remark}

	\subsection{Inner automorphisms and super-automorphisms in \texorpdfstring{$\mathbb L$}{L}}	
	In this subsection, we characterize the inner automorphisms among the elements
	$\sigma_{p,q,\mu}\in\mathbb L$, where $\mathbb L$ is the subgroup of
	$\Aut\mathcal B$ defined in \eqref{subgroup-automorfismos-B}. Using this
	description, we prove that $\mathbb L$ is a semidirect product of
	$\Inn\mathcal B\cap\mathbb L$ by the group of super-automorphisms of $\mathcal B$.
	
	\begin{lemma}\label{first-inclusion}
		Let $q(t)\in \Bbbk[t]$. Then, $\sigma=\sigma_{-tq(t)^2,q(t), 1}$ is an inner automorphism of $\mathcal B$. More precisely, we have that $\sigma=\Ad_{u_q}$,
		where $u_q=1+xq(z)$.
	\end{lemma}
	
	\begin{proof}
		Since $(xq(z))^2=0$, the element $u_q=1+xq(z)$ is invertible and also $u_q^{-1}=1-xq(z)$. Moreover, $u_qx=x=xu_q^{-1}$ and consequently $\Ad_{u_q}(x)=x$. On the other hand,
		\begin{align*}
			\Ad_{u_q}(y)
			&=
			(1+xq(z))y
			(1-xq(z))=
			y+xq(z)y
			-yxq(z)
			-xq(z)yxq(z)\\
			&=
			y-xzq(z)^2+xq(z)y-q(z)z=\sigma_{-tq(t)^2,q(t),1}(y),
		\end{align*}
		and the result follows.
	\end{proof}
	
	We determine the centralizer of $x$, which will be used to
	characterize the inner automorphisms contained in $\mathbb L$.
	Let
	\[
	u=y^2-xy,
	\qquad
	\mathcal J=\Bbbk\langle u,w\rangle.
	\]
	Both $u$ and $w$ commute with $x$ and they satisfy
	$uw-wu=w^2$. Together with the linear independence established
	in the proof of Lemma~\ref{lemma-centralizador-x1}, this
	relation identifies $\mathcal J$ with a copy of the Jordan
	plane contained in $C_{\B}(x)$.

	\begin{lemma}\label{lemma-centralizador-x1}
		The centralizer of $x$ in $\B$ is $C_{\B}(x)=\mathcal J \oplus x\mathcal J$.
	\end{lemma}
	
	\begin{proof}
		As observed above, both $u=y^2-xy$ and $w$ commute with $x$. Then, $\mathcal J\subseteq C_{\mathcal B}(x)$. Since every element of $\mathcal J$ commutes with $x$, so does every element of $x\mathcal J$. Moreover, $\mathcal J$ is concentrated in even degrees, whereas $x\mathcal J$ is concentrated in odd degrees. Therefore, the sum is direct and $\mathcal J\oplus x\mathcal J \subseteq C_{\mathcal B}(x)$.
		
		We prove the reverse inclusion by comparing the dimensions of the homogeneous components. Since $x$ is homogeneous, its centralizer is a graded subspace:
		\[
		C_{\mathcal B}(x)
		=
		\bigoplus_{n\geq0}C_{\mathcal B}(x)_n,
		\quad
		C_{\mathcal B}(x)_n
		:=
		C_{\mathcal B}(x)\cap\mathcal B_n.
		\]

		\noindent{\it Claim:} For every $m\geq 0$, one has $\dim C_{\B}(x)_{2m}\leq m+1$ and $\dim C_{\B}(x)_{2m+1}\leq m+1$. \vspace{.1cm}
		
		Consider the linear map $D_x:\mathcal B_n\longrightarrow \mathcal B_{n+1}$ given by $D_x(a)=ax-xa$, for all $a\in \B_n$.
		Then, $C_{\mathcal B}(x)_n=\ker(D_x|_{\mathcal B_n})$. Now, observe that
		\[
		\mathcal B_{2m}=\Span\{z^i y^{2m-2i}:0\leq i\leq m\}\oplus\Span\{xz^i y^{2m-2i-1}:0\leq i\leq m-1\}.
		\]
		Thus, $\dim\mathcal B_{2m}=2m+1$. Let $a_i:=z^iy^{2(m-i)}\in \B_{2m}$, with $i=1,\ldots,m$. Then, $D_x(a_i)=-xz^iy^{2(m-i)}$. So, these $m$ images are linearly independent. Therefore $\operatorname{rank}(D_x|_{\mathcal B_{2m}})\geq m$ and consequently
		\[
		\dim C_{\mathcal B}(x)_{2m}=\dim \ker(D_x|_{\mathcal B_{2m}})\leq (2m+1)-m=m+1.
		\]
		Similarly, $\mathcal B_{2m+1}$ has dimension $2m+2$. Fix $b_i:=z^iy^{2(m-i)+1}\in \B_{2m+1}$, for all $i=0,\ldots,m$. Then, the image of $b_i$ under $D_x$ has, among the PBW monomials with no left factor $x$, the leading term $z^{i+1}y^{2(m-i)}$. Thus, $\operatorname{rank}(D_x|_{\mathcal B_{2m+1}})\geq m+1$. Hence, $\dim_{\Bbbk}C_{\mathcal B}(x)_{2m+1}\leq (2m+2)-(m+1)=m+1$ and the claim is proved.
		\vspace{.1cm}
		
		On the other hand, the elements $u^{m-r}w^r$, $0\leq r\leq m$,
		belong to $C_{\mathcal B}(x)_{2m}$. Moreover, the PBW expansion of
		$u^{m-r}w^r$ contains the term $z^r y^{2(m-r)}$
		with coefficient $1$. More precisely, this is the term with the smallest
		power of $z$ among the terms not beginning with $x$. These
		terms are distinct for $0\leq r\leq m$, and hence the elements
		$u^{m-r}w^r$ are linearly independent.
		Similarly, the elements $xu^{m-r}w^r$, $0\leq r\leq m$, belong to $C_{\mathcal B}(x)_{2m+1}$. Their PBW expansions have distinct
		terms $xz^r y^{2(m-r)}$
		with coefficient $1$, and consequently they are linearly independent. Together with the dimension estimates of the claim, this gives $\dim C_{\mathcal B}(x)_{2m}=\dim C_{\mathcal B}(x)_{2m+1}=m+1$. Therefore, the desired identity holds in every homogeneous component, which completes the proof.
	\end{proof}

	Lemma \ref{first-inclusion} provides a family of inner automorphisms contained in $\mathbb L$. Using the description of the centralizer of $x$ obtained in Lemma \ref{lemma-centralizador-x1}, we prove that these are precisely the inner automorphisms belonging to $\mathbb L$.
	
	\begin{proposition}\label{prop-inner-intersection-Sigma}
		Let $\mathbb{L}$ be the subgroup of $\Aut \B$ defined in \eqref{subgroup-automorfismos-B}. Then
		\[
		\Inn\mathcal B\cap \mathbb{L}=\left\{\sigma_{-tq(t)^2, q(t), 1} \, : \,q(t)\in\Bbbk[t]\right\}.
		\]
	\end{proposition}
	
	\begin{proof}
		The inclusion from right to left follows from Lemma~\ref{first-inclusion}.
		Conversely, suppose that $\sigma_{p,q,\mu}$ is inner, that is,
		$\sigma_{p,q,\mu}=\Ad_v$ for some $v\in\B^\times$.
		Let $\pi:\B\to\B/(x)\simeq\Bbbk[y]$ be the canonical projection.
		Then, $\pi(v)=\lambda\in\Bbbk^\times$, and
		\[
		\pi(\Ad_v(b))=\lambda\pi(b)\lambda^{-1}=\pi(b),
		\quad b\in\B.
		\]
		Hence, $\pi(y)=\pi(\sigma_{p,q,\mu}(y))=\mu\pi(y)$, and
		consequently $\mu=1$. Thus $vx=xv$, so $v\in C_{\mathcal B}(x)$.
		By Lemma~\ref{lemma-centralizador-x1}, there exist
		$r,s_0\in\mathcal J$ such that $v=r+xs_0$. Notice that
		$x\mathcal J$ is a two-sided ideal of $C_{\mathcal B}(x)$ and
		$C_{\mathcal B}(x)/x\mathcal J\simeq\mathcal J$.
		Since $v^{-1}$ also commutes with $x$, the image of $v$ in this
		quotient is a unit. Under the above identification, this image is
		$r$. Since the units of $\mathcal J$ are the nonzero scalars,
		$r\in\Bbbk^\times$. Using $\Ad_v=\Ad_{r^{-1}v}$, we may therefore
		assume that $v=1+xs_0$, with $s_0\in\mathcal J$.
		
		We write $s_0$ in the PBW basis of the Jordan plane $\mathcal J$, that is, $s_0=\sum_{i,j\geq0}\lambda_{i,j}u^iw^j$, where $\lambda_{i,j}\in\Bbbk$.
		Since $u$ commutes with $x$ and $xw^j=xz^j$ for every $j\geq0$,
		we obtain $xu^iw^j=u^ixw^j=u^ixz^j=xu^iz^j$. Thus, writing $s:=\sum_{i,j\geq0}\lambda_{i,j}u^iz^j$, we have $xs_0=xs$ and hence $v=1+xs$. We claim that $\lambda_{i,j}=0$ for every $i>0$ and $j\geq0$.
		Suppose, by contradiction, that this is not the case. Choose $a>0$
		maximal such that $\lambda_{a,j}\neq0$ for some $j$, and choose $b$
		maximal such that $\lambda_{a,b}\neq0$. Since $xs=xs_0$ and $s_0$
		commutes with $x$, we have $(xs)^2=(xs_0)^2=x^2s_0^2=0$. Hence, $v^{-1}=1-xs$. Using $z=yx$, we obtain
		\begin{align}
			\Ad_v(y)-y
			&=(1+xs)y(1-xs)-y=xsy-yxs-xsyxs
			=xsy-zs-xszs. \label{identity-inner}
		\end{align}
		Observe that the first and the last summands in
		$xsy-zs-xszs$ both lie in $x\mathcal B$. Thus, they only produce
		PBW monomials of the form $xz^l y^m$. In particular, they cannot
		cancel any PBW monomial of the form $z^l y^m$. It remains to inspect the term $-zs$. The chosen summand
		$\lambda_{a,b}u^az^b$ of $s$ contributes
		$-\lambda_{a,b}zu^az^b$ to $-zs$. Since $u=y^2-xy$, induction
		on $l$ gives
		\begin{align}\label{rec-used}
			u^l=y^{2l}-l\,xy^{2l-1},
			\quad l\geq1.
		\end{align}
		In particular, since $a>0$, $y^{2a}$ is the only term in the PBW
		expansion of $u^a$ which does not belong to $x\mathcal B$.
		It also follows by induction that, for all $l,m\geq0$,
		\begin{align}\label{comu-for-impor}
			zy^{2l}z^m
			=\sum_{k=0}^{l}\binom{l}{k}(m)_k\,z^{m+1+k}y^{2l-2k},
		\end{align}
		where $(m)_0=1$ and $(m)_k=m(m+1)\cdots(m+k-1)$ for $k\geq1$.
		Taking $l=a$ and $m=b$, we see that the term $z^{b+1}y^{2a}$
		appears on the right-hand side of \eqref{comu-for-impor} with
		coefficient $1$, and all the other terms have strictly smaller
		$y$-degree. Since $zx=0$, equation \eqref{rec-used} gives
		$zu^az^b=zy^{2a}z^b$. Therefore, the contribution of
		$\lambda_{a,b}u^az^b$ to $-zs$ contains the PBW monomial
		$-\lambda_{a,b}z^{b+1}y^{2a}$. This monomial cannot be canceled by any other contribution to
		$\Ad_v(y)-y$. First, as observed above, $xsy$ and $-xszs$ belong
		to $x\mathcal B$, whereas $z^{b+1}y^{2a}$ does not. Second,
		consider another summand $\lambda_{i,j}u^iz^j$ of $s$ with
		$\lambda_{i,j}\neq0$. If $i<a$, then every term coming from
		$zu^iz^j$ has $y$-degree at most $2i<2a$. If $i=a$, then
		$j<b$ by the choice of $b$, and the term of $y$-degree $2a$ is
		$z^{j+1}y^{2a}$, which is distinct from $z^{b+1}y^{2a}$.
		Finally, there is no summand with $i>a$, by the maximality of $a$.
		Hence the coefficient of $z^{b+1}y^{2a}$ in $\Ad_v(y)-y$ is
		exactly $-\lambda_{a,b}\neq0$.
		
		On the other hand, since $\Ad_v=\sigma_{p,q,1}$, we have $\Ad_v(y)-y=xp(z)+xq(z)y-q(z)z$. The PBW monomials not belonging to $x\mathcal B$ on the
		right-hand side come only from $-q(z)z$ and are therefore pure
		powers of $z$. Thus, no monomial of the form $z^{b+1}y^{2a}$,
		with $a>0$, can occur. This contradiction proves that
		$\lambda_{i,j}=0$ for every $i>0$ and $j\geq0$.
		Consequently, $s=q_0(z)$ for some $q_0(t)\in\Bbbk[t]$.
		By Lemma~\ref{first-inclusion},
		\[
		\sigma_{p,q,1}
		=\Ad_v
		=\Ad_{1+xq_0(z)}
		=\sigma_{-tq_0(t)^2,q_0(t),1},
		\]
		which completes the proof.
	\end{proof}
	
Using the preceding result we show that $\mathbb L$ is an
extension of $\Aut_{\mathrm{super}}\B$ by
$\Inn\B\cap\mathbb L$.
	
	\begin{corollary}\label{cor-quotient-L-inner}
		Let $\mathbb L$ be as in \eqref{subgroup-automorfismos-B} and $N:=\Inn\mathcal B\cap \mathbb L$. Then there exists a split short exact sequence of groups
		\[1 \longrightarrow
		N\longrightarrow
		\mathbb L
		\overset{\Phi}{\longrightarrow}
		\Aut_{\mathrm{super}}\mathcal B
		\longrightarrow
		1,\]
		where $\Phi(\sigma_{p,q,\mu})=\sigma_{p+\mu^{-1}tq(t)^2,\mu}$. In particular, $ \mathbb L\simeq N\rtimes \Aut_{\mathrm{super}}\mathcal B$.
	\end{corollary}
	
	\begin{proof}
		Let $\sigma_{p,q,\mu},\sigma_{r,s,\nu}\in\mathbb L$. By Proposition \ref{composition}, their product is $\sigma_{\bar p,\bar q,\mu\nu}$, where
		\[
		\bar p(t)=\nu p(t)+\mu r(\mu^2t)-2\mu tq(t)s(\mu^2t),\quad\bar q(t)=\nu q(t)+\mu^2s(\mu^2t).
		\]
		A direct computation gives
		\[
		\bar p(t)+(\mu\nu)^{-1}t\bar q(t)^2=\nu\big(p(t)+\mu^{-1}tq(t)^2\big)+\mu\big(r(\mu^2t)+\nu^{-1}\mu^2t\,s(\mu^2t)^2\big).
		\]
		This is precisely the composition rule for the super-automorphisms $\sigma_{p+\mu^{-1}tq(t)^2,\mu}$ and $\sigma_{r+\nu^{-1}ts(t)^2,\nu}$. Then, $\Phi$ is a group homomorphism. The map $\Phi$ is surjective, since $\Phi(\sigma_{l,0,\mu})=\sigma_{l,\mu}$ for all $l(t)\in\Bbbk[t]$ and  $\mu\in\Bbbk^\times$. Moreover,
		$\sigma_{p,q,\mu}\in\ker\Phi$ if and only if $\mu=1$ and $p(t)+tq(t)^2=0$. Thus, $\ker\Phi=\{\sigma_{-tq(t)^2,q(t),1} \, \vert \, q(t)\in\Bbbk[t]\}$. By Proposition \ref{prop-inner-intersection-Sigma}, $\ker\Phi=N$ and the sequence is exact. Moreover, the map $\Psi:\Aut_{\mathrm{super}}\mathcal B\longrightarrow \mathbb L$ given by
		$\Psi(\sigma_{r,\mu})=\sigma_{r,0,\mu}$, $r\in \Bbbk[t], \,\,\mu\in \Bbbk^\times$,
		is a group homomorphism and satisfies
		$\Phi\circ\Psi=\operatorname{id}_{\Aut_{\mathrm{super}}\mathcal B}$. Thus, the short exact sequence splits. For the last part, denote $h_q:=\sigma_{-tq(t)^2,q(t),1}\in N$, for all $q\in \Bbbk[t]$. Then, $N\simeq (\Bbbk[t],+)$, via $q(t)\mapsto h_q$. The action of $\Aut_{\mathrm{super}}\mathcal B$ on $N$ is induced by conjugation through the section $\Psi$. Explicitly,
		\[
		\sigma_{r,\mu}\cdot h_q :=\Psi(\sigma_{r,\mu})h_q\Psi(\sigma_{r,\mu})^{-1}=\sigma_{r,0,\mu}\,h_q\,\sigma_{r,0,\mu}^{-1}=\sigma_{-t\mu^2q(\mu^2t)^2,\ \mu q(\mu^2t),\ 1}=h_{\mu q(\mu^2t)}, 
		\]
		for all $\sigma_{r,\mu}\in \Aut_{\mathrm{super}}\mathcal B$. It is straightforward to verify that 
		\[
		\Theta: N\rtimes \Aut_{\mathrm{super}}\mathcal B \longrightarrow \mathbb L,\quad (h_q,\sigma_{r,\mu})\mapsto h_q\,\sigma_{r,0,\mu},
		\]
		is an injective group homomorphism. Also, given $\sigma_{p,q,\mu}\in\mathbb L$, we set $r(t):=p(t)+\mu^{-1}tq(t)^2$ and note that $\sigma_{p,q,\mu}=
		h_{\mu^{-1}q(t)}\,\sigma_{r,0,\mu}$. Thus, $\Theta$ is a group isomorphism.
	\end{proof}
	
	\begin{remark}
		The Gelfand-Kirillov dimension of the centralizer $C_{\mathcal B}(x)$ is $2$. More precisely, its Hilbert series is $H_{C_{\mathcal B}(x)}(t)
		=
		\frac{1+t}{(1-t^2)^2}$.
		Indeed, by Lemma \ref{lemma-centralizador-x1}, $C_{\mathcal B}(x)=\mathcal J\oplus x\mathcal J$ as graded $\Bbbk$-vector spaces, where $\deg u=\deg w=2$ and $\deg x=1$, with $u=y^2-xy$. Since $\mathcal J$ is isomorphic to the Jordan plane, with both generators in degree $2$, we have $H_{\mathcal J}(t)
		=
		\frac{1}{(1-t^2)^2}$. Therefore,
		$
		H_{C_{\mathcal B}(x)}(t)
		=
		H_{\mathcal J}(t)+tH_{\mathcal J}(t)
		=
		\frac{1+t}{(1-t^2)^2}.
		$
		Since $\frac{1}{(1-t^2)^2}
		=
		\sum_{m\geq0}(m+1)t^{2m}$,
		we obtain $\dim_{\Bbbk} C_{\mathcal B}(x)_{2m}=\dim_{\Bbbk} C_{\mathcal B}(x)_{2m+1}=m+1$. Let
		$
		F_n
		:=
		\oplus_{j=0}^n C_{\mathcal B}(x)_j$
		be the filtration induced by the grading. Then, $\dim_{\Bbbk}F_{2m}=(m+1)^2$ and $\dim_{\Bbbk}F_{2m+1}=(m+1)(m+2)$. Thus, in either case, $\dim_{\Bbbk}F_n\sim \frac{1}{4}n^2$. Therefore,
		\[
		\GK C_{\B}(x)=\limsup_{n\to\infty}\frac{\log \dim_{\Bbbk}F_n}{\log n}
		=2.
		\]
	\end{remark}

\subsection{Exponentials of locally nilpotent derivations} In this subsection we explain how the automorphisms $\sigma_{p,q,1}$ arise as exponentials of explicit locally nilpotent derivations of $\mathcal B$. The case $q=0$ recovers the super-automorphisms described in \cite[Proposition 6.2]{AD}, namely $\sigma_{p,0,1}=\exp(\fe_p)$. The results below extend this construction to the whole family $\sigma_{p,q,1}$ by adding a suitable inner locally nilpotent derivation.
	
	Let $h(t)\in\Bbbk[t]$. We define a linear map $\fe_h:\Span\{x,y\}\to\mathcal B$ on the generators of $\mathcal B$ by $	\fe_h(x)=0$ and $\fe_h(y)=xh(z)$.
	
	\begin{lemma} \label{lemaeh}
		The map $\fe_h$ extends uniquely to a locally nilpotent derivation of $\mathcal B$. Moreover, $\exp(\fe_h)=\sigma_{h,0,1}$.
	\end{lemma}
	
	\begin{proof}
		It is straightforward to prove that $\fe_h$ is a well-defined derivation of $\mathcal B$. Moreover,
		$
		\fe_h(z)
		=
		\fe_h(y)x+y\fe_h(x)
		=
		xh(z)x
		=
		0.
		$		
		Thus, $\fe_h$ vanishes on both $x$ and $z$. It follows that, on a PBW monomial, repeated applications of $\fe_h$ only act on the final powers of $y$. Then, for every $a\in\B$, there exists $n\geq1$ such that
		$\fe_h^n(a)=0$. Therefore, $\fe_h$ is locally nilpotent.
		Finally, notice that $\fe_h^2(x)=0$ and $\fe_h^2(y)=0$. Thus,
		$\exp(\fe_h)(x)=x$ and $\exp(\fe_h)(y)=y+xh(z)=\sigma_{h,0,1}(y)$.
	\end{proof}
	
	\begin{remark}
		The fact that $\fe_h^2$ vanishes on the generators $x$ and $y$ does not imply that $\fe_h^2=0$ as an operator on $\mathcal B$. For example, when $h=1$, we obtain $\fe_1^2(y^3)=2xyx$. Thus, $\fe_h$ is locally nilpotent, but $\fe_h^2=0$ does not hold in general.
	\end{remark}

	\begin{lemma}\label{expuh}
		Let $h(t)\in\Bbbk[t]$ and set $u_h:=xh(z)$. Then the inner derivation $\ad(u_h)$ is nilpotent, with $\ad(u_h)^3=0$. Furthermore, its exponential is given by $	\exp(\ad(u_h))=\sigma_{-th(t)^2,h,1}$.
	\end{lemma}
	
	\begin{proof}
		
		Since $u_h^2=0$, for every $a\in\mathcal B$ we have $\ad(u_h)^2(a)=-2u_hau_h$. Thus, we have $\ad(u_h)^3(a)=0$. Also, $\ad(u_h)(x)=0$ and
		$
		\ad(u_h)(y)
		=
		xh(z)y-h(z)z.
		$
		Notice that $\ad(u_h)^2(y)=-2xz\,h(z)^2$. Since $\ad(u_h)^3=0$, we obtain
		$
		\exp(\ad(u_h))
		=
		\id+\ad(u_h)+\frac{1}{2}\ad(u_h)^2.
		$
		Thus, $\exp(\ad(u_h))(x)=x$ and
		$
		\exp(\ad(u_h))(y)=y+xh(z)y-h(z)z
		-xzh(z)^2.
		$
		This is exactly the image of $y$ under $\sigma_{-th(t)^2,h,1}$. Then,
		$
		\exp(\ad(u_h))=\sigma_{-th(t)^2,h,1}.
		$
	\end{proof}
	
	\begin{remark}
		The first assertion of Lemma~\ref{expuh} is a consequence of
		the following simple fact. If $A$ is an associative algebra
		and $w\in A$ satisfies $w^N=0$ for some integer $N\geq1$,
		then the inner derivation $\ad(w)$ is nilpotent of index
		at most $2N-1$. Indeed, for every $a\in A$ and every
		$m\geq0$, we have
		$
		\ad(w)^m(a)
		=
		\sum_{i=0}^{m}(-1)^i\binom{m}{i}w^{m-i}aw^i.
		$
		For $m=2N-1$, each summand contains either $w^i$ with
		$i\geq N$ or $w^{m-i}$ with $m-i\geq N$.
		Then, every summand is zero, and therefore
		$\ad(w)^{2N-1}=0$.
	\end{remark}

	\begin{proposition} \label{sigmapq1}
		Let $p(t),q(t)\in\Bbbk[t]$. Denote $\mathcal D_{p,q}:=\fe_{p+tq^2}+\ad(xq(z))$.
		Then, $\mathcal D_{p,q}$ is locally nilpotent and $\exp(\mathcal D_{p,q})=\sigma_{p,q,1}$.
		
	\end{proposition}
	
	\begin{proof}
		
		We denote $E=\fe_{p+tq^2}$ and $F=\ad(xq(z))$. 	By Lemma \ref{lemaeh}, the derivation $E$ is locally nilpotent. By Lemma \ref{expuh}, the derivation $F$ is nilpotent, hence locally nilpotent. We claim that $E$ and $F$ commute. Indeed,  for every derivation $E$ of $\B$ and every $u\in\B$,
		$
		[E,\ad(u)]=\ad(E(u)).
		$
		Taking $u=xq(z)$ and using $E(x)=E(z)=0$, we obtain
		$[E,F]=0$. Hence, $\mathcal D_{p,q}=E+F$ is locally nilpotent and
		$
		\exp(\mathcal D_{p,q})=\exp(F)\circ\exp(E).
		$
		By Lemma \ref{lemaeh}, $\exp(E)=\sigma_{p+tq^2,0,1}$ and, by Lemma \ref{expuh}, $\exp(F)=\sigma_{-tq^2,q,1}$. Using the composition from Proposition \ref{composition}, we obtain
		$
		\sigma_{-tq^2,q,1}\circ\sigma_{p+tq^2,0,1}
		=
		\sigma_{p,q,1}.
		$
		Therefore, $\exp(\mathcal D_{p,q})=\sigma_{p,q,1}$.
	\end{proof}
		
	\section{The automorphism group of	\texorpdfstring{$\B$}{B}}
	In this section, we prove that every automorphism
	$\varphi$ of $\B$ admits a decomposition
	$
	\varphi=\Ad_U\circ \, \sigma_{q,\mu},
	$
	where $U\in\B^\times$, $q\in\Bbbk[t]$, $\mu\in\Bbbk^\times$ and $\sigma_{q,\mu}$ is a super-automorphism. More precisely, we establish the
	semidirect product decomposition
	\[
	\Aut\B
	=
	\operatorname{Inn}\B
	\rtimes
	\Aut_{\mathrm{super}}\B.
	\]
	By \cite[Lemma~5.2 and Theorem~5.3]{AD},
	$\sigma_{q,\mu}$ are precisely the automorphisms
	preserving the $\mathbb Z_2$-grading.
	Let $\pi:\mathcal B\rightarrow\mathcal B/(x)\simeq
	\Bbbk[y]$ be the canonical projection and fix $I:=\ker\pi=(x)$ the two-sided ideal of $\B$ generated by $x$.
	
	\begin{lemma}\label{lemanilpotentkernel}
		The following assertions hold:
		\begin{enumerate}[\rm (i)]
			\setlength{\itemsep}{0.3em}
			
			\item Every nilpotent element of $\mathcal B$ belongs to $I$;
			
			\item For every $n\geq1$, $I^n=w^{n-1}I$;
			
			\item $\bigcap_{n\geq1}I^n=0$.
			
		\end{enumerate}
	\end{lemma}
	
	\begin{proof}
		Let $a\in\mathcal B$ be nilpotent. Thus, $\pi(a)$ is nilpotent in $\mathcal B/(x)\simeq\Bbbk[y]$. Since $\Bbbk[y]$ has no nonzero nilpotent elements, we obtain $\pi(a)=0$. Then, $a\in I$, proving (i). For (ii), we notice that $I$, as vector subspace of $\B$, is given by
		\begin{equation}\label{ideal-I}
			I
			=
			\Span
			\{z^by^c:b\geq1,\ c\geq0\}
			\oplus
			\Span
			\{xz^by^c:b,c\geq0\}.
		\end{equation}
		Since $x,z\in I$, this description gives $I=x\mathcal B+z\mathcal B$. Using $w=z+xy$, we obtain $I=x\mathcal B+w\mathcal B$. Moreover, $x^2=0$ implies $xw=wx=xyx$. Then,
		\[
		xI
		=x(x\mathcal B+w\mathcal B)
		=xw\mathcal B
		=wx\mathcal B
		\subseteq wI.
		\]
		Since $I$ is a two-sided ideal, $\mathcal BI=I$, and therefore $I^2=(x\mathcal B+w\mathcal B)I=xI+wI=wI$. The case $n=1$ is immediate, and the case $n=2$ was proved above. We proceed by induction. If $I^n=w^{n-1}I$, then $I^{n+1}=I^nI=w^{n-1}I^2=w^nI$. Thus, $I^n=w^{n-1}I$ for every $n\geq1$, proving (ii). 
				Finally, $I$ is generated by the homogeneous element $x$ of degree $1$, so $I\subseteq\oplus_{j\geq1}\mathcal B_j$. Since $\deg w=2$, assertion (ii) gives
		$
		I^n=w^{n-1}I
		\subseteq
		\bigoplus_{j\geq2n-1}\mathcal B_j.
		$
		Suppose that $0\neq a\in\bigcap_{n\geq1}I^n$, and let $d$ be
		the largest degree of a nonzero homogeneous component of $a$.
		Choose $n$ such that $2n-1>d$. Since $a\in I^n$, the preceding
		inclusion implies that every homogeneous component of $a$
		of degree less than $2n-1$ is zero, contradicting the choice
		of $d$, which proves (iii).
		\end{proof}
	
	\begin{corollary}\label{cor:nilpotent}
		The two-sided ideal generated by all nilpotent elements of
		$\mathcal B$ is $I=(x)$. In particular, $\varphi(I)=I$ for every $\varphi\in\operatorname{Aut}\mathcal B$.
	\end{corollary}
	
	\begin{proof}
		Let $N$ denote the two-sided ideal generated by all nilpotent
		elements of $\mathcal B$. By Lemma \ref{lemanilpotentkernel}(i), every nilpotent element belongs to $I$, so $N\subseteq I$.
		Conversely, since $x^2=0$, we have $x\in N$. As $N$ is a
		two-sided ideal, it contains $I=(x)$. Thus, $N=I$. Every automorphism $\varphi$ of $\mathcal B$ maps the set of
		nilpotent elements onto itself and therefore satisfies
		$\varphi(N)=N$. Since $N=I$, we obtain $\varphi(I)=I$.
	\end{proof}
	
	\subsection{The induced actions on \texorpdfstring{$\B/I$ and $I/I^2$}{B/I and I/I squared}}
	
	By Corollary~\ref{cor:nilpotent}, every automorphism
	$\varphi$ of $\B$ preserves $I=(x)$ and hence also $I^2$.
	Thus, $\varphi$ induces an automorphism of
	$\B/I\simeq\Bbbk[y]$ and a bijection of $I/I^2$. We record the descriptions of $I$ and $I^2$ that will be
	used below. The PBW monomials with $\epsilon+a\geq1$ form a basis
	of $I$, since the remaining monomials are the powers of $y$,
	whose classes form a basis of $\B/I\simeq\kk[y]$.
	Likewise, the PBW monomials with $\epsilon+a\geq2$ form a basis
	of $I^2$. Indeed, the reductions
	$x^2=zx=0$, $yx=z$, and $yz=xy^2+xz$ never decrease the
	number of occurrences of $x$ and $z$ in a word. Products of
	two elements of $I$ therefore involve only monomials with
	$\epsilon+a\geq2$. Conversely, each such monomial belongs to
	$I^2$, because $x,z\in I$.
	Consequently, $I/I^2$ is a free right $\kk[y]$-module with basis
	$\{\overline x,\overline z\}$, where
	$\overline x=x+I^2$ and $\overline z=z+I^2$.
	The left action of $y$ is given by
	$y\overline x=\overline z$ and
	$y\overline z=\overline x y^2$. It follows that left and right
	multiplication by $y^2$ agree on $I/I^2$.

	Let $\varphi\in\Aut\B$. Since $\varphi(I)=I$, the maps
	$\overline\varphi$ on $\B/I$ and $\Phi$ on $I/I^2$
	induced by $\varphi$ are bijective. The latter is semilinear with
	respect to $\overline\varphi$: for $v\in I/I^2$ and
	$f\in\kk[y]$, one has
	$\Phi(vf)=\Phi(v)\overline\varphi(f)$.
	In particular, $\Phi$ sends a right $\kk[y]$-basis to a right
	$\kk[y]$-basis. Indeed, generation follows from the surjectivity
	of both maps, and a relation between the image vectors becomes a
	relation between the original vectors upon applying $\Phi^{-1}$
	and $\overline\varphi^{-1}$ to its coefficients.
	
	\begin{lemma}\label{lemma-no-translation-quotient}
		Let $\varphi\in\Aut\B$. Then, there exists
		$\mu\in\kk^\times$ such that $\varphi(y)-\mu y\in I$.
	\end{lemma}
	
	\begin{proof}
		The automorphism induced on $\B/I\simeq\kk[y]$ has the form
		$\overline\varphi(y)=\lambda+\mu y$, with
		$\lambda\in\kk$ and $\mu\in\kk^\times$.
		Put $a=\varphi(x)$, $b=\varphi(y)$, and
		$\overline a=a+I^2$. Applying $\varphi$ to the cubic defining
		relation gives $b^2a-ab^2-aba=0$. Since $aba\in I^2$, its image
		in $I/I^2$ is
		\[
		0=(\lambda+\mu y)^2\overline a-
		\overline a(\lambda+\mu y)^2
		=2\lambda\mu(y\overline a-\overline a y),
		\]
		where the last equality uses the coincidence of the two actions of
		$y^2$ on $I/I^2$.
			Suppose that $\lambda\neq0$. Then
		$y\overline a=\overline a y$. On the other hand,
		$\varphi(z)+I^2=(\lambda+\mu y)\overline a
		=\overline a(\lambda+\mu y)$, since $z=yx$.
		Thus $\varphi(z)+I^2$ belongs to the right submodule generated by
		$\overline a$. This contradicts the fact that
		$\{\varphi(x)+I^2,\varphi(z)+I^2\}$ is a right
		$\kk[y]$-basis of $I/I^2$. Hence, $\lambda=0$.
	\end{proof}
	
	\begin{lemma}\label{lemma-action-I-I2}
		Let $\varphi\in\Aut\B$ and suppose that
		$\varphi(y)-\mu y\in I$, with $\mu\in\kk^\times$.
		Then $\varphi(x)-\alpha x\in I^2$ for some
		$\alpha\in\kk^\times$.
	\end{lemma}
	
	\begin{proof}
		Write $\varphi(x)+I^2=\overline x f(y)+\overline z g(y)$, with
		$f,g\in\kk[t]$. Since $\varphi(z)=\varphi(y)\varphi(x)$,
		the formulas for the left action of $y$ give
		$\varphi(z)+I^2=\mu\overline z f(y)+\mu\overline x y^2g(y)$.
		The columns of
		\[
		\begin{pmatrix}
			f(y)&\mu y^2g(y)\\
			g(y)&\mu f(y)
		\end{pmatrix}
		\]
		are therefore the coordinates of the basis
		$\{\varphi(x)+I^2,\varphi(z)+I^2\}$ relative to
		$\{\overline x,\overline z\}$. This matrix is invertible over
		$\kk[y]$, so $\mu(f(y)^2-y^2g(y)^2)\in\kk^\times$.
		It follows that both $f(y)-yg(y)$ and $f(y)+yg(y)$ are nonzero
		constants. Their difference is divisible by $y$, and hence is zero. Thus, $g=0$ and $f=\alpha\in\kk^\times$, as required.
	\end{proof}
	
	\begin{lemma}\label{lemma-alpha-equals-mu}
		Let $\varphi\in\Aut\B$. If
		$\varphi(y)-\mu y\in I$ and $\varphi(x)-\alpha x\in I^2$,
		where $\alpha,\mu\in\kk^\times$, then $\alpha=\mu$.
	\end{lemma}
	
	\begin{proof}
		Let $a=\varphi(x)$ and $b=\varphi(y)$. By hypothesis,
		we can write $a=\alpha x+r$ and $b=\mu y+s$, where $r\in I^2$ and $s\in I$. The ideals $I$ and $I^2$ are homogeneous, and
		$I\cap\B_{\leq1}=\kk x$ and $I^2\cap\B_{\leq2}=0$. Thus, every homogeneous component of $r$ has degree
		at least three, while $s$ has zero constant term
		and degree-one component $\gamma x$ for some
		$\gamma\in\kk$. Consequently, $a$ and $b$ have
		zero constant term and degree-one components
		$\alpha x$ and $\mu y+\gamma x$, respectively. Applying $\varphi$ to the cubic defining relation gives $b^2a-ab^2-aba=0$. Since $a$ and $b$ have zero constant term, the
		degree-three components of these products depend
		only on their degree-one components. Using $x^2=0$, we obtain
		$(\mu y+\gamma x)^2=\mu^2y^2+\mu\gamma(yx+xy)$. Moreover,
		$(yx+xy)x=x(yx+xy)=xyx$, so the terms involving $\gamma$ cancel when we compare the degree-three components. Therefore,
		\[
		\begin{aligned}
			0
			&=\alpha\mu^2(y^2x-xy^2)-\alpha^2\mu xyx=\alpha\mu(\mu-\alpha)xz.
		\end{aligned}
		\]
		Since $xz$ is a nonzero PBW monomial and
		$\alpha,\mu\neq0$, we conclude that $\alpha=\mu$.
	\end{proof}
	
	\begin{proposition}\label{prop-quotient}
		For every $\varphi\in\Aut\B$, there is a unique
		$\mu\in\kk^\times$ such that
		$
		\varphi(y)-\mu y\in I$ and $\varphi(x)-\mu x\in I^2$.
	\end{proposition}
	
	\begin{proof}
		Existence follows from Lemmas \ref{lemma-no-translation-quotient},
		\ref{lemma-action-I-I2} and \ref{lemma-alpha-equals-mu}.
		Uniqueness follows from $y+I\neq0$ in $\B/I$.
	\end{proof}
	
\subsection{A matrix model for the super Jordan plane}
	
	We use the presentation
	\[
	\mathbf J
	\simeq
	\kk\langle w,s\rangle/(sw-ws-w^2),
	\]
	where $s$ corresponds to $y^2$ in
	$\mathbf J=\kk\langle y^2,w\rangle$.
	Let $T=	\Mr_2(\mathbf J)$ and write  $e_{ij}$ for the matrix with $1$ in the $(i,j)$-entry
	and $0$ elsewhere. When working in $T$, we identify $a\in\mathbf J$ with the multiple $a\Id_2$
	of the $2\times2$ identity matrix. 
		Recall that $\mathbf J$ is the Ore extension
	$\kk[w][s;\delta]$, where $\delta(w)=w^2$.
	In particular, it is a domain with basis
	$\{w^a s^b:a,b\geq0\}$.
	Let $\pi:T\to 
	\Mr_2(\kk[s])$ be entrywise reduction
	modulo $w\mathbf J$.
	
	\begin{lemma}\label{lem:matrix-model}
		The map on generators of $\B$ given by
		\[
		x\longmapsto X=\begin{pmatrix}0&w\\0&0\end{pmatrix},
		\quad
		y\longmapsto Y=\begin{pmatrix}0&s\\1&0\end{pmatrix}
		\]
		defines an embedding $\iota:\B\hookrightarrow T$, with
		$\iota(w)=w\Id_2$ and $\iota(y^2)=s\Id_2$.
		Moreover, $\iota(I)=wT$ and
		$\iota(\B)=\pi^{-1}(\kk[\overline Y])$, where
		$\overline Y=\pi(Y)$.
	\end{lemma}
	
	\begin{proof}
		The equalities $X^2=0$, $Y^2=s\Id_2$, and
		$Y^2X-XY^2=w^2e_{12}=XYX$ show that the assignments define a
		homomorphism. We also have $YX=we_{22}$ and
		$XY+YX=w\Id_2$. Moreover, we have that
		\begin{equation}\label{eq:pbw-matrices}
			\begin{aligned}
				\iota(y^{2r})&=s^r\Id_2,&
				\iota(y^{2r+1})&=s^{r+1}e_{12}+s^re_{21},\\
				\iota(z^ay^{2r})&=w^as^re_{22},&
				\iota(z^ay^{2r+1})&=w^as^re_{21},&& (a\geq1),\\
				\iota(xz^ay^{2r})&=w^{a+1}s^re_{12},&
				\iota(xz^ay^{2r+1})&=w^{a+1}s^re_{11},&& (a\geq0),
			\end{aligned}
		\end{equation}
		where $r\geq0$ in every case. To prove injectivity, let
		$b\in\ker\iota$. Using the PBW basis of $\B$, write
		$b=b_1+\cdots+b_6$, where
		\[
		\begin{array}{@{}lll@{}}
			b_1=\sum_{r\geq0}\alpha_r y^{2r},
			&
			b_3=\sum_{a\geq1,\,r\geq0}\gamma_{a,r}z^ay^{2r},
			&
			b_5=\sum_{a,r\geq0}\eta_{a,r}xz^ay^{2r},
			\\[.4em]
			b_2=\sum_{r\geq0}\beta_r y^{2r+1},
			&
			b_4=\sum_{a\geq1,\,r\geq0}\delta_{a,r}z^ay^{2r+1},
			&
			b_6=\sum_{a,r\geq0}\vartheta_{a,r}xz^ay^{2r+1}.
		\end{array}
		\]
		All sums are finite. Since $\iota(b)=0$, Equations \eqref{eq:pbw-matrices} gives
				\[
		\begin{aligned}
			0=\iota(b)={}&
			\sum_{r\geq0}
			\bigl(\alpha_rs^r\Id_2
			+\beta_r(s^{r+1}e_{12}+s^re_{21})\bigr)\\
			&+\sum_{a\geq1,\,r\geq0}
			w^as^r(\gamma_{a,r}e_{22}+\delta_{a,r}e_{21})
			+\sum_{a,r\geq0}
			w^{a+1}s^r(\eta_{a,r}e_{12}+\vartheta_{a,r}e_{11}).
		\end{aligned}
		\]
		Applying $\pi$ annihilates the last two sums, since all their
		terms belong to $wT$. The $(2,2)$- and $(2,1)$-entries of the
		resulting identity give, respectively,
		$\sum_{r\geq0}\alpha_rs^r=0$ and
		$\sum_{r\geq0}\beta_rs^r=0$ in $\Bbbk[s]$.
		Hence, $\alpha_r=\beta_r=0$ for every $r$.
		Returning to the identity in $T$, we obtain
		\begin{equation}\label{eq-linear-combination}
			0=
			\sum_{a\geq1,\,r\geq0}
			w^as^r(\gamma_{a,r}e_{22}+\delta_{a,r}e_{21})
			+
			\sum_{a,r\geq0}
			w^{a+1}s^r(\eta_{a,r}e_{12}+\vartheta_{a,r}e_{11}).
		\end{equation}
Considering the four matrix entries separately, we obtain
four vanishing linear combinations of distinct monomials
in the PBW basis $\{w^ms^r:m,r\geq0\}$ of $\J$.
Therefore, all coefficients $\gamma_{a,r}$, $\delta_{a,r}$,
$\eta_{a,r}$, and $\vartheta_{a,r}$ vanish.
It follows that $b=0$, proving that $\iota$ is injective.
		Among the PBW monomials $\{x^\epsilon z^a y^b:
		\epsilon\in\{0,1\},\ a,b\geq0\}$ of $\mathcal B$, those with $\epsilon+a\geq1$ form a basis of
		$I=(x)$. Indeed, these monomials belong to $(x)$, since either
		$\epsilon=1$ or $a\geq1$ and $z=yx$. The remaining PBW
		monomials are precisely the powers $y^b$, whose residue classes form
		the standard basis of $\mathcal B/(x)\simeq\kk[y]$. It follows that the monomials with $\epsilon+a\geq1$ span exactly
		the kernel $I$ of the quotient map. The families in the sum \eqref{eq-linear-combination} are
		the matrix images of these monomials. Together they are precisely
		$
		\{w^ms^re_{ij}:m\geq1,\ r\geq0,\ 1\leq i,j\leq2\},
		$
		which is a basis of $wT$. Hence,  $\iota(I)=wT$. Finally, the restriction of $\pi$ to $\iota(\mathcal B)$ has kernel
		$\iota(I)=wT$ and $\pi(X)=0$, $\pi(Y)=\overline Y$. Thus, $\pi(\iota(\mathcal B))=\kk[\overline Y]$. Since $\ker\pi=wT\subseteq\iota(\B)$, it follows that
		$\iota(\B)=\pi^{-1}\bigl(\kk[\overline Y]\bigr)$.
	\end{proof}

We henceforth identify $\B$ with its image in $T$.
Since $\J$ is a domain and
$sw=w(s+w)$, $ws=(s-w)w$, the element $w$ is regular
and normal in $\J$, and hence in $T$.
It is also regular and normal in $\B$, as we prove directly
in Lemma~\ref{lem-x21-degree}. Its normality is noted in
\cite[\S2.3]{Lu2026SJP}. Consequently, the Ore localization $\B[w^{-1}]$ exists.
The inclusion $\B\subseteq T$ induces an identification
$\B[w^{-1}]=T[w^{-1}]$, since every $t\in T$ satisfies
$wt\in I\subseteq\B$ and thus $t=w^{-1}(wt)\in\B[w^{-1}]$.
Finally, the existence of the classical ring of quotients
$Q(\B)$ follows from \cite[Proposition~3.1(3)]{Lu2026SJP},
and we obtain the inclusions
$\B\subseteq T\subseteq\B[w^{-1}]\subseteq Q(\B)$.
	
\begin{lemma}\label{lem:idealizer}
	Inside $Q(\B)$, we have
	$
	T=\{q\in Q(\B):qI\subseteq I\}.
	$
	Every automorphism $\varphi$ of $\B$ extends uniquely to an
	automorphism $\widehat\varphi$ of $T$ satisfying
	$\widehat\varphi(wT)=wT$.
\end{lemma}

	\begin{proof}
		Since $I=wT$, the condition $qI\subseteq I$ is equivalent to
		$(w^{-1}qw)T\subseteq T$. Hence, $qI\subseteq I$ if and only if
		$w^{-1}qw\in T$. Clearly $wtw^{-1}\in T$ and $w^{-1}tw\in T$, for all $t\in T$. Thus, $qI\subseteq I$ if and only if $q\in T$.
		Every automorphism $\varphi$ of $\B$ extends uniquely to an
		automorphism $\widehat\varphi$ of $Q(\B)$, since it preserves
		the regular elements. As $I$ is characteristic and
		$\widehat\varphi|_{\B}=\varphi$, we have
		$\widehat\varphi(I)=I$. The displayed characterization of $T$
		therefore implies that $\widehat\varphi(T)=T$.
		Its restriction to $T$ preserves $wT$, because $wT=I$. Finally, if $\psi\in\Aut T$ also extends $\varphi$,
		then $wt\in I\subseteq\B$ gives $\varphi(w)\psi(t)=\psi(wt)=\varphi(wt)$, $t\in T$. Since $\varphi(w)$ is invertible in $Q(\B)$, this
		determines $\psi(t)$ uniquely for every $t\in T$.
	\end{proof}

\subsection{Descent for projective modules over the Jordan plane}
	
	Let $A$ denote the localization $\J[w^{-1}]$. Every element of $A$ has a unique finite
	expression $\sum_{i\in\mathbb Z}w^if_i(s)$, with
	$f_i\in\kk[s]$. For a nonzero element $a$, define $v(a)$
	to be the least index for which $f_i\neq0$.
	The identity $sw^j=w^js+jw^{j+1}$, valid for every integer $j$,
	implies
	\begin{equation}\label{eq:valuation}
		v(ab)=v(a)+v(b)\quad(a,b\neq0),
		\quad
		\J=\{0\}\cup\{a\in A:v(a)\geq0\}.
	\end{equation}
	For the first equality, commuting a polynomial in $s$ past
	$w^j$ leaves its initial $w$-exponent unchanged and produces
	only terms of larger $w$-exponent. Thus, if the initial terms
	of two factors are $w^if(s)$ and $w^jg(s)$, the coefficient
	of $w^{i+j}$ in their product is $f(s)g(s)\neq0$.
	The second assertion follows from the PBW basis of $\J$.
	
	\begin{lemma}\label{lem:descent}
		Let $P$ be a finitely generated projective right $\J$-module.
		If $P\otimes_{\J}A\simeq A$ as right $A$-modules, then
		$P\simeq\J$ as right $\J$-modules.
	\end{lemma}
	
	\begin{proof}
		Since $P$ is finitely generated and projective, there exist a
		right $\J$-module $Q$ and an integer $n\geq1$ such that
		$P\oplus Q\simeq\J^n$. Tensoring with $A$ gives
		\[
		\J^n\otimes_{\J}A
		\simeq (P\otimes_{\J}A)\oplus(Q\otimes_{\J}A).
		\]
		Under the canonical identification
		$\J^n\otimes_{\J}A\simeq A^n$, the natural map
		$\jmath_n:\J^n\to\J^n\otimes_{\J}A$, defined by
		$\jmath_n(v)=v\otimes1$, corresponds to the coordinatewise
		inclusion $\J^n\hookrightarrow A^n$. Hence $\jmath_n$
		is injective. Under the chosen decomposition $\J^n\simeq P\oplus Q$
		and the induced decomposition of its tensor product,
		$\jmath_n$ corresponds to the direct sum of the natural maps
		$P\to P\otimes_{\J}A$ and $Q\to Q\otimes_{\J}A$.
		Consequently, the map $P\to P\otimes_{\J}A$,
		given by $p\mapsto p\otimes1$, is injective. Choose a right
		$A$-module isomorphism $\theta:P\otimes_{\J}A\to A$. The map
		$p\mapsto\theta(p\otimes1)$ identifies $P$ with a finitely
		generated right $\J$-submodule of $A$. Under this
		identification, $PA=A$, where $PA$ denotes the right
		$A$-submodule generated by $P$. Indeed, every element of
		$P\otimes_{\J}A$ is a finite sum of tensors $p_i\otimes a_i$
		and $\theta(p_i\otimes a_i)=\theta(p_i\otimes1)a_i$. Hence
		the assertion follows from the surjectivity of $\theta$.
		
		Let $p_1,\ldots,p_r$ be nonzero generators of the right
		$\J$-module $P$. For each $1\leq i\leq r$, we write
		$p_i=\sum_{j=m_i}^{M_i}w^jf_{ij}(s)$, with
		$f_{ij}(s)\in\kk[s]$. Choose an integer $N\geq0$ such that
		$N+m_1,\ldots,N+m_r$ are all nonnegative. Then
		$w^Np_i=\sum_{j=m_i}^{M_i}w^{N+j}f_{ij}(s)\in\J$,
		because all the exponents of $w$ are nonnegative.
		Set $P'=w^NP$. Every $p\in P$ can be written as
		$p=\sum_i p_ic_i$, with $c_i\in\J$, and consequently
		$w^Np=\sum_i(w^Np_i)c_i$. It follows that
		$P'=\sum_i(w^Np_i)\J\subseteq\J$. Thus $P'$ is a
		finitely generated right ideal of $\J$. The map $F:P\to P'$,
		defined by $F(p)=w^Np$, is a homomorphism of right
		$\J$-modules, since $F(pc)=w^N(pc)=(w^Np)c=F(p)c$
		for $p\in P$ and $c\in\J$. It is surjective by the
		definition of $P'$. Moreover, $w$ is invertible in $A$
		and left multiplication by $w^{-N}$ maps $P'$ back into
		$P$ and is inverse to $F$. Consequently, $P'\simeq P$
		as right $\J$-modules, so $P'$ is projective. Finally,
		$P'A=(w^NP)A=w^N(PA)=w^NA=A$, where the last equality
		follows from the invertibility of $w^N$ in $A$.
		In particular, $P'\neq0$. We may therefore replace $P$
		by the isomorphic module $P'$, retaining the notation $P$.
		Henceforth, $P$ is a nonzero finitely generated projective
		right ideal of $\J$ satisfying $PA=A$.
		Since $P$ is a nonzero right ideal of $\J$, the set
		$\{v(p):0\neq p\in P\}$ is a nonempty subset of
		$\mathbb Z_{\geq0}$. Let $r$ be its minimum, and choose
		$p_*\in P\setminus\{0\}$ such that $v(p_*)=r$.
		Set $P'=w^{-r}P$. For every nonzero $p\in P$, we have by
		\eqref{eq:valuation} that
		$v(w^{-r}p)=-r+v(p)\geq0$. Hence $w^{-r}p\in\J$.
		Consequently, $P'$ is a right ideal of $\J$.
		Left multiplication by $w^{-r}$ defines a right
		$\J$-module isomorphism $P\to P'$, with inverse given
		by left multiplication by $w^r$. Thus $P'$ remains
		finitely generated and projective. Moreover,
		$P'A=w^{-r}(PA)=w^{-r}A=A$, since $w$ is invertible
		in $A$. The element $p_0=w^{-r}p_*$ is nonzero and
		satisfies $v(p_0)=0$. Replacing $P$ by $P'$ and retaining
		the notation $P$, we may therefore assume that
		$P\subseteq\J$, $PA=A$, and $v(p_0)=0$ for some
		$p_0\in P\setminus\{0\}$.
		
		Let $f:P\to\J$ be a homomorphism of right $\J$-modules.
		Consider the multiplication map
		$m_P:P\otimes_{\J}A\to A$, defined by
		$m_P(p\otimes a)=pa$. This map is surjective because
		$PA=A$. It is also injective. Indeed, since $w$ is normal
		in $\J$, any finite collection of elements of $A$ admits
		a common right denominator $w^\ell$. Given an element
		$\sum_i p_i\otimes a_i$ of $P\otimes_{\J}A$, write
		$a_i=b_iw^{-\ell}$ with $b_i\in\J$. Then
		\[
		\sum_i p_i\otimes a_i
		=
		\big(\sum_i p_ib_i\big)\otimes w^{-\ell}
		=
		p\otimes w^{-\ell},
		\quad\quad p:=\sum_i p_ib_i\in P.
		\]
		If $m_P(p\otimes w^{-\ell})=pw^{-\ell}=0$, then $p=0$,
		since $P\subseteq A$ and $w^{-\ell}$ is invertible in $A$.
		Hence $m_P$ is an isomorphism of right $A$-modules.
		Consider the map $f\otimes\operatorname{id}_A$ that sends
		$p\otimes a$ to $f(p)\otimes a$. Let
		$m_{\J}:\J\otimes_{\J}A\to A$ be the canonical isomorphism
		defined by $m_{\J}(j\otimes a)=ja$. Consider
		$\widetilde f:=m_{\J}\circ(f\otimes\operatorname{id}_A)
		\circ m_P^{-1}$. Then $\widetilde f$ is a right $A$-linear
		endomorphism of $A$, and $\widetilde f(p)=f(p)$ for every
		$p\in P$. Let $q=\widetilde f(1)$. Observe that
		$\widetilde f(a)=\widetilde f(1)a=qa$ for every $a\in A$.
		Consequently, $f(p)=qp$ for every $p\in P$.
		We claim that $q\in\J$. This is immediate if $q=0$.
		Otherwise, the element $p_0\in P\setminus\{0\}$ chosen
		above satisfies $v(p_0)=0$ and
		$qp_0=f(p_0)\in\J\setminus\{0\}$, since $A$ is a domain.
		By \eqref{eq:valuation},
		$v(q)=v(q)+v(p_0)=v(qp_0)\geq0$, and hence $q\in\J$.
		Thus every right $\J$-module homomorphism $P\to\J$
		is given by left multiplication by an element of $\J$.
		Choose a finite projective dual basis $p_i\in P$,
		$f_i\in\Hom_{\J}(P,\J)$, so that
		$p=\sum_i p_if_i(p)$ for every $p\in P$.
		Write $f_i(p)=q_ip$, with $q_i\in\J$, as above.
		For a nonzero $p\in P$, the dual-basis identity gives
		$(1-\sum_i p_iq_i)p=0$. Since $\J$ is a domain,
		$\sum_i p_iq_i=1$. Each summand belongs to the right ideal
		$P$, so $1\in P$ and $P=\J$.
		Since each replacement of $P$ was an isomorphism of right
		$\J$-modules, the original module is also isomorphic to $\J$.
	\end{proof}
	
	\begin{remark}\label{rem:jordan-units-centralizer}
		It is well known that the units of  $\J$ are the elements of $\kk^\times$. Also, 
		the centralizer of $w$ in $\J$ is $\kk[w]$. In fact,  write a nonzero element as
		$a=\sum_{j=0}^n a_j(w)s^j$, where $a_j(w)\in\kk[w]$
		and $a_n(w)\neq0$.
		If $n>0$, then the coefficient of $s^{n-1}$ in
		$[a,w]$ is $n a_n(w)w^2\neq 0$. Thus $[a,w]=0$ forces $n=0$.
	\end{remark}
	
	\subsection{Automorphisms of \texorpdfstring{$\Mr_2(\J)$}{M2(J)}}

	Let $\Lambda$ be a unital $\kk$-algebra. A bimodule $M$
	over $\Lambda$ is called \emph{$\kk$-central} if its
	left and right scalar actions agree, that is, $
	(c1_{\Lambda})m=m(c1_{\Lambda})$,
	for all $c\in\kk,\ m\in M$. The bimodule $M$ is called \emph{invertible} if there
	exists a $\Lambda$-bimodule $N$ such that
$	M\otimes_{\Lambda}N\simeq\Lambda$ and $N\otimes_{\Lambda}M\simeq\Lambda$
	as $\Lambda$-bimodules, where $\Lambda$ acts on itself
	by left and right multiplication. For the remainder of this subsection, we assume in addition
	that $\kk$ is algebraically closed.

	\begin{lemma}\label{lem:picard-localization}
		Every invertible $\kk$-central $A$-bimodule is free of
		rank one as a right $A$-module.
	\end{lemma}
	
	\begin{proof} Let $R=\kk[w,w^{-1}]$ and define $\partial=w^{-2}s\in A$. Since $sw-ws=w^2$, we have $[\partial,w]
		=w^{-2}(sw-ws)=1$. The commutator identity $[\partial,ab]=[\partial,a]b+a[\partial,b]$ gives $[\partial,w^n]=nw^{n-1}$ for every $n\in\mathbb Z$.
		For negative powers, one first applies this identity to
		$ww^{-1}=1$, obtaining $[\partial,w^{-1}]=-w^{-2}$, and then proceeds by induction. Consequently, $\partial f(w)=f(w)\partial+f'(w)$, where $f(w)\in R$.	
		Let $\delta=d/dw$. These identities define an algebra homomorphism $R[t;\delta]\longrightarrow A$, by
		$w\mapsto w$ and $t\mapsto w^{-2}s$. Conversely, the map $w\mapsto w$ and $s\mapsto w^2t$ respect the defining relation of $\J$, because $(w^2t)w-w(w^2t)=w^2(tw-wt)=w^2$. Since $w$ is invertible in $R[t;\delta]$, this homomorphism extends from $\J$ to $A$. The two homomorphisms are inverse to each other.
		Hence, $A\simeq R[\partial;d/dw]$. In characteristic zero, the ring $\mathcal D(R)$ of $\kk$-linear differential operators on the Laurent polynomial ring $R$ is precisely this Ore extension;
		see \cite[Example 2.1]{BCKMPV2023}. Therefore, we identify
		$A\simeq\mathcal D(R)$.	Let $M$ be an invertible $\kk$-central $A$-bimodule.
		The ring $R$ is a finitely generated Dedekind
		$\kk$-algebra with fraction field $K=\kk(w)$.
		Since $\kk$ is algebraically closed of characteristic
		zero, \cite[Corollary 1.13]{CH} applies.
		It yields a nonzero ideal $\mathfrak a\subseteq R$
		and an isomorphism $M\simeq\mathcal D(\mathfrak a,R)$ of left $A$-modules, where $\mathcal D(\mathfrak a,R)
		=\{P\in\mathcal D(K):P(\mathfrak a)\subseteq R\}$. By \cite[\S1.7]{CH},
		\[
		\mathcal D(\mathfrak a,R)=A\mathfrak a^*,
		\quad
		\mathfrak a^*
		=
		\{q\in K:q\mathfrak a\subseteq R\}.
		\]
	The ring $R$, being a localization of the principal
		ideal domain $\kk[w]$, is itself a principal ideal
		domain.
		Thus $\mathfrak a=Rh$ for some $0\neq h\in R$, and $\mathfrak a^*=Rh^{-1}$. It follows that $M\simeq Ah^{-1}$ as left $A$-modules. Right multiplication by $h^{-1}$ defines a left
		$A$-module isomorphism from $A$ to $Ah^{-1}$, $a\mapsto ah^{-1}$, whose inverse is right multiplication by $h$.
		These products are taken in $\mathcal D(K)$.
		Consequently, $M$ is free of rank one as a left
		$A$-module.
		
		We now use the invertibility of $M$ to obtain the
		corresponding statement for its right module structure.
		Write $\operatorname{End}_A({}_AM)$ for the ring of
		endomorphisms of $M$ as a left $A$-module, with
		multiplication given by composition.
	For $a\in A$, the map $R_a:M\to M$ given by $R_a(m)=m\cdot a$
	is a morphism of left $A$-modules.
		Since $R_a\circ R_b=R_{ba}$, the map $\Theta:A^{\mathrm{op}}
		\to
		\operatorname{End}_A({}_AM)$ given by $a^{\mathrm{op}}\mapsto R_a$ defines a ring morphism.
		This morphism is an isomorphism.
		Indeed, if $N$ is an inverse bimodule for $M$, then
		$M\otimes_A-$ and $N\otimes_A-$ are mutually inverse
		equivalences on the category of left $A$-modules.
		In particular, $M\otimes_A-$ induces an isomorphism
		$\operatorname{End}_A({}_AA)
		\longrightarrow
		\operatorname{End}_A({}_AM)$, using the canonical identification $M\otimes_A A\simeq M$. Under the identification
		$\operatorname{End}_A({}_AA)\simeq A^{\mathrm{op}}$
		by right multiplication, this is precisely $\Theta$.
				Choose a basis element $e$ of the free left
		$A$-module $M$.
		For each $a\in A$, there is a unique $\gamma(a)\in A$
		such that $e\cdot a=\gamma(a)e$. The right action is additive and unital, so $\gamma$ is additive and $\gamma(1)=1$.
		Moreover,
		\[
		\begin{aligned}
			\gamma(ab)e
			&=e\cdot(ab)=(e\cdot a)\cdot b=(\gamma(a)e)\cdot b=\gamma(a)\gamma(b)e.
		\end{aligned}
		\]
		Thus $\gamma(ab)=\gamma(a)\gamma(b)$.
		Since $M$ is $\kk$-central,
		$\gamma(\lambda)=\lambda$ for every $\lambda\in\kk$.
				Every left $A$-linear endomorphism of $M$ is uniquely
		determined by its value at $e$, and this value may be
		any element $ce$, with $c\in A$.
		Since $R_a(e)=\gamma(a)e$, the bijectivity of
		$\Theta$ implies the bijectivity of $\gamma$.
		Hence $\gamma$ is a $\kk$-algebra automorphism of $A$.
				Finally, consider $F:A\longrightarrow M$ where $F(a)=e\cdot a=\gamma(a)e$. This map is bijective because $\gamma$ is bijective and $e$ is a left basis element.
		It is right $A$-linear, since $F(ab)=e\cdot(ab)=(e\cdot a)\cdot b=F(a)\cdot b$. Therefore, $M\simeq A$ as right $A$-modules.
	\end{proof}

	For each $\rho\in\Aut\J$, let $\widetilde \rho\in\Aut T$
	be the automorphism obtained by applying $\rho$ to each
	matrix entry; explicitly,
	$\widetilde \rho((a_{ij}))=(\rho(a_{ij}))$ for every
	$(a_{ij})\in T=\Mr_2(\J)$.

	\begin{proposition}\label{prop:standard}
		Let $\theta\in\operatorname{Aut} T$ satisfy
		$\theta(wT)=wT$. Then there exist
		$\rho\in\operatorname{Aut} \mathbf J$ and
		$H\in\Mr_2(\mathbf J)^\times$ such that
		$
		\theta=\operatorname{Ad}_H\circ\widetilde\rho.
		$
	\end{proposition}
	
	\begin{proof}
		We identify $w$ with $w\Id_2$ and write
		$\widetilde T=\Mr_2(A)=T[w^{-1}]$. Since $w$ is regular and normal in $\mathbf J$, the same is true of $w\Id_2$ in $T$, and $T$ embeds in
		$\widetilde T$. We first extend $\theta$ to $\widetilde T$. Since
		$\theta(w)T=\theta(wT)=wT$, there exist $a,b\in T$ such that
		$\theta(w)=wa$ and $w=\theta(w)b$. Thus, $w=wab$, and cancellation gives $ab=1$. Moreover,
		if $ac=0$, then $\theta(w)c=wac=0$. The regularity of
		$\theta(w)$ implies $c=0$. Applying this to
		$a(ba-1)=0$ gives $ba=1$. Hence $a\in T^\times$,
		so $\theta(w)$ is invertible in $\widetilde T$, with
		inverse $a^{-1}w^{-1}$.
		
		The universal property of Ore localization gives a
		unique $\Bbbk$-algebra homomorphism
		$\widetilde\theta:\widetilde T\to\widetilde T$ satisfying
		$
		\widetilde\theta|_T=\theta,
		\quad
		\widetilde\theta(w^{-1})=\theta(w)^{-1}.
		$
		The same construction applies to $\theta^{-1}$, since
		$\theta^{-1}(wT)=wT$. The two extensions are inverse:
		their compositions fix $T$ and $w^{-1}$, which generate
		$\widetilde T$. Consequently,
		$\widetilde\theta\in\Aut \widetilde T$.
		Let $e=e_{11}$, $f=\theta(e)$ and $g_{ij}=\theta(e_{ij})$ where $1\leq i,j\leq2$. The images of the matrix units satisfy $g_{ij}g_{k\ell}=\delta_{jk}g_{i\ell}$, $g_{11}=f$, $g_{22}=1-f$. In particular, $f^2=f$.
		Consider $P=fTe$ as a right $\mathbf J$-module, with
		action $p\cdot j=p(je)$. The first-column of $Te$ identifies
		free right module $\mathbf J \times \mathbf J=\mathbf J^2$ and
		$P$ with $f\mathbf J^2$. Since $f$ and $1-f$ are orthogonal idempotent whose sum is one, we have $\mathbf J^2=f\mathbf J^2\oplus(1-f)\mathbf J^2$. Thus, $P$ is finitely generated and projective as a right $\mathbf J$-module. Tensoring this split
		decomposition with $A$ shows that
	$	P\otimes_{\mathbf J}A
		\simeq fA^2
		\simeq f\widetilde T e
		$
		as right $A$-modules. Indeed, under
		$\mathbf J^2\otimes_{\mathbf J}A\simeq A^2$, the tensor
		of the projection $v\mapsto fv$ is again multiplication
		by $f$. Explicitly, the preceding isomorphism is $p\otimes a\longmapsto p(ae)$. 
		
		We next equip $\widetilde P:=f\widetilde T e$ with an
		invertible $\Bbbk$-central $A$-bimodule structure. Define the
		ring isomorphisms $\kappa:A\to e\widetilde T e$ and
		$\eta:A\to f\widetilde T f$ by
		\[
		\kappa(a)=ae,
		\qquad
		\eta(a)=\widetilde\theta(ae), \qquad a\in A.
		\]
		The map $\eta$ is an isomorphism because $\widetilde\theta$
		is an automorphism and $\widetilde\theta(e)=f$.
		Set $\widetilde Q=e\widetilde T f$ and define
		\[
		a\cdot p=\eta(a)p,\qquad
		p\cdot a=p\kappa(a),\qquad
		a\cdot q=\kappa(a)q,\qquad
		q\cdot a=q\eta(a),
		\]
		for $a\in A$, $p\in\widetilde P$, and $q\in\widetilde Q$.
		Associativity shows that these are commuting left and
		right actions. They are unital because $fp=p=pe$ and
		$eq=q=qf$. Multiplication induces $A$-bimodule homomorphisms
		$\mu:\widetilde P\otimes_A\widetilde Q\to A$ and
		$\nu:\widetilde Q\otimes_A\widetilde P\to A$ given by
		\[
		\mu(p\otimes q)=\eta^{-1}(pq),
		\qquad
		\nu(q\otimes p)=\kappa^{-1}(qp).
		\]
		Also, $pq\in f\widetilde T f$ and $qp\in e\widetilde T e$.
		The balancing conditions follow from
		$(p\kappa(a))q=p(\kappa(a)q)$ and
		$(q\eta(a))p=q(\eta(a)p)$.
		Compatibility with the outer actions follows from the
		multiplicativity of $\eta$ and $\kappa$.
		For completeness, define
		\[
		L(a)=\sum_{i=1}^2 fe_{i1}\otimes e_{1i}\eta(a),
		\qquad
		N(a)=\sum_{i=1}^2 eg_{i1}\otimes g_{1i}\kappa(a).
		\]
		In each summand defining $L(a)$ and $N(a)$, the tensor factors
		belong to $\widetilde P$ and $\widetilde Q$ in the indicated order. Using $\sum_{i=1}^2e_{i1}e_{1i}=1$ and $\sum_{i=1}^2g_{i1}g_{1i}=1$, we obtain $\mu(L(a))=\eta^{-1}(f\eta(a))=a$ and $\nu(N(a))=\kappa^{-1}(e\kappa(a))=a$. For the reverse compositions, note that
		$e_{1i}p\in e\widetilde T e=\kappa(A)$ and
		$g_{1i}q\in f\widetilde T f=\eta(A)$. The balancing relations therefore give
		\[
		\begin{aligned}
			L(\mu(p\otimes q))
			&=\sum_{i=1}^2fe_{i1}\otimes e_{1i}pq
			=\sum_{i=1}^2fe_{i1}e_{1i}p\otimes q
			=p\otimes q,\\
			N(\nu(q\otimes p))
			&=\sum_{i=1}^2eg_{i1}\otimes g_{1i}qp
			=\sum_{i=1}^2eg_{i1}g_{1i}q\otimes p
			=q\otimes p.
		\end{aligned}
		\]
		Thus $\mu$ and $\nu$ are isomorphisms, and
		$\widetilde P$ is an invertible $A$-bimodule with
		inverse $\widetilde Q$.
				For $\lambda\in\Bbbk$, the $\Bbbk$-linearity of
		$\widetilde\theta$ gives $\eta(\lambda)=\lambda f$.
		Consequently,
		$
		\lambda\cdot p=(\lambda f)p=\lambda p
		=p(\lambda e)=p\cdot\lambda.
		$
		Thus, $\widetilde P$ is $\Bbbk$-central. The same
		calculation shows that $\widetilde Q$ is
		$\Bbbk$-central. Lemma \ref{lem:picard-localization} now implies that $\widetilde P$ is free of rank one as a right
		$A$-module. Hence, $P\otimes_{\mathbf J}A\simeq A$. Since $P$ is finitely generated and projective,
		Lemma \ref{lem:descent} gives $P\simeq\mathbf J$ as
		right $\mathbf J$-modules. In particular,
		$f\mathbf J^2$ is free of rank one.
		
		Choose a right $\mathbf J$-basis element $v_1$ of
		$f\mathbf J^2$ and put $v_2=g_{21}v_1$. The relations $g_{21}=(1-f)g_{21}f$, $g_{12}=fg_{12}(1-f)$, together with
		$g_{12}g_{21}=f$, $g_{21}g_{12}=1-f$, show that multiplication by $g_{21}$ and $g_{12}$ defines mutually inverse right $\mathbf J$-linear maps between $f\mathbf J^2$ and $(1-f)\mathbf J^2$. Therefore, $v_2$ is a basis element of $(1-f)\mathbf J^2$, and $v_1,v_2$ form a right $\mathbf J$-basis of
		$\mathbf J^2$. Let $H$ be the matrix with columns $v_1,v_2$.
		Multiplication by $H$ is a bijective right
		$\mathbf J$-linear map on $\mathbf J^2$, and its
		inverse is represented by a matrix over $\mathbf J$.
		Hence $H\in\Mr_2(\J)^\times$.
		For $i,j\in\{1,2\}$, the relations among the elementary matrix  $e_{ij}$ give
		\[
		\begin{aligned}
			g_{ij}v_1
			=g_{ij}fv_1
			=\delta_{j1}g_{i1}v_1
			=\delta_{j1}v_i,\qquad	g_{ij}v_2
			=g_{ij}g_{21}v_1
			=\delta_{j2}g_{i1}v_1
			=\delta_{j2}v_i.
		\end{aligned}
		\]
		Comparing columns yields $g_{ij}H=He_{ij}$. Thus
		$\theta(e_{ij})=He_{ij}H^{-1}$.	We fix $\tau:=\operatorname{Ad}_{H^{-1}}\circ \,\theta$.
		Then $\tau(e_{ij})=e_{ij}$ for every $i,j$.
		In particular, $\tau$ restricts to an automorphism
		of $eTe$. Through the identification
		$\mathbf J\simeq eTe$, this determines
		$\rho\in\operatorname{Aut} \mathbf J$ such that $\tau(je)=\rho(j)e$, with $j\in\mathbf J$. Every matrix $t=(t_{ij})\in T$ can be written as
		$
		t=\sum_{i,j=1}^2e_{i1}(t_{ij}e)e_{1j}.
		$
		Since $\tau$ fixes all elementary matrix $e_{ij}$, it follows that, for all $t\in T$,
		$
		\tau(t)
		=\sum_{i,j=1}^2e_{i1}(\rho(t_{ij})e)e_{1j}
		=\widetilde\rho(t).
		$
		Thus, $\theta=\operatorname{Ad}_H\circ\widetilde\rho$.
	\end{proof}

\subsection{Proof of Theorem A} Using the preceding results, we prove Theorem~A
first over algebraically closed fields and then,
by descent, over arbitrary fields of characteristic zero.

	\begin{theorem}\label{thm:main}
		Let $\kk$ be an algebraically closed field of characteristic
		zero. Then,
		\[
		\Aut\B=\Inn\B\rtimes\Aut_{\mathrm{super}}\B.
		\]
	\end{theorem}
	
		\begin{proof}
		Identify $\B$ with its image in $T$ under $\iota$ and let
		$\varphi\in\Aut\B$. By Proposition~\ref{prop-quotient},
		there exists $\mu\in\kk^\times$ such that
		$\varphi(y)-\mu y\in I$. By Lemma~\ref{lem:idealizer},
		$\varphi$ extends to an automorphism $\theta$ of $T$
		satisfying $\theta(wT)=wT$. Proposition~\ref{prop:standard}
		gives $\theta=\Ad_H\circ\widetilde\rho$ for some $H\in\Mr_2(\J)^\times$ and $\rho\in\Aut\J$.
		By the description of the automorphisms of the Jordan plane
		recalled in \cite[Remark~5.4]{AD}, there exist
		$c\in\kk^\times$ and $q\in\kk[t]$ such that
		\[
		\rho(w)=cw,\quad\quad \rho(s)=cs+q(w).
		\]
	For $M\in T$, write $\overline M$ for its image in
	$T/wT\simeq\Mr_2(\Bbbk[s])$, obtained by reducing each entry
	modulo $w\J$. Since $\theta(Y)-\mu Y\in I=wT$, reducing
	$\theta(Y)=H\widetilde\rho(Y)H^{-1}$ gives
	\begin{equation}\label{eq:reduced-conjugacy}
		\overline H
		\begin{pmatrix}0&cs+q(0)\\1&0\end{pmatrix}
		(\overline H)^{-1}
		=\mu\overline Y.
	\end{equation}
		Squaring both sides in $\Mr_2(\kk[s])$ gives $(cs+q(0))\Id_2=\mu^2s\Id_2$. Comparing coefficients in $\kk[s]$, we obtain $c=\mu^2$ and $q(0)=0$.
				Let $D=\operatorname{diag}(\mu,1)$. Then,
		\[
		\overline{\widetilde\rho(Y)}
		=
		\begin{pmatrix}0&\mu^2s\\1&0\end{pmatrix}
		=
		D(\mu\overline Y)D^{-1}.
		\]
		It follows from \eqref{eq:reduced-conjugacy} that
		$\overline H D$ commutes with $\overline Y$. A matrix $M=(m_{ij})\in \Mr_2(\kk[s])$ commutes with
		$\overline Y$ if and only if $m_{22}=m_{11}$ and $m_{12}=s m_{21}$. Since $\overline Y^2=s\Id_2$, the centralizer of
		$\overline Y$ is therefore
		\[
		\left\{
		\begin{pmatrix}a(s)&s b(s)\\b(s)&a(s)\end{pmatrix}
		:a,b\in\kk[s]
		\right\}
		=\kk[\overline Y].
		\]
		Moreover, $\kk[\overline Y]\simeq\kk[t]$: if
		$a(\overline Y^2)+b(\overline Y^2)\overline Y=0$,
		the $(2,2)$ and $(2,1)$-entries give $a(s)=b(s)=0$.
		Since both $\overline H D$ and its inverse commute
		with $\overline Y$, we conclude that $\overline H D=\lambda\Id_2$ for some $\lambda\in\kk^\times$.
		
		Let $U=HD$. Both $U$ and $U^{-1}$ belong to $T$, and
		their reductions are $\lambda\Id_2$ and
		$\lambda^{-1}\Id_2$, respectively. Both reductions
		belong to $\kk[\overline Y]$. Hence
		Lemma~\ref{lem:matrix-model} gives $U,U^{-1}\in\B$,
		so $U\in\B^\times$.		
		Define $\sigma=\Ad_{D^{-1}}\circ \,\widetilde\rho$.
		Then $\theta=\Ad_U\circ \,\sigma$. Since $q(0)=0$,
		write $q(t)=tr(t)$ with $r\in\kk[t]$, and let
		$p(t)=\mu^{-1}r(t)$. Direct computation gives
		\[
		\sigma(X)=\mu X,\quad\quad
		\sigma(Y)=
		\begin{pmatrix}
			0&\mu s+\mu^{-1}q(w)\\
			\mu&0
		\end{pmatrix}
		=\mu Y+p(w)X.
		\]
	Thus, $\sigma|_{\B}=\sigma_{p,\mu}$. Restricting
		$\theta=\Ad_U\circ \,\sigma$ to $\B$, we obtain $\varphi=\Ad_U\circ \,\sigma_{p,\mu}$. It remains to show that the intersection of these two
		subgroups is trivial. Since $\sigma_{p,\mu}=\sigma_{p,0,\mu}$, the group $\Aut_{\mathrm{super}}\B$ is contained in $\mathbb L$. By Corollary~\ref{cor-quotient-L-inner}, the homomorphism
		$\Phi$ restricts to the identity on
		$\Aut_{\mathrm{super}}\B$ and has kernel
		$N=\Inn\B\cap\mathbb L$. Hence,
		\[
		\Inn\B\cap\Aut_{\mathrm{super}}\B
		=
		N\cap\Aut_{\mathrm{super}}\B
		=
		\{\operatorname{id}_{\B}\}.
		\]
			Since $\Inn\B$ is normal in $\Aut\B$, the factorization
		above gives the asserted semidirect product $\Aut\B=\Inn\B\rtimes\Aut_{\mathrm{super}}\B$.
	\end{proof}
	
	The following corollary completes the proof of Theorem~A
	by extending Theorem~\ref{thm:main} to arbitrary fields
	of characteristic zero.
	
	\begin{corollary}[Theorem~A]\label{cor:descent}
		Let $\B$ be the super Jordan plane over a field $\Bbbk$ of
		characteristic zero. Then
		$
		\Aut\B=\Inn\B\rtimes\Aut_{\mathrm{super}}\B.
		$
	\end{corollary}

\begin{proof}
	Let $\overline{\kk}$ be an algebraic closure of $\kk$, set
	$G=\operatorname{Gal}(\overline{\kk}/\kk)$, and write
	$\overline{\B}=\overline{\kk}\otimes_{\kk}\B$.
	Since $\kk$ has characteristic zero,
	\cite[Proposition~7.9]{Milne} gives
	\begin{equation}\label{eq:galois-fixed-field}
		\overline{\kk}^{G}=\kk.
	\end{equation}
	We identify $\B$ with its canonical image in $\overline{\B}$.
	Each $\gamma\in G$ acts on $\overline{\B}$ through
	$\widetilde{\gamma}(\lambda\otimes b)=\gamma(\lambda)\otimes b$.
	By the PBW basis and \eqref{eq:galois-fixed-field}, the common
	fixed subring of these automorphisms is $\B$.
	Let $\varphi\in\Aut\B$. Applying Theorem~\ref{thm:main} to its
	scalar extension gives
	\begin{equation}\label{eq:descent-factor}
		\id_{\overline{\kk}}\otimes\varphi
		=\Ad_U\circ\sigma_{p,\mu},
	\end{equation}
	for some $U\in\overline{\B}^{\times}$,
	$p\in\overline{\kk}[t]$, and $\mu\in\overline{\kk}^{\times}$.
	The left-hand side commutes with every $\widetilde{\gamma}$.
	Conjugating by $\widetilde{\gamma}$ and using the uniqueness
	of the two factors in Theorem~\ref{thm:main}, we obtain
	\[
	\sigma_{\gamma(p),\gamma(\mu)}=\sigma_{p,\mu},
	\qquad
	\Ad_{\widetilde{\gamma}(U)}=\Ad_U.
	\]
	The uniqueness of the parameters of $\sigma_{p,\mu}$ now gives
	$\gamma(p)=p$ and $\gamma(\mu)=\mu$ for every $\gamma\in G$.
	Thus $p\in\kk[t]$ and $\mu\in\kk^\times$, so $\sigma_{p,\mu}$
	restricts to a super-automorphism of $\B$.
	To descend the inner factor, note that
	$U^{-1}\widetilde{\gamma}(U)$ is central in $\overline{\B}$.
	Since the center of $\overline{\B}$ is $\overline{\kk}$ by
	\cite[Proposition~3.1(1)]{Lu2026SJP}, there exists
	$c_\gamma\in\overline{\kk}^{\times}$ such that
	\begin{equation}\label{eq:descent-scalar}
		\widetilde{\gamma}(U)=c_\gamma U.
	\end{equation}
	Under the canonical projection
	$\pi:\overline{\B}\to\overline{\B}/(x)\simeq\overline{\kk}[y]$,
	the unit $U$ maps to a nonzero scalar $c=\pi(U)$.
	Applying $\pi$ to \eqref{eq:descent-scalar} yields
	$\gamma(c)=c_\gamma c$. Consequently, $V=c^{-1}U$ satisfies
	$
	\widetilde{\gamma}(V)
	=\gamma(c)^{-1}c_\gamma U=V$, $\gamma\in G$.
	Both $V$ and its inverse therefore belong to $\B$.
	Since $\Ad_V=\Ad_U$, restricting \eqref{eq:descent-factor}
	gives $\varphi=\Ad_V\circ\sigma_{p,\mu}$ with $V\in\B^\times$.
	Finally, any element of
	$\Inn\B\cap\Aut_{\mathrm{super}}\B$ extends to an element of
	$\Inn\overline{\B}\cap\Aut_{\mathrm{super}}\overline{\B}$,
	which is trivial by Theorem~\ref{thm:main}.
	The asserted semidirect product decomposition follows.
\end{proof}

The following example illustrates Theorem~\ref{thm:main}
by explicitly determining the inner and super-automorphism factors.

\begin{example}\label{ex:constructive-decomposition}
	Consider the automorphism $\varphi\in\Aut\B$ given by
	\begin{align*}
		\varphi(x)
		&=2x+2z^2-2xzy-2z^3y,\\
		\varphi(y)
		&=2y+xz+2zy^2-2xzy-2xy^3+z^3-xz^2y
		-2z^3y-2z^2y^3-z^4y.
	\end{align*}
	We determine its inner and super-automorphism factors using
	the construction in Theorem~\ref{thm:main}.
	Let $\theta$ be the extension of $\varphi$ to $T=\Mr_2(\J)$
	provided by Lemma~\ref{lem:idealizer}. In the matrix realization
	$
	X=we_{12}$\ and $
	Y=se_{12}+e_{21}$. Set $F=\theta(X)$, $G=\theta(Y)$ and $a_0=2s+w^2$.
	Substitution gives
	\[
	F=-2w^2e_{11}+2we_{12}-2w^3e_{21}+2w^2e_{22},
	\qquad
	G=-a_0we_{11}+a_0e_{12}
	+(2-wa_0w)e_{21}+wa_0e_{22}
	\]
	and $\theta(w)=FG+GF=4w\Id_2$.
	To construct $H$, we need only the images
	$g_{11}=\theta(e_{11})$ and $g_{21}=\theta(e_{21})$.
	Since $e_{11}=w^{-1}XY$ and $e_{21}=w^{-1}YXY$, we obtain
	\begin{equation}\label{eq:example-matrix-units}
		g_{11}=(4w)^{-1}FG=e_{11}+we_{21},
		\qquad
		g_{21}=(4w)^{-1}GFG=2e_{21}.
	\end{equation}
	Thus, for $f=g_{11}$, the right $\J$-module
	\[
	f\J^2=\{(a,wa)^{\mathsf T}:a\in\J\}
	\]
	has basis vector $v_1=(1,w)^{\mathsf T}$.
	Taking $v_2=g_{21}v_1=(0,2)^{\mathsf T}$ gives
	\[
	H=e_{11}+we_{21}+2e_{22},
	\qquad
	H^{-1}=e_{11}-\frac w2e_{21}+\frac12e_{22}.
	\]
	By Proposition~\ref{prop:standard},
	$\theta=\Ad_H\circ\widetilde\rho$ for some $\rho\in\Aut\J$.
	Since $H$ commutes with $w\Id_2$, we have $\rho(w)=4w$.
	Moreover,
	$
	H^{-1}GH=2a_0e_{12}+e_{21},
	$
	so $Y^2=s\Id_2$ gives $\rho(s)\Id_2=H^{-1}G^2H=2a_0\Id_2$.
	Therefore,
	$		\rho(w)=4w$ and $\rho(s)=4s+2w^2$.	Modulo the ideal $I=(x)$, the automorphism $\varphi$ sends
	$\overline y$ to $2\overline y$, so $\mu=2$.
	Following the proof of Theorem~\ref{thm:main}, set
	$D=2e_{11}+e_{22}$ and $U=HD$. Then
	\begin{equation*}
		U=2e_{11}+2we_{21}+2e_{22}
		=2(\Id_2+YXY),
		\qquad
		U^{-1}=\tfrac12(\Id_2-YXY).
	\end{equation*}
	Hence $U$ corresponds to the unit $u=2(1+yxy)$ of $\B$.
	Finally, $\sigma=\Ad_{D^{-1}}\circ\widetilde\rho$ restricts
	to $\sigma_{t,2}$ on $\B$, since
	$\rho(s)=\mu^2s+q(w)$ with $q(t)=2t^2$ and
	$p(t)=\mu^{-1}q(t)/t=t$.
	As $\Ad_u=\Ad_{1+yxy}$, the desired factorization is
	$
	\varphi=\Ad_{1+yxy}\circ\sigma_{t,2}$.
\end{example}

\section{Nilpotent elements and square-zero loci}

 In this section, we prove that every nilpotent element of $\mathcal B$ has square zero and obtain an explicit form for these elements using a matrix identification. We use this description to study the linear families
	$S_m$ and $H_m$ and their relation to the quadratic
	Veronese varieties contained in the projective square-zero
	loci. Finally, we exhibit a mixed nilpotent element of degree $5$ that lies outside the union of the two linear families and the Veronese family. For the geometric discussion, we assume that the ground field is algebraically closed.

\subsection{Nilpotent elements have square zero}
	
	We first verify that localization at the powers of $w$
	embeds $\mathcal B$ in $\mathcal B[w^{-1}]$.
	
	\begin{lemma}\label{lem-x21-degree}
		Let $0\neq a\in\mathcal{B}$ be homogeneous of degree $m$.
		Then, $w^na$ is nonzero and homogeneous of degree $m+2n$,
		for every $n\geq0$. Furthermore, $w$ is a regular normal element of $\mathcal{B}$.
	\end{lemma}
	
	\begin{proof} 
		The case $n=0$ is trivial. Assume $n>0$.
		Notice that, by \eqref{eq:x21_power}, one has $	w^n=(xy)^n+z^n=xz^{n-1}y+z^n$.
		Write the homogeneous element $a$ in the PBW basis as
		$
		a=\sum_{2j+k=m}\alpha_{j,k}z^jy^k
		+\sum_{1+2j+k=m}\beta_{j,k}xz^jy^k.
		$
		Multiplying $a$ on the left by $w^n$, we obtain
		\[
		w^na=
		\sum_{2j+k=m}\alpha_{j,k}z^{n+j}y^k
		+\alpha_{0,m}xz^{n-1}y^{m+1}
		+\sum_{1+2j+k=m}\beta_{j,k}xz^{n+j}y^k.
		\]
		Indeed, the term $xz^{n-1}yz^jy^k$ vanishes whenever
		$j>0$, and for $j=0$ it gives $xz^{n-1}y^{k+1}$.
		Also, since $z^nx=0$ for $n>0$, we have $w^nxz^jy^k=(xz^{n-1}y+z^n)xz^jy^k=xz^{n+j}y^k$.
		Furthermore, all PBW monomials appearing in $w^na$ have degree $m+2n$. Moreover, no cancellation between the three parts is possible: the first part consists of PBW
		monomials not beginning with $x$, the middle term is
		the only possible PBW monomial ending with $y^{m+1}$
		and the last part consists of PBW monomials beginning
		with $x$ and ending with a power of $y$ at most $m-1$.
		Since $a\neq0$, at least one coefficient $\alpha_{j,k}$
		or $\beta_{j,k}$ is nonzero and hence $w^na\neq0$.
		Thus, $w^na$ is nonzero and homogeneous of
		degree $m+2n$. Taking $n=1$, we obtain $wa\neq0$ for every nonzero
		homogeneous element $a$. Since multiplication by $w$ shifts
		each homogeneous degree by $2$, the images of distinct
		homogeneous components cannot cancel. By decomposing an
		arbitrary element into homogeneous components, we conclude
		that left multiplication by $w$ is injective.
		
		To prove regularity on the other side, consider the
		automorphism $\tau\in \Aut\mathcal{B}$ determined by $\tau(x)=x$ and $\tau(y)=y+x$.
		The identities $xw=wx$ and $yw=w(y+x)$ imply $bw=w\tau(b)$, for every $b\in\mathcal{B}$. Thus, if $bw=0$, the injectivity of left multiplication
		by $w$ gives $\tau(b)=0$, and hence $b=0$.
		Therefore, $w$ is regular on both sides.
		Finally, since $\tau$ is surjective, the same identity
		gives $\mathcal{B}w=w\mathcal{B}$, proving that $w$
		is normal.
	\end{proof}
	
	Since $w$ is regular and normal, its powers form a left and
	right Ore set of regular elements. We therefore identify $\mathcal B$ with a subalgebra of $\mathcal B[w^{-1}]$. In \cite{Lu2026SJP}, Tao Lu studied the prime spectrum and representation theory of the super Jordan plane. He showed, among other things, that the localization of $\mathcal B$ at the powers of $w=xy+yx$ is isomorphic to a $2\times2$ matrix algebra over a localization of the first Weyl algebra $A_1(\Bbbk)$. This leads to the following description of the classical ring of quotients of $\mathcal B$ in characteristic zero.

	\begin{theorem}\label{thm.Lu}
	Let $D_1(\Bbbk)$ denote the quotient division ring of
	the first Weyl algebra $A_1(\Bbbk)$. The classical ring of quotients $	Q(\mathcal B)$ of $\mathcal B$ is isomorphic to	$\Mr_2\bigl(D_1(\Bbbk)\bigr)$.
		In particular, the canonical map $\mathcal B\rightarrow Q(\mathcal B)$ is injective.
	\end{theorem}
	
	\begin{proof}
		The isomorphism $Q(\mathcal B)\simeq \Mr_2\bigl(D_1(\Bbbk)\bigr)$  was proved in \cite[Proposition 3.1(3)]{Lu2026SJP}. The canonical map is injective because the classical ring of quotients is obtained by localizing at regular elements.
	\end{proof}
	
	Andruskiewitsch and Dumas proved in \cite[Proposition 4.2]{AD}
	that every nilpotent element of the odd component
	$\mathcal B_{\bar{1}}$ has square zero. The following theorem shows that the same conclusion holds
	for every nilpotent element of $\mathcal B$.
	
	\begin{theorem}\label{nilpotentssquarezero}
		Let $a\in\mathcal B$ be nilpotent. Then $a^2=0$. Consequently, every nonzero nilpotent element of $\mathcal B$ has nilpotency index exactly $2$.
	\end{theorem}
	
	\begin{proof}
		By Theorem \ref{thm.Lu}, we may regard $\mathcal B$ as
		a subalgebra of $\Mr_2\bigl(D_1(\Bbbk)\bigr)$. Let $a\in\mathcal B$ be nilpotent and denote by $M$ its image in
		this matrix algebra. Then $M$ is a nilpotent endomorphism of the two-dimensional right $D_1(\Bbbk)$-vector space
		$D_1(\Bbbk)^2$. We claim that $M^2=0$. The assertion is trivial if $M=0$, so suppose
		that $M\neq0$. Since $M$ is nonzero and nilpotent, let $m\geq2$ be the least
		integer such that $M^m=0$. Then $M^{m-1}\neq0$, so there exists
		$v\in D_1(\Bbbk)^2$ such that $u:=M^{m-1}(v)\neq0$.	Moreover, since $m\geq2$, $u\in\operatorname{Im}M$.
		Thus, $0\neq u\in\ker M\cap\operatorname{Im}M$. The rank-nullity theorem over the division ring $D_1(\Bbbk)$,
		applied to $M:D_1(\Bbbk)^2\rightarrow D_1(\Bbbk)^2$, gives
		\[
		\dim_{D_1(\Bbbk)}\operatorname{Im}M
		=
		2-\dim_{D_1(\Bbbk)}\ker M
		\leq1.
		\]
		Since $\operatorname{Im}M$ contains the nonzero vector $u$,
		it follows that $\operatorname{Im}M=uD_1(\Bbbk)\subseteq\ker M$. 
		Thus, $M^2=0$. Since $\mathcal B$ embeds into
$\Mr_2(D_1(\Bbbk))$, we conclude that $a^2=0$.
		\end{proof}
	
%

\subsection{Parametrization of nilpotent elements}
	\label{subsec:matrix-nilpotents}
	
	We use Lu's matrix representation given in Theorem \ref{thm.Lu} to describe the nilpotent
	elements of $\mathcal B$. We retain the notation
	\[
	u=y^2-xy,
	\qquad\quad\quad
	\mathcal J=\Bbbk\langle u,w\rangle.
	\]
	This is the Jordan plane contained in $C_{\B}(x)$
	considered in Lemma~\ref{lemma-centralizador-x1}.
	The notation $\mathbf J$ continues to denote the distinct
	subalgebra $\Bbbk\langle y^2,w\rangle$, which is used in
	the parity description of \cite{AD}. Since $w$ is a regular normal element of $\mathcal B$, the canonical map from $\mathcal B$ to its localization $\mathcal B[w^{-1}]$ at the powers of $w$ is injective. In $\mathcal B[w^{-1}]$, denote
	$
	\partial=(yw^{-1})^2$ and $
	A=\Bbbk\langle w,w^{-1},\partial\rangle.
	$
	By \cite[Lemma 2.1]{Lu2026SJP}, with Lu's element $s$ denoted
	here by $w$, the algebra $A$ is the localization at the powers
	of $w$ of the first Weyl algebra
	$
	A_1(\Bbbk)
	=
	\Bbbk\langle w,\partial\mid \partial w-w\partial=1\rangle.
	$
	Moreover, \cite[Lemma 2.1]{Lu2026SJP} gives an isomorphism $\Phi:\mathcal B[w^{-1}]\to \Mr_2(A)$ satisfying $\Phi(x)=e_{21}$,
	$\Phi(w)=w\Id_2$,
	$\Phi(yw^{-1})=e_{12}+\partial e_{21}$.
	In particular, $\Phi(\partial)=\partial \Id_2$, so $\Phi(a)=a\Id_2$ for all $a\in A$.
	Since $A$ is a localization of the first Weyl algebra, we identify its quotient division ring $\operatorname{Frac}(A)$ with $D_1(\Bbbk)$. Restricting $\Phi$ to $\mathcal B$ and
	then using the inclusion $\Mr_2(A)\subseteq \Mr_2(D_1(\Bbbk))$ gives the
	injective homomorphism
	\begin{align}\label{injec-rho}
	\rho:\mathcal B
		\hookrightarrow\mathcal B[w^{-1}]
		\xrightarrow{\Phi}\Mr_2(A)
		\hookrightarrow \Mr_2(D_1(\Bbbk)).
	\end{align}
		We choose the identification of $Q(\mathcal B)$ with
	$\Mr_2(D_1(\Bbbk))$ induced by $\Phi$. Under this identification, $\rho$ is the canonical
	embedding of $\mathcal B$ into its classical ring
	of quotients.
	
	Since $y=(yw^{-1})w$, we obtain
	$
	\rho(y)
	=(e_{12}+\partial e_{21})w\Id_2
	=we_{12}+(\partial w)e_{21}$.
	From $wyw^{-1}=y-x$ it follows that $
	w\partial w
	=
	w(yw^{-1})^2w
	=
	(wyw^{-1})y
	=
	(y-x)y
	=
	u.
	$
Since $u=w\partial w$ and $w,\partial\in A$,
	the elements $u$ and $w$ can be identified
	with elements of $A$. Moreover, the matrix
	representation satisfies
	$
	\rho(u)=(w\partial w)\Id_2$, $
	\rho(w)=w\Id_2.
	$
	Since $\rho$ is injective, its restriction
	to $\mathcal J$ identifies $\mathcal J$
	with a subalgebra of $A$ and hence of $D_1(\Bbbk)$.
	Under this identification, we write $\rho(f)=f\Id_2$, with $f\in\mathcal J$. Recalling that $z=yx$, direct multiplication gives
	$
	\rho(z)=we_{11}$,
	$\rho(xy)=we_{22}$,
	$\rho(zy)=w^2e_{12}$	and $
	\rho(y^2)=ue_{11}+(u+w)e_{22}$.
	 For the last
	identity, we used $\partial w^2=w\partial w+w=u+w$.
	
	\begin{remark} 
		The representation $\rho$ is related to the embedding
		$\iota$ of Lemma~\ref{lem:matrix-model} by conjugation
		after identifying the two copies of the Jordan plane.
		Let
		\[
		\chi:\mathbf J\to\mathcal J,
		\qquad
		\chi(w)=w,\qquad \chi(s)=u,
		\]
		be the isomorphism determined by the Jordan relations
		$sw-ws=w^2$ and $uw-wu=w^2$.
		Using $\mathcal J\subseteq D_1(\Bbbk)$, let $\widetilde\chi:M_2(\mathbf J)\longrightarrow
		M_2(D_1(\Bbbk))$ be the homomorphism obtained by applying $\chi$ to each matrix entry. Consider $C=we_{12}+e_{21}$. Since $u=w\partial w$, direct multiplication gives 
		$C\rho(x)C^{-1}=we_{12}=\widetilde\chi(\iota(x))$
		and
		$C\rho(y)C^{-1}=ue_{12}+e_{21}=\widetilde\chi(\iota(y))$.
			As $x$ and $y$ generate $\mathcal B$, it follows that
		$\widetilde\chi(\iota(b))
		=C\rho(b)C^{-1}$, with $b\in\mathcal B$. This comparison also explains the different scalar
		matrices occurring in the two representations.
		Indeed, $\rho(u)=u\Id_2$, whereas $\iota(y^2)=s\Id_2$.
		Since $wuw^{-1}=u-w$, we have $C\rho(u)C^{-1}=(u-w)e_{11}+ue_{22}$
		and $C\rho(y^2)C^{-1}=u\Id_2$.
		The relation $uw-wu=w^2$ gives $wuw^{-1}=u-w$. Thus, conjugation by $C$ sends the scalar matrix $u\Id_2$
		to the nonscalar diagonal matrix $\operatorname{diag}(u-w,u)$. In particular, being represented by a scalar matrix
		depends on the chosen matrix coordinates.
	\end{remark}

	\begin{lemma}\label{lem:matrix-image-nilpotent-ideal}
		The algebra $\mathcal J=\Bbbk\langle u,w\rangle$ has PBW basis
		$\{w^ju^k:j,k\geq0\}$, 
		$uw-wu=w^2$, $w\mathcal J=\mathcal Jw$ and
		$\operatorname{Frac}(\mathcal J)=D_1(\Bbbk)$.
		Moreover,
		$\rho(I)=\mathcal Jw e_{11}
		+\mathcal Jw^2 e_{12}
		+\mathcal J e_{21}
		+\mathcal Jw e_{22}$.
	\end{lemma}
	
	\begin{proof}
		From $u=w\partial w$ and $\partial w-w\partial=1$ follow 
		$
		uw-wu
		=
		w(\partial w-w\partial)w
		=
		w^2.
		$
		This relation shows that the monomials $w^ju^k$ span
		$\mathcal J$. Since $u=w^2\partial+w$, ordering in the Weyl algebra gives, for $j,k\geq0$,
		\[
		w^ju^k
		=
		w^{j+2k}\partial^k
		+
		\sum_{\ell<k}f_\ell(w)\partial^\ell,
		\quad
		f_\ell(w)\in\Bbbk[w].
		\]
		Comparing the highest powers of $\partial$ in the Weyl PBW
		basis proves that these monomials are linearly independent. The identities $uw=w(u+w)$ and $wu=(u-w)w$ give $w\mathcal J=\mathcal Jw$. Since $\mathcal J$ is a Noetherian domain contained in $D_1(\Bbbk)$, its quotient division ring is naturally a division subring of $D_1(\Bbbk)$. Moreover,
		$
		\partial=w^{-1}uw^{-1}
		\in\operatorname{Frac}(\mathcal J).
		$
		Thus, $\operatorname{Frac}(\mathcal J)$ contains both $w$
		and $\partial$, which generate $D_1(\Bbbk)$ as a division ring. Consequently, $\operatorname{Frac}(\mathcal J)=D_1(\Bbbk)$. Finally, the PBW basis of $I$ and the preceding formulas yield,
		for all $j,k\geq0$,
		\[
		\begin{aligned}
			\rho(z^{j+1}y^{2k})&=w^{j+1}u^k e_{11},
			&\quad
			\rho(z^{j+1}y^{2k+1})&=w^{j+1}u^kw e_{12},\\
			\rho(xz^jy^{2k})&=w^ju^k e_{21},
			&
			\rho(xz^jy^{2k+1})&=w^ju^kw e_{22}.
		\end{aligned}
		\]
		These four families are supported in distinct matrix entries
		and span, respectively, $w\mathcal J$, $w\mathcal Jw$, $\mathcal J$, $\mathcal Jw$. Their coefficients can therefore be chosen independently. Using that 
		$w\mathcal J=\mathcal Jw$ and $w\mathcal Jw=\mathcal Jw^2$,	we obtain the stated description of $\rho(I)$.
	\end{proof}

	\begin{lemma}\label{prop:matrix-criterion-nilpotent}
		Let $a\in\mathcal B$ and 
		$\rho(a)=\alpha e_{11}+\beta e_{12}
		+c e_{21}+\delta e_{22}$,
		where $\alpha,\beta,c,\delta\in D_1(\Bbbk)$.
		Then $a$ is nilpotent if and only if one of the
		following conditions holds:
		\begin{enumerate}[\rm (i)]
			\item $c=0$ and $\alpha=\delta=0$;
			\item $c\neq0$, $\beta=-\alpha^2c^{-1}$
			and $\delta=-c\alpha c^{-1}$.
		\end{enumerate}
	\end{lemma}
		\begin{proof}
		The result follows by a direct computation.
	\end{proof}

	\begin{proposition}\label{prop:nilpotent-rational-normal-form}
		Every nilpotent element of $\mathcal B$ has exactly one
		of the following forms:
		\begin{enumerate}[\rm (i)]
			\item $gzy$, where $g\in\mathcal J$;
			\item
			$a(r,c)=xc+(rcw^{-1})z-(crw^{-1})xy
			-(rcrw^{-2})zy$, where
			$r\in D_1(\Bbbk)$, $0\neq c\in\mathcal J$,
			$rc,cr\in\mathcal Jw$ and
			$rcr\in\mathcal Jw^2$.
		\end{enumerate}
		Conversely, every element of either form has square zero
		and the parameters $g$, $r$ and $c$ are unique in each case.
	\end{proposition}

	\begin{proof}
		Let $a$ be nilpotent. By Lemma \ref{lemanilpotentkernel}(i),
		$a\in I$. Write $\rho(a):=\alpha e_{11}+\beta e_{12}
		+c e_{21}+\delta e_{22}$. If $c=0$, Lemma \ref{prop:matrix-criterion-nilpotent}
		gives $\alpha=\delta=0$. Since $\beta\in\mathcal Jw^2$,
		there is a unique $g\in\mathcal J$ such that $\beta=gw^2$.
		Thus, $\rho(a)=\rho(gzy)$, and $a=gzy$. If $c\neq0$, denote $r=\alpha c^{-1}$. The Lemma \ref{prop:matrix-criterion-nilpotent} gives 
		$\rho(a)=rc\,e_{11}-rcr\,e_{12}+c\,e_{21}-cr\,e_{22}$.
		Moreover, Lemma \ref{lem:matrix-image-nilpotent-ideal} yields precisely the stated conditions on $r,c$. 
		Also, write $rc=f_1w$, $cr=f_2w$, and $rcr=f_3w^2$,
		where $f_1,f_2,f_3\in\mathcal J$.
		Then $a(r,c)=xc+f_1z-f_2xy-f_3zy\in\mathcal B$
		and the formulas of Lemma~\ref{lem:matrix-image-nilpotent-ideal}
		give
		$
		\rho(a(r,c))
		=rc\,e_{11}-rcr\,e_{12}+c\,e_{21}-cr\,e_{22}
		=\rho(a).
		$
		By the injectivity of $\rho$, we conclude that $a=a(r,c)$.

		Conversely, $\rho(gzy)=gw^2e_{12}$ has square zero.
		For admissible $r,c$, we have
		\[
		\rho(a(r,c))
		=rc\,e_{11}-rcr\,e_{12}+c\,e_{21}-cr\,e_{22}.
		\]
		A direct computation shows that $\rho(a(r,c))^2=0$.
		Since $\rho$ is injective, both forms yield square-zero
		elements of $\mathcal B$.
		The two cases are disjoint, as the $(2,1)$-entry of
		$\rho(a)$ is zero in the first case and nonzero in the
		second. Moreover, the parameters are uniquely determined
		by $\rho(a)$: in the first case,
		$g=\rho(a)_{12}w^{-2}$, while in the second,
		$c=\rho(a)_{21}$ and $r=\rho(a)_{11}c^{-1}$.
		This completes the proof.
		\end{proof}

\subsection{Linear subspaces of nilpotent elements}
	
	In this subsection, we introduce two linear families of square-zero elements of $\mathcal B$ and study the corresponding affine and projective square-zero locus. We also relate these families to the description of the even and odd components of the super Jordan plane
	obtained by Andruskiewitsch and Dumas in \cite{AD}.
	For $m\geq0$, consider the vector space
	\[
	S_m
	:=
	\Span
	\left\{
	xz^jy^{2k}
	:
	j,k\in\mathbb N_0,\ j+k\leq m
	\right\}
	\subseteq
	\mathcal B_{\leq2m+1},
	\]
	and, for $m\geq1$, let
	\[
	H_m
	:=
	\Span
	\left\{
	z^jy^{2k+1}
	:
	j,k\in\mathbb N_0,\ j\geq1,\ j+k\leq m
	\right\}
	\subseteq
	\mathcal B_{\leq2m+1}.
	\]
These subspaces consist entirely of square-zero elements. Indeed, for any two PBW monomials appearing in $S_m$, Lemma \ref{prop:product} gives $(xz^jy^{2k})(xz^{j'}y^{2k'})=0$. Then, $S_mS_m=0$. Similarly, since $j,j'\geq1$, Lemma \ref{prop:product} gives $(z^jy^{2k+1})(z^{j'}y^{2k'+1})=0$ and therefore $H_mH_m=0$. Consequently, every element of $S_m\cup H_m$ has square zero.
The PBW monomials spanning $S_m$ and $H_m$ are linearly independent. Moreover, $S_m\cap H_m=0$,
	$
	\dim S_m=\binom{m+2}{2}$ and $\dim H_m=\binom{m+1}{2}$. Recall that a regular element is neither
	a left nor a right zero divisor.

	\begin{proposition}\label{prop:sums-of-nilpotents}
		Let $a,b\in\mathcal B$ be nilpotent elements.
		Then $a+b$ is nilpotent if and only if $ab=0$,
		or equivalently, $ba=0$. Otherwise, $a+b$ is regular.
	\end{proposition}
	
	\begin{proof}
		The assertion is immediate if $a=0$ or $b=0$.
		Otherwise, since $a^2=b^2=0$, the nonzero matrices
		$\rho(a)$ and $\rho(b)$ have one-dimensional images
		and kernels over $D_1(\Bbbk)$, with
		$H_a:=\operatorname{Im}\rho(a)=\ker\rho(a)$ and
		$H_b:=\operatorname{Im}\rho(b)=\ker\rho(b)$.
		Consequently,
		\[
		ab=0\iff H_b\subseteq H_a
		\iff H_a=H_b\iff ba=0.
		\]
		If these conditions hold, then $(a+b)^2=0$.
		
		Suppose now that $H_a\neq H_b$. Since these are distinct
		one-dimensional subspaces of $D_1(\Bbbk)^2$, we have
		$H_a\cap H_b=0$. If $\rho(a+b)(v)=0$, then
		$\rho(a)(v)=-\rho(b)(v)\in H_a\cap H_b$, so
		$\rho(a)(v)=\rho(b)(v)=0$. Hence
		$v\in H_a\cap H_b=0$, proving that $\rho(a+b)$
		is invertible. Since $\rho$ is injective, $a+b$
		is regular and cannot be nilpotent.
	\end{proof}

	For a nonzero nilpotent element $a\in\mathcal B$, we fix
	$
	\mathcal L_a:=\{b\in\mathcal B:ab=ba=0\}.
	$

\begin{proposition}\label{prop:maximal-linear-nilpotent-subspaces}
	Let $a\in\mathcal B$ be a nonzero nilpotent element.
	Then $\mathcal L_a$ is the unique maximal
	$\Bbbk$-linear subspace consisting of nilpotent elements
	and containing $a$. Moreover, $\mathcal L_a^2=0$.
	Every maximal such subspace is of this form, and
	any two distinct ones intersect trivially.
\end{proposition}

\begin{proof}
	By definition, $\mathcal L_a$ is a linear subspace
	containing $a$. For every $b\in\mathcal L_a$, the
	defining annihilation conditions yield
	$\operatorname{Im}\rho(b)\subseteq H_a
	\subseteq\ker\rho(b)$, so $b^2=0$.
	Similarly, for any $b,c\in\mathcal L_a$, these
	inclusions imply $bc=0$. Thus $\mathcal L_a^2=0$.
	Let $L$ be any linear subspace consisting of nilpotent
	elements and containing $a$. For every $b\in L$,
	the sum $a+b$ is nilpotent. By
	Proposition~\ref{prop:sums-of-nilpotents},
	$ab=ba=0$, whence $b\in\mathcal L_a$.
	Therefore, $L\subseteq\mathcal L_a$, proving
	maximality and uniqueness.
		Every maximal such subspace contains a nonzero
	element and hence has the form $\mathcal L_a$.
	Finally, if $0\neq c\in\mathcal L_a\cap\mathcal L_b$,
	uniqueness applied to $c$ yields
	$\mathcal L_a=\mathcal L_c=\mathcal L_b$.
	Thus distinct maximal subspaces intersect trivially.
\end{proof}

\begin{corollary}\label{cor:linear-family-maximality}
Let $\mathcal S:=\bigcup_{m\geq0}S_m$ and
$\mathcal H:=\bigcup_{m\geq1}H_m$. Then
	\[
	\mathcal S=x\mathcal J=\mathcal L_x,
	\qquad
	\mathcal H=\mathcal Jzy=\mathcal L_{zy}.
	\]
	Consequently, $\mathcal S$ and $\mathcal H$ are
	maximal linear subspaces consisting of nilpotent
	elements. Moreover, for every $m\geq0$, $S_m$
	is maximal among linear subspaces of
	$\mathcal B_{\leq 2m+1}$ consisting of square-zero
	elements and the same holds for $H_m$ when $m\geq1$.
\end{corollary}

\begin{proof}
	By Lemma~\ref{lem:matrix-image-nilpotent-ideal},
	we have
	$\rho(\mathcal S)=\mathcal J e_{21}=\rho(x\mathcal J)$
	and
	$\rho(\mathcal H)=\mathcal Jw^2e_{12}
	=\rho(\mathcal Jzy)$.
	By Proposition~\ref{prop:maximal-linear-nilpotent-subspaces},
	every element of $\mathcal L_x$ and $\mathcal L_{zy}$
	is nilpotent and hence belongs to $I$.
	The description of $\rho(I)$, together with the
	annihilation conditions defining these subspaces, gives
	$\rho(\mathcal L_x)=\mathcal J e_{21}$ and
	$\rho(\mathcal L_{zy})=\mathcal Jw^2e_{12}$.
	Since $\rho$ is injective, the asserted equalities follow.
	Their maximality follows from
	Proposition~\ref{prop:maximal-linear-nilpotent-subspaces}.
	
	Finally, let $L\subseteq\mathcal B_{\leq 2m+1}$
	be a linear subspace consisting of square-zero elements
	and containing $S_m$. Since $x\in S_m$, maximality of
	$\mathcal L_x$ implies
	$L\subseteq\mathcal L_x=\mathcal S$.
	Hence
	$L\subseteq\mathcal S\cap\mathcal B_{\leq 2m+1}=S_m$.
	The same argument applies to $H_m$ for $m\geq1$,
	using $zy\in H_m$ and $\mathcal L_{zy}=\mathcal H$.
\end{proof}

		Recall that $\mathbf J=\Bbbk\langle y^2,w\rangle$ is the copy of the Jordan plane considered in \cite{AD}. Notice that the element denoted by $z$ in \cite{AD} corresponds to  $w=xy+yx$. The argument in the proof of \cite[Proposition 3.5]{AD} gives the
	following.

	\begin{lemma}\label{evensquarezero}
		The even component $\mathcal B_{\bar 0}$ contains
		no nonzero square-zero elements.
	\end{lemma}
	
	\begin{proof}
		By \cite[Proposition 3.1]{AD}, we have
		$\mathcal B_{\bar 0}=\mathbf J\oplus\mathbf Jxy$.
		Thus, every $a\in\mathcal B_{\bar 0}$ can be
		written uniquely as $a=f+gxy$, with $f,g\in\mathbf J$.
		Moreover, by \cite[Lemma 2.4]{AD}, there is an
		automorphism $\eta$ of $\mathbf J$ satisfying
		$\eta(w)=w$, $\eta(y^2)=y^2-w$, and
		$xyh=\eta(h)xy$ for all $h\in\mathbf J$.
		Since $(xy)^2=wxy$, the product formula in
		\cite[Lemma 3.4]{AD} yields
		\[
		a^2=f^2+\bigl(fg+g\eta(f+gw)\bigr)xy.
		\]
		Suppose that $a^2=0$. By the direct sum decomposition,
		$f^2=0$, and hence $f=0$, since $\mathbf J$ is a domain.
		Consequently, $g\eta(g)w=0$. As $\eta$ is an
		automorphism and $w\neq0$, the fact that $\mathbf J$
		is a domain implies $g=0$. Therefore, $a=0$.
	\end{proof}

	\begin{lemma}\label{squarezeroodddegree}
		Let $0\neq a\in\mathcal B$ satisfy $a^2=0$.
		Then $\deg a$ is odd.
	\end{lemma}
	
	\begin{proof}
		Write $a=\sum_{i=0}^d a_i$, where
		$a_i\in\mathcal B_i$ and $a_d\neq0$.
		Since $a^2=0$, its homogeneous component of
		degree $2d$ gives $a_d^2=0$.
		If $d$ were even, then
		$a_d\in\mathcal B_{\bar 0}$, and
		Lemma~\ref{evensquarezero} would imply $a_d=0$,
		a contradiction. Hence $d=\deg a$ is odd.
	\end{proof}

By Lemma~\ref{squarezeroodddegree}, every nonzero
square-zero element of $\mathcal B$ has odd degree.
Consequently, for every $m\geq0$,
$
\{a\in\mathcal B_{\leq 2m+2}:a^2=0\}
=
\{a\in\mathcal B_{\leq 2m+1}:a^2=0\}.
$
Thus, increasing the truncation degree from $2m+1$
to $2m+2$ introduces no new square-zero elements. For the geometric discussion that follows, we assume that
$\Bbbk$ is algebraically closed. We therefore define the affine square-zero locus
\begin{equation}\label{affine-zero-locus}
	\mathcal Z_m
	:=\{a\in\mathcal B_{\leq 2m+1}:a^2=0\}
\end{equation}
and its projectivization
\begin{equation}\label{zero-locus}
	\overline{\mathcal Z}_m
	:=\{[a]\in\mathbb P(\mathcal B_{\leq 2m+1}):a^2=0\}.
\end{equation}
		
\begin{proposition}\label{prop:closed-square-zero-locus}
		The subset $\overline{\mathcal Z}_m$ is a Zariski-closed algebraic subset of
		$\mathbb P(\mathcal B_{\leq2m+1})$ defined by homogeneous quadratic
		equations.
	\end{proposition}
	
	\begin{proof}
		Let $\{b_1,\dots,b_N\}$ be a PBW basis of $\mathcal B_{\leq2m+1}$. Since $\mathcal B$ is
		graded, $\mathcal B_{\leq2m+1}\mathcal B_{\leq2m+1}
		\subseteq\mathcal B_{\leq4m+2}$. Let $\{u_1,\dots,u_M\}$
		be a PBW basis of $\mathcal B_{\leq4m+2}$, and write
		$a=\sum_{i=1}^N\lambda_i b_i$. By bilinearity of multiplication, $a^2=\sum_{i,j=1}^N
		\lambda_i\lambda_j b_i b_j$. Expressing each product $b_i b_j$ in the PBW basis $\{u_1,\dots,u_M\}$, we obtain homogeneous quadratic forms $F_1,\dots,F_M \in
		\Bbbk[\lambda_1,\dots,\lambda_N]$ such that
		\[
		a^2=\sum_{k=1}^M
		F_k(\lambda_1,\dots,\lambda_N)u_k.
		\]
		Since the $u_k$ are linearly independent, $a^2=0$ if and only if $F_k(\lambda_1,\dots,\lambda_N)=0$, for every $k=1,\dots,M$. Therefore, $\overline{\mathcal Z}_m
		=
		V(F_1,\dots,F_M)
		\subseteq
		\mathbb P(\mathcal B_{\leq2m+1})$. Furthermore, since the $F_k$ are homogeneous, their common zero set is Zariski closed in projective space.
	\end{proof}

Since every element of $S_m$ and $H_m$ has square zero,
their projectivizations are contained in
$\overline{\mathcal Z}_m$. Moreover,
\[
\mathbb P(S_m)\simeq
\mathbb P^{\binom{m+2}{2}-1},
\qquad
\mathbb P(H_m)\simeq
\mathbb P^{\binom{m+1}{2}-1},
\]
where the second identification holds for $m\geq1$.
Since $S_m\cap H_m=\{0\}$, these projective linear
subspaces are disjoint. Consequently, for $m\geq1$,
$
\mathbb P(S_m)\sqcup\mathbb P(H_m)
\subseteq\overline{\mathcal Z}_m.
$
In particular,
$\dim\overline{\mathcal Z}_m
\geq\binom{m+2}{2}-1$
for every $m\geq0$.
This motivates the following conjecture.

	\begin{Conjecture}\label{dimensionZm}
		For every $m\geq0$,
		$
		\dim\overline{\mathcal Z}_m
		=
		\binom{m+2}{2}-1.
		$
	\end{Conjecture}
	
	Since $\mathcal Z_m$ is the affine locus over the nonempty projective variety $\overline{\mathcal Z}_m$, one has
	$\dim\mathcal Z_m=\dim\overline{\mathcal Z}_m+1$. Notice that if Conjecture \ref{dimensionZm} holds, then
	$\mathbb P(S_m)$ is an irreducible component of maximal dimension of $\overline{\mathcal Z}_m$.

\begin{remark}
	The linear families $S_m$ and $H_m$ are closely related
	to the classification of square-zero elements in the
	odd component of $\mathcal B$ obtained by Andruskiewitsch
	and Dumas \cite[Proposition 4.2]{AD}.
	Recall that $\mathbf J=\Bbbk\langle y^2,w\rangle$
	is a copy of the Jordan plane contained in
	$\mathcal B_{\bar 0}$. Every element
	$f\in\mathcal B_{\bar 1}$ can be written uniquely
	as $f=f^+x+f^-y$, with $f^+,f^-\in\mathbf J$.
	By \cite[Proposition 4.2]{AD}, we have $f^2=0$
	if and only if either $f^-=0$ or
	$f^-y^2=-f^+w$.
	
	In the first case, $f\in\mathbf Jx$.
	In the second, reduction modulo $(w)$ gives
	$f^-=0$ in $\mathbf J/(w)\simeq\Bbbk[y^2]$.
	Thus $f^-=hw$ for some $h\in\mathbf J$.
	Since $wy^2=(y^2-w)w$, the relation
	$f^-y^2=-f^+w$ implies
	$f^+=-h(y^2-w)$.
	Consequently,
	$f=h\bigl(-(y^2-w)x+wy\bigr)=hzy$.
		Hence the classification in \cite{AD} takes the form
	$
	\{f\in\mathcal B_{\bar 1}:f^2=0\}
	=\mathbf Jx\cup\mathbf Jzy.
	$
	The PBW expansions give
	$\mathbf Jx\cap\mathcal B_{\leq 2m+1}=S_m$
	and
	$\mathbf Jzy\cap\mathcal B_{\leq 2m+1}=H_m$
	for $m\geq1$. Therefore,
	\[
	\overline{\mathcal Z}_m
	\cap\mathbb P(\mathcal B_{\bar 1}
	\cap\mathcal B_{\leq 2m+1})
	=\mathbb P(S_m)\sqcup\mathbb P(H_m).
	\]
	Thus, within the odd component, the projective
	square-zero locus consists of exactly two disjoint
	projective linear subspaces.
\end{remark}

	\medbreak
	
	Combining the preceding description of the odd component
	with Proposition~\ref{prop:nilpotent-rational-normal-form},
	we obtain the following three disjoint families of nonzero nilpotent elements
	of $\mathcal B$.
	The first two are the linear families
	$\mathcal S=x\mathcal J$ and $\mathcal H=\mathcal Jzy$,
	consisting of elements of the form $xc$ and $gzy$,
	respectively, with $0\neq c,g\in\mathcal J$.
	The third is the family
	\[
	\mathcal M
	:=\{a(r,c):r,c\neq0
	\text{ satisfy the conditions of
		Proposition~\ref{prop:nilpotent-rational-normal-form}}\},
	\]
	which we call the family of \emph{mixed nilpotent elements}.
		The families $\mathcal S$ and $\mathcal H$ consist of
	odd elements, since
	$\mathcal J\subseteq\mathcal B_{\bar 0}$.
	For $a=a(r,c)\in\mathcal M$, the parity decomposition
	$a=a_{\bar 0}+a_{\bar 1}$ is given by
	\[
	\begin{aligned}
		a_{\bar 0}
		&=(rcw^{-1})z-(crw^{-1})xy,\qquad a_{\bar 1}
		&=xc-(rcrw^{-2})zy.
	\end{aligned}
	\]
	The admissibility conditions ensure that all coefficients
	belong to $\mathcal J$. Moreover,
	\[
	\rho(a_{\bar 0})=rc\,e_{11}-cr\,e_{22},
	\qquad
	\rho(a_{\bar 1})=-rcr\,e_{12}+c\,e_{21}.
	\]
	Since $D_1(\Bbbk)$ is a division ring and $r,c\neq0$,
	both parity components are nonzero, justifying the
	terminology \emph{mixed}. In the next subsection, we construct a quadratic
	Veronese family of square-zero elements and determine
	its intersections with the linear families $S_m$ and
	$H_m$, as well as its mixed elements.

\subsection{Veronese varieties in the square-zero loci}
	
	Let $V$ be a vector space of dimension $m+1$ over $\Bbbk$, and let $d\geq 1$. Recall that $\operatorname{Sym}^d(V)$ is the quotient of $V^{\otimes d}$
	by the subspace generated by
	\[
	v_1\otimes\cdots\otimes v_d-
	v_{\sigma(1)}\otimes\cdots\otimes v_{\sigma(d)},
	\quad \sigma\in S_d.
	\]
	The class of $v_1\otimes\cdots\otimes v_d$ is denoted by
	$v_1\odot\cdots\odot v_d$. In particular, $v^{\odot d}$ denotes the
	class of $v^{\otimes d}$ in $\operatorname{Sym}^d(V)$. The {\it Veronese map} of degree $d$ (associated to $V$) is the morphism
	\[
	\overline \nu_d:\mathbb P(V)\longrightarrow\mathbb P(\operatorname{Sym}^d(V)),
	\quad
	\overline \nu_d([v])=[v^{\odot d}], \quad 0\neq v\in V.
	\]
	The map $\overline \nu_d$ is well defined because replacing
	$v$ by $\lambda v$, with $\lambda\in\Bbbk^\times$, replaces
	$v^{\odot d}$ by $\lambda^d v^{\odot d}$ and hence does not change its
	projective class. Fix a basis $v_0,\ldots,v_m$ of $V$ and write
	$v=t_0v_0+\cdots+t_mv_m$. After choosing suitably normalized homogeneous coordinates on
	$\mathbb P(\operatorname{Sym}^d(V))$, the map $\overline \nu_d$ is given by
	\[
	\overline \nu_d([t_0:\cdots:t_m])
	=
	\bigl[t_0^{\alpha_0}\cdots t_m^{\alpha_m}\bigr]_{
		\substack{\alpha_0+\cdots+\alpha_m=d\\ }} \, \, \, .
	\]
	Thus, the coordinates of $\overline \nu_d([v])$ are all the monomials of degree
	$d$ in $t_0,\ldots,t_m$. Since
	$\dim\operatorname{Sym}^d(V)=\binom{m+d}{d}$, after choosing
	coordinates, $\overline \nu_d$ has the form
	$\overline \nu_d:\mathbb P^m\longrightarrow
	\mathbb P^{\binom{m+d}{d}-1}$. It is a closed embedding, and its image
	is called the {\it  Veronese variety} of degree $d$; see more details in Lecture 2 of \cite{HarrisAG}. In the case $d=2$ and $m=1$, one
	recovers the quadratic Veronese map, that is,
	\[
	\overline \nu_2:\mathbb P^1\longrightarrow\mathbb P^2,
	\quad
	[\lambda:\mu]\longmapsto
	[\lambda^2:\lambda\mu:\mu^2].
	\]
	
	In order to relate the quadratic Veronese map with square-zero elements of the super Jordan plane, consider the vector space $V_m:=\{p(y)\in\Bbbk[y]:\deg p\leq m\}$ and the map $\nu_2:V_m\longrightarrow\operatorname{Sym}^2(V_m)$  given by $\nu_2(p)=p\odot p$, for all $p\in V_m$. Also, we define the map
	\[
	Q_m:V_m\longrightarrow\mathcal B_{\leq 2m+1},
	\quad
	Q_m(p):=p(-y)xp(y).
	\]
	\medbreak
	\begin{proposition}\label{prop-veronese}
		The map $Q_m$ satisfies the following properties:
		\begin{enumerate}[\rm (i)]
			\setlength{\itemsep}{0.3em}
			\item $Q_m(p)$ is a square-zero element of $\mathcal B$ of degree at most
			$2m+1$, for every $p\in V_m$.
			\item There exists a linear map
			$
			L_m:\operatorname{Sym}^2(V_m)\longrightarrow
			\mathcal B_{\leq2m+1}
			$
			such that $Q_m=L_m\circ\nu_2$.
		\end{enumerate}
	\end{proposition}
	
	\begin{proof}
		Let $p\in V_m$. Since $p(y)p(-y)$ is invariant under the
		substitution $y\mapsto -y$, it is an even polynomial. Then, there exists $r(t)\in\Bbbk[t]$ such that
		$p(y)p(-y)=r(y^2)$. Therefore, $Q_m(p)^2
		=p(-y)x\,r(y^2)\,xp(y)$. Writing $r(t)=\sum_{k\geq0}c_kt^k$ and using the identity $xy^{2k}x=0$, for every $k\geq0$, we obtain
		$
		xr(y^2)x=\sum_{k\geq0}c_kxy^{2k}x=0.
		$		
		Thus, $Q_m(p)^2=0$. Moreover, since $\deg p\leq m$, $\deg Q_m(p)\leq 2m+1$ whenever $Q_m(p)\neq0$. This proves (i).
		
	Define the linear map $L_m: \operatorname{Sym}^2(V_m)
		\longrightarrow \mathcal B_{\leq2m+1}$ on the basis
		$\{y^r\odot y^s:0\leq r\leq s\leq m\}$ by
	\begin{equation}\label{map-L}
		L_m(y^r\odot y^s)
		=
		\begin{cases}
			(-1)^ry^rxy^r, & \text{if }r=s,\\[2mm]
			\dfrac{1}{2}
			\left(
			(-1)^ry^rxy^s+(-1)^sy^sxy^r
			\right),
			& \text{if }r<s.
		\end{cases}
\end{equation}
	Notice that $\deg(y^rxy^s)=r+s+1\leq2m+1$, so $L_m$ indeed takes values in $\mathcal B_{\leq2m+1}$. Let
		$p(y)=\sum_{r=0}^m a_ry^r$. Then,
		$
		p\odot p
		=
		\sum_{r=0}^m a_r^2(y^r\odot y^r)
		+
		2\sum_{0\leq r<s\leq m}
		a_ra_s(y^r\odot y^s)
		$ and consequently,
		\[
		\begin{aligned}
			L_m(p\odot p)
			&=
			\sum_{r=0}^m
			a_r^2(-1)^ry^rxy^r+
			\sum_{0\leq r<s\leq m}
			a_ra_s
			\left(
			(-1)^ry^rxy^s
			+
			(-1)^sy^sxy^r
			\right)\\
			&=
			p(-y)xp(y)
			=
			Q_m(p).
		\end{aligned}
		\]
	Therefore, $Q_m=L_m\circ\nu_2$, which proves (ii).
	\end{proof}
	
	In order to prove that $L_m$ is injective, we distinguish the images of the basis elements $y^r\odot y^s$ by means of suitable PBW monomials. Taking $n=1$ in \eqref{eqpropy1yx} and \eqref{eqpropy2yx},
	using $z=yx$, and multiplying on the right by $y^\ell$, we obtain
	the following formulas for $k,\ell\geq0$; the first formula
	for $k=0$ is immediate:
	\begin{align}
		y^{2k}xy^\ell
		&=
		\sum_{i=0}^k
		\binom{k}{i}i!\,
		xz^iy^{2(k-i)+\ell},
		\label{eq:even-yxy}\\
		y^{2k+1}xy^\ell
		&=
		\sum_{i=0}^k
		\binom{k}{i}i!\,
		z^{i+1}y^{2(k-i)+\ell}.
		\label{eq:odd-yxy}
	\end{align}
	For $0\leq r\leq s\leq m$, denote $X_{rs}:=L_m(y^r\odot y^s)$. By the definition of $L_m$, each $X_{rs}$ is homogeneous of degree $r+s+1$.
	
	\begin{lemma}\label{lem:distinguished-PBW-monomial}
		Let $m$ be a nonnegative integer and $0\leq r\leq s\leq m$. 
		\begin{enumerate}[\rm (i)]
			\setlength{\itemsep}{0.3em}
			\item If $s=2k$, then the PBW monomial $xz^ky^r$ occurs in $X_{rs}$ with a nonzero coefficient.
			\item If $s=2k+1$, then the PBW monomial $z^{k+1}y^r$ occurs in $X_{rs}$ with a nonzero coefficient.
		\end{enumerate}
		
		Moreover, let $r',s'$ satisfy
		$r'+s'=r+s$ and $r'\leq s'<s$.
		
		\begin{enumerate}[\rm (i)]
			\setcounter{enumi}{2}
			\setlength{\itemsep}{0.3em}
			\item If $s=2k$, then the PBW monomial $xz^ky^r$
			does not occur in $X_{r's'}$.
			\item If $s=2k+1$, then the PBW monomial $z^{k+1}y^r$
			does not occur in $X_{r's'}$.
		\end{enumerate}
	\end{lemma}
	
	\begin{proof}
		Suppose that $s=2k$. By \eqref{eq:even-yxy},
		$y^sxy^r$ contains the PBW monomial $xz^ky^r$ with coefficient $k!$. If $r<s$, then the other term occurring in $X_{rs}$ is $y^rxy^s$. Since $r<s=2k$, formulas \eqref{eq:even-yxy} and \eqref{eq:odd-yxy} show that $y^rxy^s$ cannot contain the monomial $xz^ky^r$. Consequently, the coefficient of $xz^ky^r$ in $X_{rs}$
		is $k!/2$ when $r<s$, and $k!$ when $r=s$. In either case, it is nonzero.
		Suppose that $s=2k+1$. By \eqref{eq:odd-yxy},
		$y^sxy^r$ contains the PBW monomial $z^{k+1}y^r$ with coefficient $k!$. If $r<s$, then $y^rxy^s$ cannot contain this monomial, since $r<s=2k+1$ and formulas \eqref{eq:even-yxy} and \eqref{eq:odd-yxy} show that a left exponent strictly smaller than $2k+1$ cannot produce the power $z^{k+1}$ without a left factor $x$. Thus, the coefficient of $z^{k+1}y^r$ in $X_{rs}$ is $-k!/2$ when
		$r<s$, and $-k!$ when $r=s$. Again, it is nonzero.
			Every PBW monomial occurring in $X_{r's'}$ comes from either $y^{r'}xy^{s'}$ or $y^{s'}xy^{r'}$. Since
		$r'\leq s'<s$, the exponent of $y$ to the left of $x$ in both expressions is strictly smaller than $s$.
		
		If $s=2k$, an even left exponent $2j<2k$ satisfies $j<k$, and \eqref{eq:even-yxy} produces only PBW monomials of the form $xz^iy^c$ with $i\leq j<k$. An odd left exponent produces, by \eqref{eq:odd-yxy}, only monomials of the form $z^{i+1}y^c$. Then, $xz^ky^r$ cannot occur in $X_{r's'}$.
		Similarly, if $s=2k+1$, an odd left exponent $2j+1<2k+1$ satisfies $j<k$, and \eqref{eq:odd-yxy} produces only PBW monomials of the form $z^{i+1}y^c$ with $i\leq j<k$. An even left exponent produces only monomials of the form $xz^iy^c$. Then, $z^{k+1}y^r$ cannot occur in
		$X_{r's'}$.
	\end{proof}
	
	\begin{proposition}\label{prop:Lm-injective}
		For every $m\geq0$, the linear map
		$
		L_m:
		\operatorname{Sym}^2(V_m)
		\longrightarrow
		\mathcal B_{\leq2m+1}
		$ defined in \eqref{map-L}
		is injective.
	\end{proposition}
	
	\begin{proof}
		Suppose that $\sum_{0\leq r\leq s\leq m}
		c_{rs}X_{rs}=0$. Since $X_{rs}\in\mathcal B_{r+s+1}$ and the grading of $\mathcal B$ is a direct sum, each homogeneous component of this relation must vanish separately. Thus, for every
		$0\leq n\leq2m$, $Y_n:=\sum_{\substack{0\le r\le s\le m\\r+s=n}}c_{rs}X_{rs}=0$.
				Suppose, by contradiction, that for some $n$ there exists a pair $(r,s)$ with $r+s=n$ and $c_{rs}\neq0$. Choose such a pair with $s$ maximal.
		Suppose that $s=2k$. By Lemma \ref{lem:distinguished-PBW-monomial}, the PBW monomial $xz^ky^r$ occurs in $X_{rs}$ with a nonzero coefficient, $\alpha_{rs}$, and does not occur in any $X_{r's'}$ with $r'+s'=n$, $s'<s$. On the other hand, by the maximality of $s$, $c_{r's'}=0$ whenever $r'+s'=n$, $s'>s$. Moreover, if $s'=s$, then $r'=r$, since $r'+s'=r+s=n$. Consequently, the coefficient of $xz^ky^r$ in $Y_n$ is precisely $c_{rs}\alpha_{rs}$. Since $Y_n=0$ and $\alpha_{rs}\neq0$, we obtain $c_{rs}=0$, contradicting the choice of $(r,s)$.
		The case $s=2k+1$ is analogous. In this case,
		Lemma \ref{lem:distinguished-PBW-monomial} shows that the PBW monomial $z^{k+1}y^r$ occurs in $X_{rs}$ with a nonzero coefficient and does not occur in
		any $X_{r's'}$ with $r'+s'=n$, $s'<s$. The same maximality argument again gives $c_{rs}=0$, a contradiction. Therefore, $c_{rs}=0$, for all $0\leq r\leq s\leq m$ and, then, the map $L_m$ is injective.
	\end{proof}
	
	We show that the quadratic family determined by $Q_m$ has a natural projective interpretation: after projectivization, it defines a closed embedding whose image is a quadratic Veronese variety contained in the square-zero locus $\overline{\mathcal Z}_m$.
	
	\begin{theorem}\label{thm:Veronese-square-zero-locus}
		For every $m\geq0$, the map $\overline Q_m:\mathbb P(V_m)\longrightarrow\mathbb P(\mathcal B_{\leq2m+1})$ given by
		\[\overline Q_m([p])=[Q_m(p)], \quad 0\neq p\in V_m,\]
		is a closed embedding. Moreover, $\overline Q_m=\mathbb P(L_m)\circ\overline \nu_2$
		and $\overline Q_m\big(\mathbb P(V_m)\big)$ is a copy of the quadratic Veronese variety $\overline \nu_2(\mathbb P(V_m))$ contained in $\overline{\mathcal Z}_m$.
	\end{theorem}
	
	\begin{proof}
		Let $0\neq p\in V_m$. Since $p\odot p\neq0$ and $L_m$ is injective, 
		$Q_m(p)=L_m(p\odot p)\neq0$.
		Moreover, $Q_m(\alpha p)=\alpha^2Q_m(p)$, for every $\alpha\in\Bbbk^\times$. Then, $[Q_m(\alpha p)]=[Q_m(p)]$
		 and therefore the map $\overline Q_m$ is well defined. By Proposition \ref{prop:Lm-injective}, the linear map $L_m$ is injective and then its projectivization
		$
		\mathbb P(L_m):
		\mathbb P(\operatorname{Sym}^2(V_m))
		\to
		\mathbb P(\mathcal B_{\leq2m+1})
		$
		is a linear closed embedding. Since $Q_m=L_m\circ\nu_2$, we obtain $\overline Q_m
		=\mathbb P(L_m)\circ\overline \nu_2$. The quadratic Veronese map $\overline \nu_2$ is a closed embedding. Therefore,
		$\overline Q_m$ is also a closed embedding. Moreover,
		$\mathbb P(L_m)$ identifies $\overline \nu_2(\mathbb P(V_m))$ with
		$\overline Q_m(\mathbb P(V_m))$, so the latter is projectively equivalent to the quadratic Veronese variety
		$\overline \nu_2(\mathbb P(V_m))$. Finally, Proposition \ref{prop-veronese} gives $Q_m(p)^2=0$, for every $p\in V_m$. Then, $
		\overline Q_m(\mathbb P(V_m))\subseteq\overline{\mathcal Z}_m$.
	\end{proof}
	
We now relate the Veronese construction to the
	description of the nilpotent elements of $\B$. Write
	\[
	Q(p)=p(-y)xp(y), \quad p\in\Bbbk[y],
	\quad
	\mathcal V=\{Q(p):p\in\Bbbk[y]\}.
	\]
	The maps $Q_m$ of the Veronese construction are the restrictions
	of $Q$ to the polynomials of degree at most $m$.
	By Proposition \ref{prop-veronese}, every element of $\mathcal V$
	has square zero.
	
	\begin{lemma}\label{lem:degree-veronese-elements}
		For every nonzero $p\in\Bbbk[y]$,
		$
		\deg Q(p)=2\deg p+1.
		$
	\end{lemma}
	
	\begin{proof}
		If $p(y)=\alpha_dy^d+\cdots+\alpha_0$ with $\alpha_d\neq0$,
		the component of degree $2d+1$ of $Q(p)$ is
		$(-1)^d\alpha_d^2y^dxy^d$.
		The PBW expansion of $y^dxy^d$ contains $xy^{2d}$ with coefficient $1$
		when $d$ is even and $zy^{2d-1}$ with coefficient $1$
		when $d$ is odd, and hence the result follows.
	\end{proof}
	\begin{proposition}\label{prop:veronese-normal-form}
		Let $0\neq p\in\Bbbk[y]$, and write uniquely
		$p(y)=f(y^2)+yg(y^2)$ with $f,g\in\Bbbk[t]$.
		Then:
		\begin{enumerate}[\rm (i)]
			\item $Q(p)\in\mathcal S$ if and only if $g=0$.
			
			\item $Q(p)\in\mathcal H$ if and only if $f=0$.
			
			\item If $f,g\neq0$, both parity components of $Q(p)$
			are nonzero.
			
			\item If $f\neq0$, then $Q(p)=a(r,c)$, where $r=-g(u)f(u)^{-1}w$, $c=f(u+w)f(u)$
			and $a(r,c)$ is defined in
			Proposition~\ref{prop:nilpotent-rational-normal-form}(ii).
		\end{enumerate}
	\end{proposition}
	
	\begin{proof}
Let $\rho$ be the map defined in \eqref{injec-rho}. A direct
computation gives
\[
\rho(Q(p))=
\begin{pmatrix}-g(u)w\\ f(u+w)\end{pmatrix}
\begin{pmatrix}f(u)&g(u)w\end{pmatrix}.
\]
		Thus the lower-left entry is $f(u+w)f(u)$ and the
		upper-right entry is $-g(u)wg(u)w$.
		Since $\mathcal J$ is a domain, these vanish precisely when
		$f=0$ and $g=0$, respectively.
		If both polynomials are nonzero, the upper-left entry
		$-g(u)wf(u)$ is also nonzero.
		The assertions about the families and parity follow.
		If $f\neq0$, the lower-left entry is nonzero, and the
		uniqueness formulas in Proposition
		\ref{prop:nilpotent-rational-normal-form} give
		\[
		\begin{aligned}
			r
			&=-g(u)wf(u)\bigl(f(u+w)f(u)\bigr)^{-1}=-g(u)wf(u+w)^{-1}=-g(u)f(u)^{-1}w.
		\end{aligned}\]
		The last equality follows from $f(u)w=wf(u+w)$.
	\end{proof}

	\begin{example}\label{ex:mixed-veronese-element}
		Taking $p(y)=1-y$, we obtain
		\[
		a:=Q(1-y)=(1+y)x(1-y)
		=x+z-xy-zy=a(w,1).
		\]
		Its parity components are $a_{\bar0}=z-xy$
		and $a_{\bar1}=x-zy$, both of which are
		nonzero. Thus $a$ is a mixed nilpotent element
		belonging to neither $\mathcal S$ nor $\mathcal H$.
	\end{example}
	
\begin{proposition}\label{propZ1}
	The projective square-zero locus $\overline{\mathcal Z}_1$
	has exactly two irreducible components:
	$
	\overline{\mathcal Z}_1=\mathbb P(S_1)\cup\mathcal C,
	$
	where $\mathbb P(S_1)\simeq\mathbb P^2$ and $\mathcal C$
	is the nonsingular conic
	\[
	\mathcal C=
	\left\{
	[Ax+B(xy-z)+Fzy]:AF+B^2=0
	\right\}.
	\]
	The conic $\mathcal C$ coincides with
	$\overline Q_1(\mathbb P(V_1))$ and is projectively
	equivalent to $\overline\nu_2(\mathbb P^1)$.
	The two components intersect precisely at $[x]$.
	In particular, $\dim\overline{\mathcal Z}_1=2$.
\end{proposition}

\begin{proof}
	Recall that $I=(x)=\ker\bigl(\mathcal B\longrightarrow
	\mathcal B/(x)\simeq\Bbbk[y]\bigr)$.
	By Lemma~\ref{lemanilpotentkernel}(i),
	$\mathcal Z_1\subseteq I\cap\mathcal B_{\leq3}$.
	Using the PBW basis, we obtain
	\[
	I\cap\mathcal B_{\leq3}
	=
	\Bbbk x\oplus\Bbbk xy\oplus\Bbbk z
	\oplus\Bbbk xz\oplus\Bbbk xy^2\oplus\Bbbk zy.
	\]
	Thus, every $a\in\mathcal Z_1$ has the form
	$a=Ax+Bxy+Cz+Dxz+Exy^2+Fzy$ for unique
	$A,B,C,D,E,F\in\Bbbk$.
	Write $a=a_1+a_2+a_3$, where $a_1=Ax$,
	$a_2=Bxy+Cz$ and $a_3=Dxz+Exy^2+Fzy$.
	Since $a_1^2=0$, the equation $a^2=0$ is equivalent
	to the vanishing of its homogeneous components of
	degrees $3$, $4$, $5$ and $6$.
	Using the defining relations of $\mathcal B$, we obtain
	\[
	\begin{aligned}
		a_1a_2+a_2a_1
		&=A(B+C)xz,\\
		a_2^2+a_1a_3+a_3a_1
		&=(AF+B^2)xzy+(AF+C^2)z^2,\\
		a_2a_3+a_3a_2
		&=\bigl(D(B+C)+CE\bigr)xz^2\\
		&\quad+E(B+C)xzy^2+F(B+C)z^2y,\\
		a_3^2
		&=F(D+E)xz^2y+EFxzy^3\\
		&\quad+DFz^3+EFz^2y^2.
	\end{aligned}
	\]
	By the linear independence of the PBW monomials,
	$a^2=0$ if and only if
\[
\begin{aligned}
	A(B+C)&=0, & AF+B^2&=0,\\
	F(B+C)&=0, & AF+C^2&=0,\\
	E(B+C)&=0, & D(B+C)+CE&=0,\\
	DF&=0, & EF&=0.
\end{aligned}
\]
	We claim that $B+C=0$. Otherwise, the first and
	second lines imply $A=F=0$ and hence $B=C=0$,
	a contradiction. Substituting $C=-B$ into the
	preceding equations, we conclude that $a^2=0$
	if and only if
	\[
	C=-B,\qquad AF+B^2=0,\qquad BE=DF=EF=0.
	\]
	If $F=0$, then $B=C=0$, with no further restrictions
	on $A,D,E$. Thus, the elements of $\mathcal Z_1$
	with $F=0$ are precisely the elements of
	$S_1=\Bbbk x\oplus\Bbbk xz\oplus\Bbbk xy^2$.
	If $F\neq0$, then $D=E=0$, so
	$a=Ax+B(xy-z)+Fzy$ with $AF+B^2=0$.
	Conversely, every element of this form has square
	zero. Allowing $F=0$ in the latter description
	adds only scalar multiples of $x$. Therefore,
	after projectivization,
	$\overline{\mathcal Z}_1=\mathbb P(S_1)\cup\mathcal C$,
	where
	\[
	\mathcal C=
	\left\{
	[Ax+B(xy-z)+Fzy]:AF+B^2=0
	\right\}.
	\]
	
	The linear independence of $x$, $xy-z$ and $zy$
	identifies $\mathcal C$ with the conic $AF+B^2=0$
	in $\mathbb P^2$, with homogeneous coordinates
	$[A:B:F]$. This conic is nonsingular, since the
	partial derivatives $F$, $2B$ and $A$ do not vanish
	simultaneously at any projective point.
	Moreover, for $p(y)=t_0+t_1y\in V_1$, we have
	\[
	Q_1(p)=(t_0-t_1y)x(t_0+t_1y)
	=t_0^2x+t_0t_1(xy-z)-t_1^2zy.
	\]
	Thus, in the above coordinates, $\overline Q_1$
	is given by $[t_0:t_1]\mapsto[t_0^2:t_0t_1:-t_1^2]$.
	This is the quadratic Veronese embedding up to
	a projective change of coordinates and its image
	is precisely the conic $AF+B^2=0$. Hence
	$\mathcal C=\overline Q_1(\mathbb P(V_1))
	\simeq\overline\nu_2(\mathbb P^1)$, in agreement
	with Theorem~\ref{thm:Veronese-square-zero-locus}.
	Finally, by the linear independence of the PBW
	monomials, we have
	$S_1\cap\Span\{x,xy-z,zy\}=\Bbbk x$.
	Since $[x]\in\mathcal C$, it follows that
	$\mathbb P(S_1)\cap\mathcal C=\{[x]\}$.
	Both $\mathbb P(S_1)\simeq\mathbb P^2$ and
	$\mathcal C\simeq\mathbb P^1$ are irreducible,
	and neither is contained in the other.
	Therefore, they are exactly the two irreducible
	components of $\overline{\mathcal Z}_1$ and
	$\dim\overline{\mathcal Z}_1=2$.
\end{proof}

\begin{remark}\label{rem:veronese-linear-intersections}
	Let $m\geq1$ and write
	$V_m=V_m^{\mathrm{ev}}\oplus V_m^{\mathrm{odd}}$,
	where $V_m^{\mathrm{ev}}$ and $V_m^{\mathrm{odd}}$
	are the subspaces of even and odd polynomials,
	respectively.
	By Proposition~\ref{prop:veronese-normal-form},
	\[
	Q_m(V_m)\cap S_m=Q_m(V_m^{\mathrm{ev}}),
	\qquad
	Q_m(V_m)\cap H_m=Q_m(V_m^{\mathrm{odd}}).
	\]
	After projectivization, we obtain
	\[
	\begin{aligned}
		\overline Q_m(\mathbb P(V_m))\cap\mathbb P(S_m)
		&=\overline Q_m(\mathbb P(V_m^{\mathrm{ev}})),\qquad\overline Q_m(\mathbb P(V_m))\cap\mathbb P(H_m)
		&=\overline Q_m(\mathbb P(V_m^{\mathrm{odd}})).
	\end{aligned}
	\]
	By Theorem~\ref{thm:Veronese-square-zero-locus},
	these intersections are copies of the quadratic
	Veronese varieties
	$\overline\nu_2(\mathbb P^{\lfloor m/2\rfloor})$
	and $\overline\nu_2(\mathbb P^{\lceil m/2\rceil-1})$,
	respectively.
	Moreover, if $p\in V_m$ has nonzero even and odd
	parts, then both parity components of $Q_m(p)$
	are nonzero and $Q_m(p)\notin S_m\cup H_m$.
\end{remark}

	\medbreak

	Proposition~\ref{propZ1} gives a complete description of
	the projective square-zero locus in degree at most $3$.
	Notice that $\mathbb P(H_1)=\{[zy]\}$ is already contained
	in the Veronese conic, since $Q(y)=-zy$.
	It is therefore natural to ask whether the two linear
	families $\mathbb P(S_m)$ and $\mathbb P(H_m)$, together
	with the Veronese variety, exhaust $\overline{\mathcal Z}_m$.
	The following example shows that this is not the case
	for $m\geq2$.
	
	\begin{example}\label{prop:mixed-element-outside-veronese}
		Consider the element
		\[
		b:=(1+y)xz(1-y)
		=xz+z^2-xzy-z^2y=a(w,w).
		\]
		By Proposition~\ref{prop:nilpotent-rational-normal-form},
		$b^2=0$. Its degree is $5$, and its parity components
		$b_{\bar0}=z^2-xzy$ and
		$b_{\bar1}=xz-z^2y$ are both nonzero.
		Thus, $b\notin\mathcal S\cup\mathcal H$.
		
		We claim that $b$ does not belong to the Veronese
		family, even up to a nonzero scalar multiple.
		Indeed, suppose that $b=\lambda Q(p)$ for some
		$\lambda\in\Bbbk^\times$ and $p\in\Bbbk[y]$.
		By Lemma~\ref{lem:degree-veronese-elements},
		$\deg p=2$. Write $p(y)=\alpha+\beta y+\gamma y^2$,
		with $\gamma\neq0$. Comparing the homogeneous
		components of degree $5$ and using
		$y^2x=xy^2+xz$, we obtain
		\[
		-z^2y=\lambda\gamma^2(xy^4+xzy^2),
		\]
		which contradicts the linear independence of the
		PBW monomials. Hence
		$b\notin\mathcal S\cup\mathcal H\cup\mathcal V$.
	\end{example}
	
	\section*{Declaration on the use of AI tools}
	
	During the preparation of this manuscript, we used ChatGPT
	to explore approaches to the proof of Theorem~A, including
	the use of results on invertible bimodules over rings of
	differential operators. This led us to the work of
	Cannings and Holland~\cite{CH}, whose results we
	independently studied and applied. ChatGPT was also used
	to improve the English and exposition of the manuscript.
	All mathematical arguments and assumptions were verified
	by the authors, who take full responsibility for the
	content of this work.

\end{document}